\documentclass[12pt]{amsart}
\usepackage[english]{babel}
\usepackage{amsfonts,amssymb}

\newtheorem{thm}{Theorem}

\newtheorem{lemma}{Lemma}
\newtheorem{cor}{Corollary}
\newtheorem{prop}{Proposition}
\theoremstyle{definition}
\newtheorem{defn}{Definition}

\newtheorem{exam}{Example}

\newcommand{\set}[1]{\left\{#1\right\}}
\newcommand{\norm}[1]{\left\Vert#1\right\Vert}

\newcommand{\cX}{{\mathcal{X}}}
\newcommand{\cA}{{\mathcal{A}}}
\newcommand{\cB}{{\mathcal{B}}}
\newcommand{\cBp}{{\widehat{\cB}^+}}
\newcommand{\cBd}{{\widehat{\cB}}}
\newcommand{\cC}{{\mathcal{C}}}
\newcommand{\cD}{{\mathcal{D}}}

\newcommand{\cF}{{\mathcal{F}}}

\newcommand{\cH}{{\mathcal{H}}}
\newcommand{\cI}{{\mathcal{I}}}
\newcommand{\cJ}{{\mathcal{J}}}
\newcommand{\cK}{{\mathcal{K}}}
\newcommand{\cN}{{\mathcal{N}}}
\newcommand{\cO}{{\mathcal{O}}}

\newcommand{\cT}{{\mathcal{T}}}

\newcommand{\cU}{{\mathcal{U}}}
\newcommand{\cW}{{\mathcal{W}}}

\newcommand{\dC}{{\mathbb{C}}}
\newcommand{\dN}{{\mathbb{N}}}
\newcommand{\dR}{{\mathbb{R}}}
\newcommand{\dT}{{\mathbb{T}}}
\newcommand{\dZ}{{\mathbb{Z}}}

\newcommand{\gA}{\mathfrak{A}}
\newcommand{\gB}{\mathfrak{B}}
\newcommand{\gX}{\mathfrak{X}}
\newcommand{\gog}{\mathfrak{g}}
\newcommand{\goh}{\mathfrak{h}}

\newcommand{\Ind}{{\mathrm{Ind}}}
\newcommand{\Aut}{{\mathrm{Aut}}}
\newcommand{\Res}{{\mathrm{Res}}}

\newcommand{\Orb}{{\mathrm{Orb}}}
\newcommand{\St}{{\mathrm{St}}}
\newcommand{\Lin}{{\mathrm{Lin}}}
\newcommand{\Rep}{\mathrm{Rep}}
\newcommand{\Ran}{\mathrm{Ran}}

\newcommand{\mn}{\medskip}

\newcommand{\edex}{\circ}
\title{Unbounded Induced Representations of $*$-Algebras.}

\begin{document}
\author{Yurii Savchuk}
\address{Universit\"at Leipzig, Mathematisches Institut Johannisgasse 26, 04103 Leipzig, Germany}
\email{savchuk@math.uni-leipzig.de}
\thanks{The first author was supported by the International Max Planck Research School for Mathematics in the Sciences (Leipzig)}

\author{Konrad Schm\"udgen}
\address{Universit\"at Leipzig, Mathematisches Institut Johannisgasse 26, 04103 Leipzig, Germany}
\email{schmuedgen@math.uni-leipzig.de}

\dedicatory{Dedicated to the memory of A.U. Klimyk
(14.04.1939-22.07.2008)}


\date{\today.}


\maketitle
\begin{abstract}
Induced representations of $*$-algebras by unbounded operators in
Hilbert space are investigated. Conditional expectations of a
$*$-algebra $\cA$ onto a unital $*$-subalgebra $\cB$ are introduced
and used to define inner products on the corresponding induced
modules. The main part of the paper is concerned with group graded
$*$-algebras $\cA=\oplus_{g\in G}\cA_g$ for which the $*$-subalgebra
$\cB:=\cA_e$ is commutative. Then the canonical projection
$p:\cA\to\cB$ is a conditional expectation and there is a partial
action of the group $G$ on the set $\cBp$ of all characters of $\cB$
which are nonnegative on the cone $\sum\cA^2\cap\cB.$ The complete
Mackey theory is developed for $*$-representations of $\cA$ which
are induced from characters of $\cBp.$ Systems of imprimitivity are
defined and two versions of the Imprimitivity Theorem are proved in
this context. A concept of well-behaved $*$-representations of such
$*$-algebras $\cA$ is introduced and studied. It is shown that
well-behaved representations are direct sums of cyclic well-behaved
representations and that induced representations of well-behaved
representations are again well-behaved. The theory applies to a
large variety of examples. For important examples such as the Weyl
algebra, enveloping algebras of the Lie algebras $su(2),$ $su(1,1)$,
and of the Virasoro algebra, and $*$-algebras generated by dynamical
systems our theory is carried out in great detail.
\end{abstract}

\section{Introduction}
Induced representations are a fundamental tool in representation
theory of groups and algebras. They were first defined in 1898 for
finite groups by G. Frobenius and in 1955 for algebras
by D.G. Higman. If $\cB$ is a subalgebra of an algebra $\cA$ and $V$
is a left $\cB$-module, then the left $\cA$-module $\cA\otimes_\cB V$ with
action defined by $a_0(a\otimes v):=a_0a\otimes v$ is called
\textit{induced module} of $V$.

In his seminal paper \cite{rief} M. Rieffel introduced induced
representations for $C^*$-algebras and developed a major part of
Mackey's theory in this context. Another pioneering
paper is due to J.M.G. Fell \cite{fe}. A detailed treatment of this
theory is given in the monograph \cite{fd}.
An essential step in
Rieffel's inducing process is the definition of an inner product on
the algebraic tensor product $\cA\otimes_\cB V$.
That is, if there exists a conditional expectation $p$ from a $C^*$-algebra $\cA$
onto its $C^*$-subalgebra $\cB$ and if a Hilbert space $(V,\langle
\cdot, \cdot \rangle)$ is a Hermitian $\cB$-module (that is, $\langle
bx,y\rangle=\langle x,b^*y\rangle$ for $x,y\in V$ and $b\in\cB$),
then there exists a pre-inner product $\langle\cdot,\cdot\rangle_0$
on $\cA\otimes_\cB V$ such that
\begin{align}\label{innpro}
\langle a_1\otimes v_1,a_2\otimes v_2\rangle_0:=\langle
p(a_2^*a_1)v_1,v_2\rangle
\end{align} and the quotient space of
$\cA\otimes_\cB V$ by the null space of the form
$\langle\cdot,\cdot\rangle_0$ is a Hermitian $\cA$-module.

The aim of the present paper is to develop the basics of a theory of \textit{unbounded} induced $*$-representations for complex unital $*$-algebras. In contrast to the case of $C^*$-algebras there are various notions of positivity for general $*$-algebras that lead to different definitions of conditional expectations. The subtleties of positivity play a central role for our theory. We shall define (see Definition \ref{defn_cond_exp} below) a \textit{conditional expectation} from a unital $*$-algebra $\cA$ to a unital $*$-subalgebra $\cB$ to be a $\cB$-linear projection $p$ of $\cA$ onto $\cB$ which preserves involution and units and satisfies the following positivity condition:
$$p(\sum \cA^2)\subseteq \cB\cap\sum \cA^2.$$
Then a cyclic Hermitian $\cB$-module $V$ is "inducible" to $\cA$ via $p$ if and only if every element of $\cB\cap\sum\cA^2$ is represented by a positive symmetric operator on $V.$

Many bounded or unbounded $*$-representations of $*$-algebras $\cA$ are induced from appropriate $*$-subalgebras $\cB$ in our setting.
In Sections \ref{sect_env_alg}--\ref{sect_further_exam} we shall see that for a number of important $*$-algebras the "nice" irreducible $*$-representations are precisely those representations which are induced from characters which are non-negative on the cone $\cB\cap\sum\cA^2$. Among them are the $*$-algebras of the quantum group $SU_q(2)$ and of the Podles' spheres which have only bounded representations. This underlines the crucial role of positivity and it shows that our theory might be useful for general countably generated group graded $*$-algebras. It should be emphasized for all our examples neither the theory in \cite{fd} nor induction of $C^*$-algebras applies.

Let us briefly explain the basic idea for the Weyl algebra. We do not carry out all details of proofs, because this is just the special case $f(t)=1+t$ of the $*$-algebra treated in Section \ref{sect_dyn_sys}.

\begin{exam}\label{exam_weyl_alg} Let $\cA$ be the Weyl algebra $\dC\langle a,a^*|aa^*-a^*a=1\rangle$ and let $\cB$ be the unital $*$-subalgebra $\dC[N]$ of polynomials in $N:=a^*a.$ Each element $x\in\cA$ can be written as
$$x=\sum_{r=0}^k a^rf_r(N)+\sum_{s=1}^l a^{*s}f_{-s}(N)$$
with polynomials $f_j\in\dC[N]$ uniquely determined by $x$. Defining $p(x)=f_0(N)$, we obtain a conditional expectation $p$ from $\cA$ to $\cB$. It can be proved (see \cite{fs} or formula (\ref{eq_cone_in_weyl_alg}) below) that an element $f(N)\in\dC[N]$ belongs to $\cB\cap\sum\cA^2$ if and only if there are polynomials $g_0,\dots,g_k\in\dC[N]$ such that
\begin{align}\label{bsuma2}
f(N)=g_0(N)^*g_0(N)+Ng_1(N)^*g_1(N)+\cdots+N(N-1)\cdots(N-k+1)g_k(N)^*g_k(N).
\end{align}
For $\lambda\in\dR$, let $V_\lambda=\dC$ be the one-dimensional $\cB$-module given by $N=\lambda$. It is not difficult to show that $f(N)=f(\lambda)\geq 0$ for each polynomial $f(N)$ of the form (\ref{bsuma2}) if and only if $\lambda \in \dN_0$.

Now suppose that $\lambda\in\dN_0$. Let $\cH_\lambda$ denote the Hilbert space obtained from the pre-inner product (\ref{innpro}) on $\cA\otimes_\cB V_\lambda$. Clearly, the vectors $a^r\otimes 1,a^{*(r+1)}\otimes 1$, where $r\in\dN_0,$ form a base of the vector space $\cA\otimes_\cB V_\lambda.$ From the relation $aa^*-a^*a=1$ it follows that
\begin{gather}\label{aN}
a^ra^{*r} =(N+1)\dots(N+r),\ a^{*r}a^r =N(N-1)\dots(N-r+1)
\end{gather}
for $r \in \dN_0$. If $r>\lambda$, then $p(a^{*r}a^r)(\lambda)=0$, so $a^r\otimes 1$ belongs to the kernel of the form (\ref{innpro}). Set
$$
e_k:=\sqrt{k!\lambda!^{-1}}\ a^{\lambda-k}\otimes 1\ \mbox{for}\
k=0,\dots,\lambda\ \mbox{and}\ e_{k+\lambda}:=
\sqrt{\lambda!(\lambda+k)!^{-1}}\ a^{*k}\otimes 1\ \mbox{for}\ k\in\dN.%
$$
From (\ref{innpro}) and (\ref{aN}) we easily compute that $\langle e_k,e_n \rangle_0 = \delta_{kn}$ for $k,n\in \dN_0$. Hence $\{e_k;k \in \dN_0\}$ is an orthonormal base of $\cH_\lambda$. From the definition of $e_k$ we immediately obtain that
$$
a^*e_k=\sqrt{k+1} e_{k+1}\ \mbox{and}\ ae_k=\sqrt{k}e_{k-1}\
\mbox{for}\ k\in\dN_0,\ \mbox{where}\ e_{-1}:=0.
$$
This shows that for each $\lambda \in \dN_0$ the Hermitian $\cA$-module induced from the $\cB$-module $V_\lambda$ via $p$ is nothing but the \textit{Bargman-Fock representation} of the Weyl algebra.

If $\lambda \notin \dN_0$, the form (\ref{innpro}) is not positive semi-definite. Indeed, by (\ref{aN}) we have $\langle a\otimes1,a\otimes 1\rangle_0 =\lambda < 0$ if $\lambda <0$ and $\langle a^{k+1}\otimes 1, a^{k+1} \otimes 1 \rangle_0 =\lambda \cdots (\lambda{-} k+1)(\lambda -k)<0$ if $k-1<\lambda <k$ for $k\in\dN$.

Summarizing, we have shown that the $\cB$-module $V_\lambda$ is inducible to a Hermitian $\cA$-module if and only if $f(N)=f(\lambda)\geq 0$ for all $f \in \cB \cap \sum \cA^2$ or equivalently if $\lambda\in\dN_0$. \hfill $\Box$
\end{exam}

Our paper is organized in the following way. In Section
\ref{sect_rigged_induced} we study induced $*$-representations
defined by rigged modules. We follow mainly the approach given in
Chapter XI of \cite{fd} with some necessary modifications needed for
unbounded representations. As an application we show that the
well-behaved representations of $*$-algebras defined in \cite{s2} by
means of compatible pairs are induced representations coming from
certain rigged modules. Section \ref{sect_cond_exp} is concerned
with conditional expectations of general $*$-algebras. We give
various definitions depending on the corresponding positivity
conditions and develop a number of examples for these notions.
Section \ref{sect_grad_alg} is devoted to $G$-graded $*$-algebras
$\cA=\oplus_{g\in G}\cA_g$ for a discrete group $G.$ If $H$ is a
subgroup of $G$, then there exists a canonical conditional
expectation of $\cA$ on the $*$-subalgebra $\cA_H=\oplus_{h \in H}~
\cA_h$. Hence $*$-representations of $\cA_H$ can be induced to a
$*$-representations of $\cA.$ From Section \ref{sect_grad_alg_com}
on we are dealing with $G$-graded $*$-algebras $\cA=\oplus_{g\in
G}\cA_g$ for which the $*$-subalgebra $\cB:=A_e$ is commutative.
There is a large variety of $G$-graded $*$-algebras (Weyl algebra,
enveloping algebras of $su(2)$ and $su(1,1),$ quotients of the
enveloping algebra of the Virasoro algebra, $*$-algebras associated
with dynamical systems, quantum disc algebras, Podles' quantum
spheres, quantum algebras, and many others) that have this property.
In Section \ref{sect_sys_impr} we study systems of imprimitivity and
prove our first Imprimitivity Theorem. In Section
\ref{sect_grad_alg_com} we show that there is a partial action of
the group $G$ on the set $\cBp$ of all characters of the commutative
$*$-algebra $\cB$ which are nonnegative on the cone $\cB\cap\sum
\cA^2$. This partial action is used for a detailed study of the
inducing process from characters of the set $\cBp.$ In particular,
we characterize irreducible representations and equivalent
representations in terms of stabilizer groups of characters.

A fundamental problem in unbounded representation theory is to
define and characterize well-behaved representations of a general
$*$-algebra. In Section \ref{sect_well_beh} we develop a new concept of
well-behaved representations for $G$-graded $*$-algebras
$\cA=\oplus_{g\in G}\cA_g$ with commutative $*$-subalgebra $\cA_e.$
Among others we prove that well-behaved representations decompose
into direct sums of cyclic well-behaved representations. This theorem is technically rather involved and it is probably the deepest result of our paper. In Section
\ref{sect_well_beh_sys_impr} we define well-behaved systems of
imprimitivity and prove an Imprimitivity Theorem for well-behaved
representations. The next two sections of the paper are devoted to
detailed treatments of important examples. In Section
\ref{sect_env_alg} we study the enveloping algebras of three Lie
algebras. For the real Lie algebras $su(2)$ and $su(1,1)$ we prove
that the induced representations from characters of $\cBp$ are
precisely the representations $dU$, where $U$ is an irreducible
unitary representation of the Lie group $SU(2)$ resp. of the
universal covering group of $SU(1,1).$ For the enveloping algebra of
the Virasoro algebra we characterize irreducible $*$-representations
with finite-dimensional weight spaces as induced representations
from characters of $\cBp$. In Section \ref{sect_dyn_sys} we
investigate $*$-algebras associated with some dynamical systems.
For all these examples well-behaved representations
according to our definition in Section \ref{sect_well_beh} coincide
with distinguished "nice" representations of these
$*$-algebras thereby showing the usefulness of our concept of well-behavedness and emphasizing the role of positivity. In
Section \ref{sect_further_exam} we mention a number of other
examples for which our theory applies.

We close this introduction by collecting some definitions and
notations.

By a $*$-\textit{algebra} we mean a complex associative algebra
$\cA$ equipped with a mapping $a\mapsto a^*$ of $\cA$ into itself,
called the \textit{involution} of $\cA$, such that $(\lambda a+\mu
b)^* = \bar{\lambda}a^*+ \bar{\mu} b^*, (ab)^* = b^* a^*$ and
$(a^*)^*=a$ for $a,b\in \cA$ and $\lambda, \mu\in \dC$. The unit of
$\cA$ (if it exists) will be denoted by $\textbf{1}_\cA$ and the
group of all $*$-automorphisms of $\cA$ by $\Aut \cA$. We shall say
that a group $G$ \textit{acts as automorphism group on $\cA$} if
there is a group homomorphism $g\mapsto\alpha_g$ of $G$ into
$\Aut\cA.$ A subset $\cC$ of $\cA_h:=\{a\in\cA: a=a^*\}$ is called a
\textit{pre-quadratic module} if $\cC+\cC\subseteq\cC,\
\dR_+{\cdot}\cC\subseteq\cC$, and $a^*\cC a\in\cC$ for all
$a\in\cA$. A \textit{quadratic module} of $\cA$ is a pre-quadratic
module $\cC$ such that $\mathbf{1}_\cA\in\cC$ (see e.g. \cite{s4}).
The wedge
$$
\sum\cA^2:=\left\{\sum^n_{j=1} a_j^* a_j;\ a_1,{\dots},a_n\in\cA,
n\in\dN\right\}
$$
of all finite sums of squares is obviously the smallest quadratic
module of $\cA$.

Throughout this paper we use some terminology and results from
unbounded representation theory in Hilbert space (see e.g. in
\cite{s1}). In particular, we shall speak about $*$-representations
rather than Hermitian modules. Let us repeat some basic notions and
facts.

Let $\cD$ be a dense linear subspace of a Hilbert space $\cH$ with
scalar product $\langle\cdot,\cdot\rangle.$ A
$*$-\textit{representation} of a $*$-algebra $\cA$ on $\cD$ is an
algebra homomorphism $\pi$ of $\cA$ into the algebra $L(\cD)$ of
linear operators on $\cD$ such that
$\langle\pi(a)\varphi,\psi\rangle=\langle\varphi,\pi(a^*)\psi\rangle$
for all $\varphi,\psi\in\cD$ and $a\in \cA$. We call $\cD(\pi):=\cD$
the \textit{domain} of $\pi$ and write $\cH_\pi:=\cH$. Two
$*$-representation $\pi_1$ and $\pi_2$ of $\cA$ are
\textit{(unitarily) equivalent} if there exists an isometric linear
mapping $U$ of $\cD(\pi_1)$ onto $\cD(\pi_2)$ such that $\pi_2(a)=
U\pi_1(a)U^{-1}$ for $a\in \cA$. The \textit{direct sum
representation} $\pi_1\oplus\pi_2$ acts on the domain
$\cD(\pi_1)\oplus \cD(\pi_2)$ by $(\pi_1\oplus \pi_2)(a)=
\pi_1(a)\oplus \pi_2(a)$, $a\in\cA$. A $*$-representation $\pi$ is
called \textit{irreducible} if a direct sum decomposition $\pi
=\pi_1\oplus\pi_2$ is only possible when $\cD(\pi_1)=\{0\}$ or
$\cD(\pi_2)=\{0\}$. If $T$ is a Hilbert space operator, $\cD(T),
\Ran T,\ \overline{T}$ and $T^*$ denote its domain, its range, its
closure and its adjoint, respectively.

Suppose that $\pi$ is a $*$-representation of $\cA$. If $\cC$ is a
pre-quadratic module of $\cA$, $\pi$ is called $\cC$-positive if
$(\pi(c)\varphi,\varphi)\geq 0$ for all $c\in\cC$ and $\varphi\in
\cD(\pi)$. We denote by $\Res_\cB\pi$ the restriction of $\pi$ to
a $*$-subalgebra $\cB$. The \textit{graph topology} of $\pi$ is
the locally convex topology on the vector space $\cD(\pi)$ defined
by the norms $\varphi\mapsto\norm{\varphi}+\norm{\pi(a)\varphi},$
where $a\in\cA$. If $\cD(\overline{\pi})$ denotes the completion
the $\cD(\pi)$ in the graph topology of $\pi,$ then
$\overline{\pi}(a):=\overline{\pi(a)}\upharpoonright\cD(\overline{\pi}),$
$a\in\cA,$ defines a $*$-representation of $\cA$ with domain
$\cD(\overline{\pi}),$ called the \textit{closure} of $\pi.$ In
particular, $\pi$ is \textit{closed} if and only if $\cD(\pi)$ is
complete in the graph topology of $\pi.$ By a \textit{core for
$\pi$} we mean a dense linear subspace $\cD_0$ of $\cD(\pi)$ with
respect to the graph topology of $\pi$. A $*$-representation $\pi$
is called \textit{non-degenerate} if $\pi(\cA)\cD(\pi):= {\rm
Lin}~\{\pi(a)\varphi; a \in \cA, \varphi \in \cD(\pi)\}$ is dense
in $\cD(\pi)$ in the graph topology of $\pi$. If $\cA$ is unital
and $\pi$ is non-degenerate, then we have
$\pi(1_\cA)\varphi=\varphi$ for all $\varphi \in \cD(\pi)$. We say
that $\pi$ is \textit{cyclic} if there exists a vector $\varphi\in
\cD(\pi)$ such that $\pi(\cA)\varphi$ is dense in $\cD(\pi)$ in
the graph topology of $\pi$. Further, $\pi$ is called
\textit{self-adjoint} if $\cD(\pi)$ is the intersection of all
domains $\cD(\pi(a)^*)$, where $a \in \cA$. The \textit{(strong)
commutant} $\pi(\cA)^\prime$ consists of all bounded operators $T$
on $\cH_\pi$ such that $T\cD(T) \subseteq \cD(T)$ and
$\pi(a)T\varphi=T\pi(a)\varphi$ for $a \in \cA$. If $\pi$ is
self-adjoint, $\pi(\cA)^\prime$ is a von Neumann algebra. A closed
$*$-representation $\pi$ of a commutative $*$-algebra $\cB$ is
called \textit{integrable} if $\overline{\pi(b^*)}=\pi(b)^*$ for all
$b \in \cB$.

\section{Rigged Modules and Induced Representations}\label{sect_rigged_induced}%
\subsection{ }
Let $\cB$ be a $*$-algebra. From \cite{fd}, p. 1078, we repeat the
following
\begin{defn}\label{defn_right_B_rig}
\textit{A right $\cB$-rigged module} is a right $\cB$-module $\gX$
equipped with a map $\langle\cdot,\cdot\rangle_{\cB}:\gX\times\gX\to
\cB$ which is $\dC$-linear in the first variable and
$\dC$-anti-linear in the second variable and satisfies the following
conditions:
\begin{enumerate}
  \item[$(i)_{\ }$] $\langle x,y\rangle_\cB=(\langle y,x\rangle_\cB)^*$ for $x,y\in\gX,$
  \item[$(ii)_1$]  $\langle xb,y\rangle_\cB=\langle x,y\rangle_\cB b$ for $x,y\in\gX$ and $b\in \cB.$
\end{enumerate}
\end{defn}

\noindent Clearly, $(i)$ and $(ii)_1$ are equivalent to the
conditions $(i)$ and $(ii)_2,$ where
\begin{enumerate}
  \item[$(ii)_2$] $\ \ \langle x,yb\rangle_\cB=b^*\langle x,y\rangle_\cB\ \mbox{for}\ x,y\in\gX\ \mbox{and}\ b\in \cB.$
\end{enumerate}

\noindent Suppose that $(\gX,\langle\cdot,\cdot\rangle)$ is a right
$\cB$-rigged module. By $(ii)_1$ and $(ii)_2$ we have
\begin{enumerate}
\item[$(ii)_{\ }$] $\ \langle xb_1,yb_2\rangle_\cB=b_2^*\langle
x,y\rangle_\cB b_1\ \mbox{for}\ x,y\in\gX\ \mbox{and}\ b_1,b_2\in
\cB.$
\end{enumerate}

Suppose that $\rho$ is a $*$-representation of $\cB$ on
$(\cD(\rho),\langle\cdot,\cdot\rangle)$. Let
$\gX\otimes_\cB\cD(\rho)$ denote the quotient of the tensor product
$\gX\otimes\cD(\rho)$ over $\dC$ by the subspace
$$
\cN_{\rho}=\left\{\sum_{k=1}^r x_kb_k\otimes\varphi_k- \sum_{k=1}^r
x_k\otimes\rho(b_k)\varphi_k;\ x_k\in\gX,\ b_k\in \cB,\
\varphi_k\in\cD(\rho),\ r\in \dN\right\}.
$$

\begin{lemma}\label{lemma_pre_inner_prod}
\begin{gather}\label{eq_welldef1}
\langle\sum_k x_k\otimes\varphi_k,\sum_l
y_l\otimes\psi_l\rangle_0:=\sum_{k,l}(\rho(\langle
x_k,y_l\rangle_\cB)\varphi_k,\psi_l),
\end{gather}
where $x_k,y_l\in\gX$ and $\varphi_k, \psi_l\in\cD(\rho),$ is a
well-defined Hermitian sesquilinear form $\langle \cdot,\cdot
\rangle_0$ on the tensor products $\gX\otimes\cD(\rho)$ and
$\gX\otimes_\cB\cD(\rho).$
\end{lemma}
\noindent\textbf{Proof.} Obviously, $\langle\cdot,\cdot\rangle_0$ is
well-defined on the tensor product $\gX\otimes\cD(\rho)$ over $\dC$.
To prove that $\langle\cdot,\cdot\rangle_0$ is also well-defined on
the tensor
product $\gX\otimes_\cB\cD(\rho)$ 
it suffices to show that $\langle \zeta,\eta\rangle_0=0$ and
$\langle\eta,\zeta \rangle_0=0$ for arbitrary vectors $\eta=\sum
y_j\otimes\psi_j\in \gX\otimes\cD(\rho)$ and $\zeta=\sum_k
x_kb_k\otimes\varphi_k-\sum_k x_k\otimes\rho(b_k)\varphi_k\in
\cN_\rho.$ From $(ii)_1$ we obtain
$$\sum_{k,l}(\rho(\langle x_kb_k,y_l\rangle_\cB)
\varphi_k,\psi_l)=\sum_{k,l}(\rho(\langle
x_k,y_l\rangle_\cB)\rho(b_k) \varphi_k,\psi_l).
$$
Using condition $(i)$ it follows from the latter that $\langle
\zeta,\eta\rangle_0=0.$ Similarly, $(i)$ and $(ii)_2$ yield
$\langle\eta,\zeta\rangle_0=0.$ Condition $(i)$ implies that
$\langle\cdot,\cdot\rangle_0$ is Hermitian (that is $\langle
\zeta,\eta\rangle_0=\overline{\langle\eta,\zeta\rangle_0}$ for all
$\zeta,\eta\in\gX\otimes\cD(\rho)$ resp.
$\zeta,\eta\in\gX\otimes_\cB\cD(\rho).$) \hfill $\Box$\mn

Let $\cC$ be the set of finite sums of elements $\langle
x,x\rangle_\cB,$ where $x\in \cX.$ Then $\cC$ is a pre-quadratic
module of the $*$-algebra $\cB.$ Indeed, condition $(ii)$ implies
that $b^*cb\in \cC$ for $b\in \cB$ and $c\in\cC.$

Let $\Rep_c\cB$ denote the family of all direct sums of cyclic
$*$-representations of $\cB.$ Note that each cyclic
$*$-representation is obviously non-degenerate.

\begin{lemma}\label{lemma_nonneg_sesq_form}
If $\rho\in\Rep_c\cB$ and $\rho$ is $\cC$-positive, then
$\langle\cdot,\cdot\rangle_0$ is a nonnegative sesquilinear form on
$\gX\otimes_\cB\cD(\rho).$
\end{lemma}
\noindent\textbf{Proof.} Assume first that $\rho$ is a cyclic
representation with a cyclic vector $\xi\in\cD(\rho).$ Take
$\eta=\sum_{k=1}^n x_k\otimes\psi_k\in\gX\otimes_\cB\cD(\rho)$ and
fix $\varepsilon>0.$ Since $\xi$ is cyclic, there exist
$b_1,\dots,b_n \in \cB$ such that
$\norm{\rho(b_k)\xi-\psi_k}<\varepsilon$ and
$\norm{\rho(\langle{}x_k,x_l\rangle_\cB)(\rho(b_k)\xi-\psi_k)}<\varepsilon$
for all $k,l=1,\dots,n.$ Then for $k,l=1,\dots,n$ we get
\begin{gather*}
|\langle\rho(\langle{}x_k,x_l\rangle_\cB)\psi_k,\psi_l\rangle-\langle\rho(\langle{}x_k,x_l\rangle_\cB)\rho(b_k)\xi,\rho(b_l)\xi\rangle|\leq\\
\leq|\langle\rho(\langle{}x_k,x_l\rangle_\cB)\psi_k,\psi_l-\rho(b_l)\xi\rangle|+
|\langle\rho(\langle{}x_k,x_l\rangle_\cB)(\rho(b_k)\xi-\psi_k),\rho(b_l)\xi\rangle|\leq\\
\leq\norm{\rho(\langle{}x_k,x_l\rangle_\cB)\psi_k}\varepsilon+\norm{\rho(b_l)\xi}\varepsilon\leq
\norm{\rho(\langle{}x_k,x_l\rangle_\cB)\psi_k}\varepsilon+ \norm{\psi_l}\varepsilon + \varepsilon^2.
\end{gather*}
Therefore
$\langle\eta,\eta\rangle_0=\sum_{k,l=1}^n\langle\rho(\langle{}x_k,x_l\rangle_\cB)\psi_k,\psi_l\rangle$
can be approximated as small as we want by
\begin{gather*}
\sum_{k,l=1}^n\langle\rho(\langle{}x_k,x_l\rangle_\cB)\rho(b_k)\xi,\rho(b_l)\xi\rangle=
\sum_{k,l=1}^n\langle\rho(\langle{}x_kb_k,x_lb_l\rangle_\cB)\xi,\xi\rangle=
\langle\rho(\langle{}\sum_{k=1}^nx_kb_k,\sum_{k=1}^nx_kb_k\rangle_\cB)\xi,\xi\rangle,
\end{gather*}
which is nonnegative. This implies that $\langle\eta,\eta\rangle_0$ is
also nonnegative.

In the case when $\rho$ is a direct sum of cyclic representations
$\rho_i$ use the equality
$\gX\otimes_\cB\cD(\rho)=\sum_i\gX\otimes_\cB\cD(\rho_i).$ \hfill
$\Box$\mn

\noindent\textbf{Remark.} There is a counter-part of Lemma
\ref{lemma_nonneg_sesq_form} for $*$-representations $\rho$ of $\cB$
which are not necessarily direct sums of cyclic $*$-representations.
If $\rho$ is \textit{non-degenerate} and \textit{completely
positive} with respect to the corresponding matrix ordering (see
\cite{s1}, 11.1 and 11.2, for this concept), then the sesquilinear
form $\langle\cdot,\cdot\rangle_0$ is nonnegative on
$\gX\otimes\cD(\rho)$ resp. $\gX\otimes_\cB\cD(\rho).$

\subsection{ } Now let $\cA$ be another $*$-algebra.
\begin{defn}\label{defn_right_B_left_A}
\textit{A right $\cB$-rigged left $\cA$-module} is a right
$\cB$-rigged module $(\gX,\langle\cdot,\cdot\rangle_\cB)$ which is a
left $\cA$-module such that
\begin{enumerate}
  \item[$(iii)$] $\ \langle ax,y\rangle_\cB=\langle x,a^*y\rangle_\cB\ \mbox{for}\ a\in \cA,\ x,y\in\gX.$
\end{enumerate}

A \textit{right $\cB$-rigged $\cA{-}\cB$-bimodule} is a right
$\cB$-rigged left $\cA$-module satisfying
\begin{enumerate}
  \item[$(iv)$] $\ (ax)b=a(xb)\ \mbox{for}\ a\in \cA,\ b\in \cB,\ x\in\gX.$
\end{enumerate}
\end{defn}

\begin{lemma}\label{lemma_eq2_eq3}
Suppose $(\gX,\langle\cdot,\cdot\rangle_\cB)$ is a right
$\cB$-rigged left $\cA$-module (resp. $\cA-\cB$-bimodule). Then
\begin{gather}\label{eq2}
\pi_0(a)(\sum_k x_k\otimes\varphi_k)=\sum_k ax_k\otimes\varphi_k,\
a\in \cA,
\end{gather}
where $x_k\in\gX,\ \varphi_k\in\cD(\rho),$ is a well-defined
homomorphism of $\cA$ into the linear mappings of the vector space
$\gX\otimes\cD(\rho)$ (resp. $\gX\otimes_\cB\cD(\rho)$) such that
\begin{gather}\label{eq3}
\langle\pi_0(a)\zeta,\eta\rangle_0=\langle\zeta,\pi_0(a^*)\eta\rangle_0\ \mbox{for}\ a\in \cA,\ %
\zeta,\eta\in\gX\otimes\cD(\rho)\ \mbox{resp.}\
\zeta,\eta\in\gX\otimes_\cB\cD(\rho).%
\end{gather}
\end{lemma}
\noindent\textbf{Proof.} Since $\gX$ is a left $\cA$-module, $\pi_0$
is an algebra homomorphism into $L(\gX\otimes\cD(\rho)).$ Equation
(\ref{eq3}) follows then immediately by combining
(\ref{eq_welldef1}), (\ref{eq2}) and Definition
\ref{defn_right_B_left_A}, $(iv).$

If $\gX$ is an $\cA-\cB$-bimodule, $\pi_0$ is well-defined on
$\gX\otimes_\cB\cD(\rho),$ since by $(iv)$ we have
\begin{align*}
\pi_0(a)(\sum_k x_kb_k\otimes\varphi_k) &=\sum_k
a(x_kb_k)\otimes\varphi_k=\sum_k (ax_k)b_k\otimes\varphi_k\\
&=\sum_k ax_k\otimes\rho(b_k)\varphi_k=\pi_0(a)(\sum_k
x_k\otimes\rho(b_k)\varphi_k).
\end{align*}
\hfill $\Box$

\begin{lemma}\label{lemma4}
Suppose $\gX$ is a right $\cB$-rigged left $\cA$-module and $\rho$
is a $*$-representation of $\cB$ such that the sesquilinear form
$\langle\cdot,\cdot\rangle_0$ on $\gX\otimes_\cB \cD(\rho)$ is
nonnegative. Let $\langle\cdot,\cdot\rangle$ be the scalar product
on the quotient space
$\cD(\pi_0):=(\gX\otimes_\cB\cD(\rho))/\cK_\rho$ defined by
$\langle[\eta],[\zeta]\rangle=\langle\eta,\zeta\rangle_0$, where
$\cK_\rho:=\{\eta:\langle \eta,\eta\rangle_0=0\}$ and $[\eta]:=\eta
+\cK_\rho$. Then
$$\pi_0(a)[\eta]=[\pi_0(a)\eta],\ a\in \cA,\ \eta\in
\gX\otimes\cD(\rho),
$$
defines a $*$-representation $\pi_0$ of $\cA$ on the pre-Hilbert
space $(\cD(\pi_0),\langle\cdot,\cdot\rangle).$
\end{lemma}
\noindent\textbf{Proof.} Because of Lemma \ref{lemma_eq2_eq3} it
suffices to check that $\pi(a)$ is well-defined on $\cD(\pi_0),$
that is, $\pi_0(a)\cK_\rho\subseteq\cK_\rho.$ Let $\eta\in\cK_\rho.$
Using (\ref{eq3}) and the Cauchy-Schwarz inequality for the
nonnegative sesquilinear form $\langle\cdot,\cdot\rangle_0$ we
obtain
\begin{gather*}
\langle\pi_0(a)\eta,\pi_0(a)\eta\rangle_0=
\langle\eta,\pi_0(a^*)\pi_0(a)\eta\rangle_0=\langle\eta,\pi_0(a^*a)\eta\rangle_0\leq\\
\leq\langle\eta,\eta\rangle_0^{1/2}\langle\pi_0(a^*a)\eta,\pi_0(a^*a)\eta\rangle_0^{1/2}=0.
\end{gather*}
That is, $\pi_0(a)\eta\in\cK_\rho.$ \hfill $\Box$\mn

Let $\pi$ denote the closure of the $*$-representation $\pi_0$ from
Lemma \ref{lemma4}.

\begin{defn}
We say the $*$-representation $\pi$ of $\cA$ \textit{is induced from
the $*$-representation $\rho$ of $\cB$ via the right $\cB$-rigged
left $\cA$-module $\gX$} or simply $\pi$ is \textit{induced from}
$\rho.$ A $*$-representation $\rho$ of $\cB$ is called
\textit{inducible} (from $\cB$ to $\cA$) if the sesquilinear form
(\ref{eq_welldef1}) is nonnegative.
\end{defn}

We denote $\pi$ by $\Ind_{\cB\uparrow\cA}\rho$ or simply by
$\Ind\rho$ if no confusion can arise. The main assertions of the
preceding lemmas are summarized by the following proposition.
\begin{prop}\label{propind}
Suppose that $\cA$ and $\cB$ are $*$-algebras and $\gX$ is a right
$\cB$-rigged left $\cA$-module. If $\rho$ is a $*$-representation of
$\cB$ such that the sesquilinear form $\langle\cdot, \cdot
\rangle_0$ on $\gX\otimes\cD(\rho)$ given by (\ref{eq_welldef1}) is
nonnegative, then $\Ind\rho$ is a closed $*$-representation of $\cA$
defined on the core $(\gX\otimes\cD(\rho))/\cK_\rho$ by
$$
\Ind\rho (a)\left[\sum_k x_k\otimes\varphi_k\right]= \left[\sum_k
ax_k \otimes \varphi_k\right],\ \mbox{where}\ a\in\cA, x_k\in\gX,
\varphi_k \in \cD(\rho).
$$
If $\rho$ is a $\cC$-positive $*$-representation from $\Rep_c \cB$,
then the form $\langle\cdot, \cdot \rangle_0$ is nonnegative and
hence the induced representation $\Ind \rho$ exists. If $\gX$ is a
right $\cB$-rigged $\cA-\cB$--bimodule, then the core $(\gX \otimes
\cD(\rho))/\cK_\rho$ is a quotient of the tensor product $\gX
\otimes_\cB \cD(\rho)$.
\end{prop}

For applications the following proposition is convenient.
\begin{prop}\label{propind1}
Let $\cA$ and $\cB$ $*$-algebras and let $\cX$ be a right $\cB$-rigged left $\cA$-module. Let $\rho$ be a $*$-reprsentation of $\cB$. Assume that there exists a Hilbert space $(\cH_1,\langle{\cdot},{\cdot}\rangle_1)$ and a (well-defined) linear mapping $\Phi:\cX\otimes \cD(\rho) \to \cH_1$ such that $\cD_1:=\Phi (\cX\otimes \cD(\rho))$ is dense in $\cH_1$ and
\begin{align}\label{phicond}
\langle \Phi(x\otimes \varphi),\Phi(y \otimes \psi)\rangle_1 =(\rho(\langle x,y\rangle_\cB)\varphi, \psi),~~ x,y \in \cX,\varphi,\psi \in \cD(\rho).
\end{align}
Then $\rho$ is inducible and $\Ind\rho$ is unitarily equivalent to the closure of the $*$-representation $\pi_1$ on $\cD_1$ defined by $\pi_1(a)(\Phi(x\otimes \varphi))= \Phi(ax \otimes \varphi)$, where $a\in \cA, x \in\cX,\varphi \in \cD(\rho)$.
\end{prop}
\noindent\textbf{Proof.} Define a the linear mapping $U$ of $\cX\otimes \cD(\rho)$ onto $\cD_1$ by
$U(\eta)=\Phi(\sum_k x_k\otimes \varphi_k)$ for $\eta =\sum_k x_k\otimes \varphi_k$. Comparing (\ref{eq_welldef1}) and (\ref{phicond}) we see that the form (\ref{eq_welldef1}) is nonnegative, so $\rho$ is inducible. Further it follows that $\eta \in \cK_\rho$
if and only if $\Phi(\sum_k x_k\otimes \varphi_k)=0$. Hence $U$ yields an isometric linear mapping, denoted again by $U$, of the unitary space $ ((\gX\otimes\cD(\rho))/\cK_\rho,\langle{\cdot},{\cdot}\rangle) $ onto the unitary space $(\cD_1,\langle{\cdot},{\cdot}\rangle_1)$ such that $\pi_1(a)=U~ \Ind\rho(a) U^{-1}$, $a\in \cA$.
\hfill $\Box$

\mn
\noindent\textbf{Remark.} Above we have defined induced representations for a right $\cB$-rigged left $\cA$-module $\cX$.
However, except for Example \ref{compair} in all applications below $\cX$ is even a right $\cB$-rigged $\cA{-}\cB$-bimodule. Moreover, if $\cX$ is a right $\cB$-rigged left $\cA$-module, then using the axioms $(ii)_1$ and $(iii)$ we compute
$$
\langle (ax)b-a(xb), y\rangle_\cB= \langle ax,y\rangle_\cB b - \langle xb,a^*y\rangle_\cB= \langle x,a^* y\rangle_\cB b-\langle x,a^*y\rangle _\cB b=0.
$$
for arbitrary $a\in \cA$, $b \in\cB$ and $x,y \in \cX$. That is, all elements $(ax)b{-}a(xb)$ are annihilated by $\cX$ with respect to the $\cB$-valued form $\langle{\cdot},{\cdot}\rangle_\cB$. In particular, if this form is nondegenerate, then the right $\cB$-rigged left $\cA$-module $\cX$ is a right $\cB$-rigged $\cA{-}\cB$-bimodule.

\mn
The following lemma is needed in Section \ref{sect_well_beh} below.

\begin{lemma}\label{lemma_cyc_vec_of_ind_rep}
Suppose $\gX$ is a right $\cB$-rigged left $\cA$-module (resp.
$\cA-\cB$-bimodule) and $\rho$ is an inducible cyclic
$*$-representation of $\cB$ with cyclic vector $v\in\cD(\rho).$ Then
the linear subspace of vectors $[x\otimes v]$, where $x\in\gX$, is a
core of $\pi=\Ind\rho.$
\end{lemma}
\noindent\textbf{Proof.} It suffices to show that for arbitrary
$\varepsilon>0,$ $a \in \cA$, $x\in\gX$, and $w \in \cD(\rho)$ there
exists $b\in\cB$ such that $\norm{\pi(a)([x\otimes w]-[x\otimes
\rho(b)v])}<\varepsilon$. Since $v$ is cyclic, there is a $b\in\cB$
such that $\norm{\rho(\langle
ax,ax\rangle_\cB)(\rho(b)v-w)}<\varepsilon$ and
$\norm{\rho(b)v-w}<\varepsilon.$ Using the Cauchy-Schwarz inequality we
get
\begin{gather*}
\norm{\pi(a)([x\otimes w]-[x\otimes \rho(b)v])}^2=\norm{[ax\otimes(w-\rho(b)v)]}^2\\
=\langle\rho(\langle ax,ax\rangle_\cB)(w-\rho(b)v),(w-\rho(b)v)\rangle_0<\varepsilon ^2.
\end{gather*}
\hfill $\Box$\mn

The next lemma is a standard fact about induced
representations. We omit its simple proof.
\begin{lemma}\label{lemma_ind_of_sum}
Suppose $\gX$ is a right $\cB$-rigged left $\cA$-module (resp.
$\cA-\cB$-bimodule) and $\rho$ is a $*$-representation of $\cB.$
Assume that $\rho$ is a direct sum of representations $\rho_i, i\in
I.$ Then $\rho$ is inducible if and only if each $\rho_i$ is
inducible. Moreover, $\Ind\rho=\oplus_{i\in{}I}\Ind\rho_i.$
\end{lemma}

We close this section by showing that the considerations of \cite{s2} fit
nicely into the theory of induced representations.
\begin{exam}\label{compair}
\textit{Compatible pairs in the sense of ~\cite{s2}}\\ Let $\cA$
and $\cB$ be two $*$-algebras. Following \cite{s2}, we call
$(\cA,\cB)$ a \textit{compatible pair} if $\cB$ is a left
$\cA$-module, with a left action denoted by $\rhd,$ such that
\begin{gather}\label{eq_action}
(a\rhd b)^*c=b^*(a^*\rhd c)\ \mbox{for}\ a\in \cA\ \mbox{and}\ b\in
\cB.
\end{gather}

Now let $(\cA,\cB)$ be such a compatible pair. We equip $\gX=\cB$
with the $\cB$-valued sesquilinear form $\langle
b,c\rangle_\cB:=c^*b,\ b,c\in \cB,$ and with the right $\cB$-action
given by the multiplication. Then
$(\gX,\langle\cdot,\cdot\rangle_\cB)$ is a right $\cB$-rigged left
$\cA$-module. Indeed, axioms $(i)$ and $(ii)_2$ are obvious. Axiom
$(iii)$ follows from (\ref{eq_action}), since for arbitrary $a\in
\cA$ and $b,c\in \cB$ we have
$$
\langle a\rhd b,c\rangle_\cB=c^*(a\rhd b)=(a^*\rhd c)^*b=\langle
b,a^*\rhd c\rangle_\cB.
$$
Suppose that $\rho \in {\rm Rep}_c~\cB$. Since bounded
$*$-representations acting on the whole Hilbert space are obviously
in ${\rm Rep}_c~\cB$, this covers all representations of $\cB$
considered in \cite{s2}. Since the pre-quadratic module $\cC$ for
the form $\langle\cdot,\cdot \rangle_\cB$ is $\sum \cB^2$, $\rho$ is
$\cC$-positive. Therefore, by Proposition \ref{propind}, $\rho$
induces a $*$-representation $\pi{=}{\rm Ind}\rho$ of $\cA$. We
shall give a more explicit description of this representation $\pi$
expressed by formula (\ref{desc}) below.

Clearly, an element $\zeta=\sum b_k\otimes\varphi_k\in\gX\otimes
\cD(\rho)$ belongs to the kernel $\cK_{\rho}$ of the sesquilinear
form $\langle\cdot,\cdot\rangle_0$ if and only if
$$
\langle\zeta,\zeta\rangle_0=\sum_{k,l}\langle\rho(\langle
b_k,b_l\rangle_\cB)\varphi_k,\varphi_l\rangle=\langle\sum_k\rho(b_k)\varphi_k,\sum_l\rho(b_l)\varphi_l\rangle=0
$$
or equivalently if $\sum_k\rho(b_k)\varphi_k=0.$ Hence $\cK_\rho$ is
the kernel of the mapping
$$
\cB\otimes\cD(\rho)\ni\sum_k
b_k\otimes\varphi_k\mapsto\sum_k\rho(b_k)\varphi_k\in\rho(\cB)\cD(\rho),
$$
so we have an isomorphism of vector spaces
$\cD(\pi_0)=(\cB\otimes\cD(\rho))/\cK_\rho$ and
$\rho(\cB)\cD(\rho)$. If we identify $\cD(\pi_0)$ and
$\rho(\cB)\cD(\rho)$ by identifying $b\otimes\varphi$ and
$\rho(b)\varphi$, then we have
\begin{align}\label{desc}
\pi(a)(\sum_k\rho(b_k)\varphi_k)=\pi_0(a)(\sum_k\rho(b_k)\varphi_k)=\sum_k\rho(a\rhd
b_k)\varphi_k
\end{align}
for $a\in \cA.$ This formula shows that the $*$-representation $\pi_0$ and its closure $\pi =\Ind \rho$ as defined above are precisely the $*$-representations $\tilde{\rho}$ and $\rho^\prime$ as defined in \cite{s2}, Proposition 1.1. That is, \textit{all well-behaved $*$-representations $\rho^\prime$ of $\cA$ associated with the compatible pair $(\cA,\cB)$ in the sense of \cite{s2} are induced $*$-representations $\Ind~\rho$}. Note that the well-behaved $*$-representations in the sense of \cite{s2} are closely related to representations constructed from unbounded $C^*$-seminorms (see \cite{apt}, Chapter 8, for details).

In \cite{s2} a number of examples of compatible pairs are developed. A typical example of a compatible pair $(\cA,\cB)$ is obtained as follows: $\cB$ is the $*$-algebra $C_0^\infty (G)$ of a Lie group $G$ with convolution multiplication, $\cA$ is the enveloping algebra $\cU(g)$ of the Lie algebra $\gog$ of $G$ and $x\rhd f$ is the action of $x \in \cU(\gog)$ as a right-invariant differential operator on $f \in C_0^\infty(G)$. Note that as in all other examples of compatible pairs treated in \cite{s2} the $*$-algebra $\cB$ has no unit.

Moreover, all examples described in \cite{s2} are of the following form: $\cA$ and $\cB$ are $*$-subalgebras of a common unital
$*$-algebra $\gA$ and the left action of $a\in\cA$ on $b\in\cB$ is just the multiplication in the larger algebra $\gA.$ In this case it
follows at once from the $*$-algebra axioms that condition (\ref{eq_action}) is valid and that $(\gX,\langle\cdot,\cdot\rangle_\cB)$ is a right $\cB$-rigged $\cA-\cB$--bimodule. \hfill $\edex$
\end{exam}

\section{Conditional expectations}\label{sect_cond_exp}
In the rest of this paper we assume that $\cB$ is a unital $*$-subalgebra of a unital $*$-algebra $\cA$.

Most examples of rigged modules are derived from conditional expectations. This is a fundamental concept for this paper. Since positivity will play a crucial role in what follows, we require various versions of this notion.

\begin{defn}\label{defn_cond_exp}
A linear map $p:\cA\to\cB$ is called a \textit{conditional expectation} of $\cA$ onto $\cB$ if
\begin{enumerate}
  \item[$(i)$] $p(a^*)=p(a)^*,\ p(b_1ab_2)=b_1p(a)b_2\ \mbox{for all}\ a\in\cA,\ b_1,b_2\in\cB,\ p(\mathbf{1}_\cA)=\mathbf{1}_\cB,$\\
  and $p$ is positive in the sense that
  \item[$(ii)$] $p(\sum\cA^2)\subseteq\sum\cA^2\cap\cB.$
\end{enumerate}

A linear map $p$ satisfying only condition $(i)$ is called a $\cB$-\textit{bimodule projection} of $\cA$ onto $\cB.$

A conditional expectation $p$ will be called a \textit{strong conditional expectation} if
\begin{enumerate}
  \item[$(ii)_1$] $p(\sum\cA^2)\subseteq\sum\cB^2.$
\end{enumerate}

Let $\cC_\cA$ and $\cC_\cB$ be pre-quadratic modules of $\cA$ resp.
$\cB.$ A $\cB$-bimodule projection $p$ will be called
$(\cC_\cA,\cC_\cB)$-\textit{conditional expectation} of $\cA$ onto
$\cB$ if
\begin{enumerate}
  \item[$(ii)_2$] $p(\cC_\cA)\subseteq\cC_\cB.$
\end{enumerate}
\end{defn}

Note that axiom $(i)$ implies that any $\cB$-\textit{bimodule
projection} of $\cA$ onto $\cB$ is indeed a projection of $\cA$ onto
$\cB.$

The bridge of these notions to rigged modules is given by the
following simple lemma.
\begin{lemma}
Suppose that $p:\cA\to \cB$ is a $\cB$-bimodule projection of $\cA$
onto $\cB$ and define $\langle b,c\rangle_\cB:=p(c^*b)$ for $b,c\in
\cB$ and $\gX:=\cA$. Then $(\gX,\langle\cdot,\cdot\rangle_\cB)$ is a
right $\cB$-rigged $\cA-\cB$-bimodule with left and right actions
given by the multiplications in $\cA$.
\end{lemma}
\noindent\textbf{Proof.} Conditions $(i)$, $(ii)_1$, $(iii)$ and
$(iv)$ in Definitions \ref{defn_right_B_rig} and
\ref{defn_right_B_left_A} follow at once from $(i)$ in Definition
\ref{defn_cond_exp} and the $*$-algebra axioms. For instance, we
verify $(ii)_1$. If $x,y \in \gX(=\cA)$ and $b\in\cB$, then using
axiom $(i)$ in Definition \ref{defn_cond_exp} we have $\langle
xb,y\rangle_\cB=p(y^*xb)=p(y^*x)b=\langle x,y\rangle_\cB b.$ \hfill
$\Box$

\begin{defn}
A $\cB$-bimodule projection $p$ of $\cA$ onto $\cB$ is called
\textit{faithful} if $p(x^*x)=0$ for some $x\in\cA$ implies that
$x=0.$
\end{defn}

The next lemma illustrates the importance of this notion.
\begin{lemma}
Suppose that $p$ is a faithful $\cB$-bimodule projection of $\cA$ onto $\cB$. Let $\pi_i,\ i\in I,$ be a family of inducible $*$-representations of $\cB$ which separates the elements of $\cB$. Then the family $\Ind\pi_i,\ i\in I,$ separates the elements of $\cA.$
\end{lemma}
\noindent\textbf{Proof.} Let $a\in \cA$, $a\neq 0$. Since $p$ is
faithful, $p(a^*a)\neq 0.$ Since the family $\pi_i,i \in I,$
separate the elements of $\cB$, there exist a representation
$\pi_{i_0}, i_0\in I,$ and a vector $\varphi\in\cD(\pi_{i_0})$ such
that $\pi_{i_0}(p(a^*a))\varphi \neq 0.$ Then we have
$\norm{\Ind\pi_{i_0}(a)[1\otimes\varphi]}=
\norm{\pi_{i_0}(p(a^*a))\varphi}\neq 0.$ \hfill $\Box$\mn

The following simple proposition is taken from \cite{ver}. It characterizes a $\cB$-bimodule
projection in terms of its kernel.
\begin{prop}\label{charbim}
There exists a $\cB$-bimodule projection from $\cA$ onto $\cB$ if
and only if there exists a $*$-invariant subspace $\cT\subseteq\cA$
such that $\cA=\cB\oplus\cT$ and
\begin{gather}\label{eq_BTB_in_T}
\cB\cT\cB\subseteq\cT.
\end{gather}
If this is true, the $\cB$-bimodule projection $p$ is uniquely
defined by the requirement $\ker p=\cT$ and we have $p(\sum\cA^2)=
\sum\cB^2+p(\sum\cT^2).$
\end{prop}
\noindent\textbf{Proof.} Let $p$ be a $\cB$-bimodule projection from
$\cA$ onto $\cB$ and put $\cT=\ker p.$ For $t\in\cT$ and
$b_1,b_2\in\cB$ we have $p(b_1tb_2)=b_1p(t)b_2=0$ and
$p(t^*)=p(t)^*=0$, so that $\cT$ satisfies (\ref{eq_BTB_in_T}) and
is $*$-invariant. For arbitrary $a\in\cA$ we have $p(a)\in\cB$ and
$a-p(a)\in\cT,$ so that $\cA=\cB\oplus\cT.$

Conversely, if $\cT$ is given, one easily checks that the linear map
$p$ defined by $p(b)=b,\ b\in\cB,$ and $p(t)=0,\ t\in\cT,$ is indeed
a $\cB$-bimodule projection. \hfill $\Box$\mn

In the remaining part of this section we develop a number of
examples. In the first example we use Proposition \ref{charbim} to
show that there is no $\cB$-bimodule projection.

\begin{exam}
Let $\cA$ be the Weyl algebra from Example \ref{exam_weyl_alg}. As
it is well-known, the Hermitian elements $P=\frac{1}{\sqrt{2}} i
(a^*{-}a)$ and $Q=\frac{1}{\sqrt{2}} (a^*{+}a)$ satisfy the
commutation relation $PQ-PQ=-i$.

We show that \textit{there is no $\cB$-bimodule projection of $\cA$
onto $\cB:=\dC[P].$} Assume to the contrary there is such a
projection $p$ and let $\cT$ be its kernel. Then, since
$\cA=\cB\oplus\cT,$ there exists a polynomial $f\in\dC[t]$ such that
$Q+f(P)\in\cT.$ By (\ref{eq_BTB_in_T}) we have $PQ+Pf(P)$ and
$QP+f(P)P\in\cT$ which implies that $PQ-QP=-i\in\cT.$ Hence
$\mathbf{1}_\cA\in\cT$ and so $p=0$ which is a contradiction.

Using Proposition \ref{charbim} one can check that the map
$p$ defined in Example \ref{exam_weyl_alg} is the unique
$\cB$-bimodule projection from $\cA$ onto $\cB:=\dC[N].$ \hfill
$\edex$
\end{exam}

\begin{exam}\label{exam_dec_un}
Let $q_1,\dots,q_n \in \cA$ be a decomposition of unit of the
unital $*$-algebra $\cA$, that is, $q_1+ \cdots +q_n=1$ and
$q_i=q_i^2=q_i^*$ for $i=1,{\dots},n$. It is not difficult to show
that $q_iq_j=0$ for all $i\neq j$ and that the map
$$
p:a\mapsto q_1aq_1+\dots+q_naq_n
$$
is a conditional expectation of $\cA$ onto the $*$-subalgebra
$\cB=\{b\in \cA:b=p(b)\}.$ If $\cA$ is an $O^*$-algebra, then $p$ is
faithful. \hfill $\edex$
\end{exam}

\begin{exam}\label{exam_cond_exp_CG_CH}
Suppose that $G$ is a discrete group and $H$ is a subgroup of $G$.
Let $\cA=\dC[G]$ and $\cB=\dC[H]$ be the group algebras of $G$ and
$H$, respectively. Recall that the group algebra $\dC[G]$ of a
discrete group $G$ is a unital $*$-algebra with multiplication given
by the convolution and involution determined by the inversion of
group elements. More precisely, $\dC[G]$ is a complex vector space
with basis given by the group elements of $G$ and the product of two
base element $g$ and $h$ is just the group product $g h$ and $g^*$
is the inverse $g^{-1}.$ Let $p$ be the canonical projection of
$\dC[G]$ onto $\dC[H]$ defined
by $p(g)=g$ if $g \in H$ and $p(g)=0$ if $g \notin H$. 
\begin{prop}\label{prop_cond_exp_CG_CH}
$p$ is a faithful strong conditional expectation of $\dC[G]$ onto
$\dC[H].$
\end{prop}
\noindent\textbf{Proof.} It is clear from its definition that $p$
satisfies condition $(i)$ of the Definition \ref{defn_cond_exp}, so
$p$ is a $\dC[H]$-bimodule projection.

We shall prove that $p(\sum\dC[G]^2)\subseteq\sum \dC[H]^2.$ Let us
fix precisely one element $k_t\in G$ in each left coset $t\in G/H.$
Take an arbitrary element $a=\sum_{g\in G}\theta_g g$ of the group
algebra $\dC[G].$ Then there exist elements $a_i\in\dC[H]$, $i\in
G/H,$ such that $a=\sum_{g\in G}\theta_g g=\sum_{i\in G/H}k_ia_i.$
If $i,j\in G/H$ and $i{\neq} j$, then $k_i^{-1}k_j \notin H$ and
hence
$p(k_i^{-1}k_j)=0.$ Using this fact 
we obtain
\begin{gather*}
p(a^*a)=p\left(\left(\sum_{i\in G/H}k_ia_i\right)^*\left(\sum_{j\in
G/H}k_ja_j\right)\right)=p\left(\sum_{i,j\in
G/H}a_i^*k_i^{-1}k_ja_j\right)=\\
\sum_{i,j\in G/H}p(a_i^*k_i^{-1}k_ja_j)=\sum_{i,j\in
G/H}a_i^*p(k_i^{-1}k_j)a_j=\sum_{i \in G/H}a_i^*a_i,
\end{gather*}
so $p(a^*a)\in\sum\dC[H]^2.$ That is, $p$ is a strong conditional
expectation.

From the preceding equality it follows also that $p$ is faithful.
Indeed, if $p(a^*a)=0,$ then $\sum_i a_i^*a_i^{}=0$ which implies
that $a_i=0$ for all $i\in G/H$ and hence $a=0$. \hfill $\Box\
\edex$
\end{exam}

A large source of conditional expectations is obtained from groups
of $*$-automorphisms. The idea is taken from the
following standard construction of conditional
expectations of $C^*$-algebras reproduced from \cite{rief}, Example
1.5.
\begin{exam}\label{exam_cond_exp_cstar_alg}
Suppose that $\cA$ is a $C^*$-algebra and $G$ is a compact group
such that there is a continuous action $g\mapsto\alpha_g$ of $G$
as automorphism group of $\cA$. Let $dg$ denote the normalized
Haar measure of $G$. Then the map
$$a\mapsto\int_{G}\alpha_g(a)dg,\ a\in \cA,
$$
is a \textit{strong conditional expectation} of $\cA$ onto the
$C^*$-subalgebra $\cB$ of stable elements. \hfill $\edex$
\end{exam}

We now generalize this example to the case of
general $*$-algebras.

\begin{exam}\label{exam_loc_fin_comp_gr}
Suppose that $G$ is a compact group which acts by $*$-automorphisms
$\alpha_g,\ g\in G,$ on a $*$-algebra $\cA.$ Assume in addition that
the action is \textit{locally finite-dimensional}, that is, for
every $a\in\cA$ there exists a finite-dimensional linear subspace
$V\subset\cA$ such that $a\in V$, $\alpha_g(V)\subseteq V$ for all
$g\in G$, and the map $g \to \alpha_g(a)$ of $G$ into $V$ is
continuous. Then the mapping $p$ given by
\begin{align}\label{defpG}
p(a)=\int_{G}\alpha_g(a)dg,\ a\in\cA,
\end{align}
is well-defined. One easily verifies that $p$ is a $\cB$-bimodule
projection from $\cA$ onto the $*$-subalgebra
$\cB:=\{a\in\cA:\alpha_g(a)=a$ for all $g\in G\}$ of stable
elements.

Every $G$-invariant finite-dimensional subspace $V\subseteq\cA$ is
a unitarizable $G$-module. Since $G$ is compact, $\cA$
is a direct sum of submodules $\cA_t,\ t\in\widehat{G},$ where
$\cA_t$ denotes the direct sum of submodules in $\cA$ isomorphic
to $t\in\widehat{G}.$ In the case when $\cA$ is a $C^*$-algebra,
the subspaces $\cA_t,\ t\in\widehat{G},$ are called
\textit{spectral subspaces,} see e.g. \cite{stormer} and
\cite{evans}. The mapping $p$ is nothing but the projection
of the direct sum $\cA=\oplus_{t\in\widehat{G}}\cA_t$ onto the spectral
subspace $\cA_0$ corresponding to the trivial representation.

An analogue of the map $p$ was considered in \cite{cimpric}. Suppose
$R$ is a real closed field, $R[V]$ is the coordinate ring of an
affine variety $V$ and $G$ is a linear algebraic group over $R$
acting on $R[V].$ If $G$ is reductive, there is a canonical
projection $\rho$ from $R[V]$ onto the subring $R[V]^G$ of
$G$-invariants called \textit{Reynolds operator} (see
\cite{cimpric} for details). In the case when $G(R)$
semi-algebraically compact, Corollary 3.6 in
\cite{cimpric} states that $\rho(\sum
R[V]^2)\subseteq\sum R[V]^2.$

\begin{prop}\label{locfin}
The map $p$ defined by (\ref{defpG}) is a conditional expectation of
$\cA$ onto $\cB$.
\end{prop}
\noindent\textbf{Proof.} It remains to show that
$p(\sum\cA^2)\subseteq\sum\cA^2.$ Let $a\in\cA.$ Then there is a
finite-dimensional $G$-invariant subspace $V$ of $\cA$ containing
$a$. Then $V$ is a finite direct sum of submodules $V^{(t)},$ where
$V^{(t)}$ is multiple of $t\in\widehat{G}.$ Fix $t \in \widehat{G}$
and let $V^{(t)}=\oplus_i^{}V_i^{(t)}$ be a decomposition of
$V^{(t)}$ into a direct sum of irreducible $G$-modules. We can
choose an orthonormal base $a_{ij}^{(t)}$ in each space $V_i^{(t)}$
such that the matrices corresponding to $\alpha_g$ are unitary and
equal for all $i$, i.e. we have
\begin{gather*}
\alpha_g(a_{ij}^{(t)})=\sum_k u_{kj}^{(t)}(g)a_{ik}^{(t)},\ g\in
G,\ t\in\widehat{G}.
\end{gather*}
Let us fix elements $a_{i_1j_1}^{(t)},\ a_{i_2j_2}^{(s)}\in
V\subseteq\cA.$ Using the orthogonality relations of matrix
elements $u_{kj_1}^{(t)}$ and $u_{mj_2}^{(s)}$ on the compact
group $G$ we compute
\begin{gather*}
p((a_{i_1j_1}^{(t)})^*a_{i_2j_2}^{(s)})=\int\left(\sum_k
\overline{u_{kj_1}^{(t)}(g)}(a_{i_1k}^{(t)})^*\right)\cdot\left(\sum_m u_{mj_2}^{(s)}(g)a_{i_2m}^{(s)}\right)dg=\\
=\sum_{k,m}\int\overline{u_{kj_1}^{(t)}(g)}u_{mj_2}^{(s)}(g)dg\cdot\left(a_{i_1k}^{(t)}\right)^*a_{i_2m}^{(s)}=\\
=\frac {\delta_{ts}~\delta_{j_1j_2}}{\dim t}\sum_k\left(a_{i_1k}^{(t)}\right)^*a_{i_2k}^{(t)}.%
\end{gather*}
Since $a\in V$, we can write $a$ as a finite sum $
a=\sum_{i,j,t}\lambda_{ij}^{(t)}a_{ij}^{(t)}$, where
$\lambda_{ij}^{(t)}\in\dC.$ Applying the preceding equality we
obtain
\begin{align*}
p(a^*a)&=p\left(\sum_{i,j,t}\overline{\lambda_{ij}^{(t)}}\left(a_{ij}^{(t)}\right)^*\cdot\sum_{k,l,s}\lambda_{kl}^{(s)}a_{kl}^{(s)}\right)=
\sum_{j,t}p\left(\sum_{i}\overline{\lambda_{ij}^{(t)}}\left(a_{ij}^{(t)}\right)^*\cdot
\sum_{k}\lambda_{kj}^{(t)}a_{kj}^{(t)}\right)=\\
&=\sum_{j,t}\frac 1{\dim
t}\left(\sum_{i}\lambda_{ij}^{(t)}a_{ij}^{(t)}\right)^*\cdot
\left(\sum_{k}\lambda_{kj}^{(t)}a_{kj}^{(t)}\right)\in\sum\cA^2.\ \ \ \ \ \ \ \ \ \ \ \ \ \ \ \ \ \ \ \ \ \ \ \ \ \ \ \ \ \ \ \ \ \hfill\Box%
\end{align*}

\noindent In general this conditional expectation $p$ is not strong,
i.e. $p(\sum\cA^2)$ is not contained in $\sum\cB^2$. \hfill $\edex$
\end{exam}

\section{Group graded $*$-Algebras}\label{sect_grad_alg}
The algebraic representation theory of group graded algebras has
been extensively studied, see e.g. the books \cite{nastas} and
\cite{markus}. The monograph \cite{fd} deals with
$*$-algebraic bundles which can be considered as generalizations
of $G$-graded $*$-algebras to the case when $G$ is a topological
group. However, in \cite{fd} only
 bounded Hilbert space representations are treated. As we shall see below,
there are a plenty of important $G$-graded $*$-algebras (Weyl
algebra, enveloping algebras etc.) for which most
$*$-representations are unbounded.

\begin{defn}\label{defn_graded_nc}
Let $G$ be a (discrete) group. A \textit{$G$-graded $*$-algebra}
is a $*$-algebra $\cA$ which is a direct sum $\cA=\bigoplus_{g\in
G}\cA_g$ of vector spaces $\cA_g$, $g \in G$, such that
\begin{gather}\label{eq_Ggrad}
\cA_g\cdot\cA_h\subseteq \cA_{g\cdot h}\ \mbox{and}\
(\cA_g)^*\subseteq\cA_{g^{-1}}\ \mbox{for}\ g,h\in G.
\end{gather}
\end{defn}

From the two conditions in (\ref{eq_Ggrad}) it follows that a
$G$-grading of a $*$-algebra $\cA$ is completely determined if we
know the corresponding components for a set of generators of the
algebra $\cA$. In what follows we shall describe most of our $G$-gradings
of $*$-algebras in this manner.
\begin{exam}\label{semienv}
In this example we use some basics from the theory of semi-simple
Lie algebras. All facts we need can be found in the monograph
\cite{dix1}, 7.0 and 7.4.1. Suppose that $\gog$ is a semi-simple
complex Lie algebra. We denote by $\goh$ a Cartan subalgebra, by $Q$
the root lattice and by $H_1,\dots,H_l,
X_{-\alpha_1},\dots,X_{-\alpha_n},X_{\alpha_1},\dots,X_{\alpha_n}$ a
Cartan-Weyl basis of the Lie algebra $\gog$. If we consider the
complex universal enveloping algebra $\cU(\gog)$ of $\gog$ as a
$\gog$-module and so as an $\goh$-module by the adjoint
representation, we obtain a direct sum decomposition $
\cU(\gog)=\sum_{\lambda\in Q} \cU(\gog)_\lambda. $ This means that
$\cU(\gog)$ is a $G$-graded algebra, where $G$ is the abelian group
$Q.$ If $\cU(\gog)$ is equipped with an involution such that
$(X_{\alpha_j})^*{=}\varepsilon_j X_{-\alpha_j}$ and $(H_k)^*=H_k$
for all $j,k$, where $\varepsilon_j \in \{1,-1\}$, then we have
$(\cU(\gog)_\lambda)^*=\cU(\gog)_{-\lambda}$ and hence $\cU(\gog)$
is a $Q$-graded $*$-algebra. The algebra $\cU(\gog)_0$ is just the
commutant of the Cartan algebra $\goh$ in $\cU(\gog)$. Its structure
is descibed in \cite{dix1}, 7.4.2. \hfill $\edex$
\end{exam}

\begin{exam}\label{freepol}
Let $\cF=\dC\langle z_1,\dots,z_d,w_1,\dots,w_d \rangle$ be the free
polynomial algebra with generators $z_1,\dots,z_d,w_1,\dots,w_d$ and
involution determined by $(z_j)^*{=}w_j$, $j=1,\dots,d$. Then $\cF$
is a $\dZ$-graded $*$-algebra with $\dZ$-grading given by $z_j \in
\cF_1$. \hfill $\edex$
\end{exam}

To derive further examples we shall use the following lemma. We omit
its simple proof.
\begin{lemma}\label{lemma_factorgrad}
If $\cF=\bigoplus_{g\in G}\cF_g$ is a $G$-graded $*$-algebra and
$\cJ$ is a two-sided $*$-ideal of $\cF$ generated by subsets of
$\cF_g$, $g\in G,$ then the quotient $*$-algebra $\cF/\cJ$ is also
$G$-graded.
\end{lemma}


The proofs of the existence of gradings for all examples occuring
in this paper follow by the same pattern: We first define the
corresponding grading on the free $*$-algebra (Example
\ref{freepol}). If the polynomials of the defining relations
belong to single components of this grading, Lemma
\ref{lemma_factorgrad} applies and gives the grading of the
$*$-algebra. We illustrate this by a number of examples in the
last section.

Throughout the rest of this section $G$ is a discrete group with
unit element $e,\ H$ denotes a subgroup of $G$ and
$\cA=\bigoplus_{g\in G}\cA_g$ is a unital $G$-graded $*$-algebra.
The subspace $\cA_e$ is a $*$-subalgebra of $\cA$ which will be
denoted by $\cB.$ Clearly, $\mathbf{1}_\cA\in\cB,$ so that
$\mathbf{1}_\cA=\mathbf{1}_\cB$.

For a subset $X\subseteq G$ we denote by $\cA_X$ the linear subspace
$\bigoplus_{g\in X}\cA_g$ of $\cA.$ From (\ref{eq_Ggrad}) we
conclude that $\cA_H$ is a $*$-subalgebra of $\cA$ for the subgroup
$H$ of $G$.
%

\begin{prop}\label{prop_gradexp}
Let $p_H$ be the canonical projection of $\cA$ onto $\cA_H$, that
is, $p_H(a)=\sum_{g\in H}a_g$ for $a=\sum_{g\in G}a_g,$ where
$a_g\in\cA_g.$ Then $p_H$ is a conditional expectation of $\cA$ onto
$\cA_H.$
\end{prop}
\noindent\textbf{Proof.} Condition $(i)$ of Definition
\ref{defn_cond_exp} follows at once from (\ref{eq_Ggrad}). Our proof
is complete once we have shown that $p_H(\sum \cA^2)\subseteq\sum
\cA^2$.

We choose one element $k_i\in G,\ i\in G/H,$ in each left coset of
$H$ in $G.$ Let $a=\sum_{i\in G/H}b_i,$ where $b_i\in\cA_{k_iH}.$ If
$i,j\in G/H$, then $b_j^*b_i\in\cA_{Hk_j^{-1}k_iH},$ hence we have
$p_H(b_i^*b_i)=b_i^*b_i$ and $p_H(b_j^*b_i)=0$ if $i\neq j.$ Using
the latter facts we obtain
\begin{gather}\label{eq_p_H(a^*a)}
p_H(a^*a)=p_H(\sum_{i\in G/H}\sum_{j\in G/H}b_j^*b_i^{})=\sum_{i\in
G/H}b_i^*b_i^{}\in\sum\cA^2.
\end{gather}
\hfill $\Box$\mn

The the map $p_H$ from Proposition \ref{prop_gradexp} is called
the \textit{canonical conditional expectation} of the $G$-graded
$*$-algebra $\cA$ onto the $*$-subalgebra $\cA_H.$

Equation (\ref{eq_p_H(a^*a)}) shows that $p_H$ is faithful when
$\sum_{k=1}^n a_k^*a_k =0$ for arbitrary $a_1,\dots,a_n \in \cA$
implies that $a_1=\dots=a_n=0$. In particular, $p_H$ is faithful
when $\cA$ is an $O^*$-algebra.

Another immediate consequence of (\ref{eq_p_H(a^*a)}) is stated as
\begin{cor}\label{corconeb}
An element $a\in\cA$ belongs to the cone $\sum\cA^2\cap\cA_H$ if and
only if it can be presented as a finite sum of squares $\sum
b_i^*b_i^{},$ where each $b_i$ belongs to some $\cA_{gH},\ gH\in
G/H.$
\end{cor}

\begin{exam}\label{exam_weyl_cone}
Let $\cA=\langle a,a^*|aa^*-a^*a=1\rangle$ be the Weyl algebra (see
Example \ref{exam_weyl_alg}). Then $\cA$ is a $\dZ$-graded
$*$-algebra with $\dZ$-grading defined by $a\in\cA_1,\
a^*\in\cA_{-1}$ and we have $\cB=\dC[N]$, where $N=a^*a.$ We now use
Corollary \ref{corconeb} to describe the cone $\sum\cA^2\cap\cB.$

Suppose $k\in\dN$. Let $a_k\in\cA_k.$ Then $a_k$ is of the form
$a_k=a^kp_k$, where $ p_k\in\dC[N]$, and
$$a_k^*a_k^{}=p_k^*a^{k*}a^kp_k=N(N-1)\dots(N-k+1)p_k^*p_k^{}.$$
For $a_{-k}\in\cA_{-k}$ we have $a_{-k}=a^{*k}p_{-k}$, where
$p_{-k}\in\dC[N]$, and
$$a_{-k}^*a_{-k}^{}=p_{-k}^*a^{k}a^{*k}p_{-k}=(N+1)(N+2)\dots(N+k)p_{-k}^*p_{-k}^{}.$$
One easily verifies that $a_{-k}^*a_{-k}$ belongs to
$\sum\cB^2+N\sum\cB^2$. Hence from Corollary \ref{corconeb} we
obtain
\begin{gather}\label{eq_cone_in_weyl_alg}
\sum\cA^2\cap\cB=\sum\cB^2+N\sum\cB^2+N(N-1)\sum\cB^2+\dots
\end{gather}
This result was derived in \cite{fs} by other methods. Among others
it shows that $\sum\cA^2\cap\cB\neq \sum\cB^2$ and that the
canonical conditional expectations $p:\cA \to\cB$ is not strong.
\hfill $\edex$
\end{exam}

\begin{exam} Let $G$ be a discrete group and $H$ a normal
subgroup of $G$. Then the group algebra $\dC[G]$ becomes a
$G/H$-graded $*$-algebra in canonical manner. The canonical
conditional expectation coincides with the one from the Example
\ref{exam_cond_exp_CG_CH}, so by Proposition
\ref{prop_cond_exp_CG_CH} it is strong. In particular, we have
$\sum\dC[G]^2\cap\dC[H]=\sum\dC[H]^2.$ \hfill $\edex$
\end{exam}

\begin{exam}\label{exam_crossed_prod}
Let $A$ be a unital $*$-algebra. Let $G$ be a (discrete) group which acts as $*$-automorphism group $g\mapsto\alpha_g$ on $A.$ Recall that the crossed product $*$-algebra $\cA=A\times_{\alpha}G$ is defined as follows. As a linear space $\cA$ is the tensor product $A\otimes\dC[G]$ or equivalently the vector space of $A$-valued functions on $G$ with finite support. Product and involution on $\cA$ are determined by $(a\otimes g)(b\otimes h)=a\alpha_g(b)\otimes gh$ and $(a\otimes g)^*=\alpha_{g^{-1}}(a^*)\otimes g^{-1},$ respectively. If we identify $b$ with $b\otimes e$ and $g$ with $1 \otimes g$, then the $*$-algebra $A\times_{\alpha}G$ can be considered as the universal $*$-algebra generated by the two $*$-subalgebras $A$ and $\dC[G]$ with cross commutation relations $gb=\alpha_g(b)g$ for $b\in A$ and $g\in G$. Set $\cA_g:=A\otimes g$ for $g \in G$. Then $\cA$ becomes a $G$-graded $*$-algebra with canonical conditional expectation $p$ onto $\cB=\cA_e$ given by $p(a\otimes g)=\delta_{g,e}a\otimes e.$

\begin{prop}
The canonical conditional expectation $p:\cA\times_{\alpha}G\to\cB$
is strong.
\end{prop}
\noindent\textbf{Proof.} Let $x=\sum_{g\in G}a_g\otimes g,\ a_g\in A,$ be an element of the $\cA\times_{\alpha}G.$ Then
\begin{gather*}
p(xx^*)=p\left(\sum_{g\in G}\sum_{h\in G}(a_g\otimes g)(a_h\otimes
h)^*\right)=p\left(\sum_{g\in G}\sum_{h\in
G}a_g\alpha_{gh^{-1}}(a_h^*)\otimes gh^{-1}\right)=\\
=\sum_{g\in G}a_ga_g^*\otimes e=\sum_{g\in G}(a_g\otimes
e)(a_g\otimes e)^*\in\sum\cB^2.
\end{gather*}
\hfill $\Box\ \edex$
\end{exam}

\begin{exam}\label{exam_graded}
Let $G$ be a compact abelian group. Then the dual group $\widehat{G}$
is a discrete abelian group. We now establish a duality between
actions of $G$ and $\widehat{G}$-gradings on a $*$-algebra $\cA$
(cf. Example \ref{exam_loc_fin_comp_gr}).

Suppose that an action $\alpha:G\to Aut(\cA)$ is given. Assume, in
addition, that the action is locally finite-dimensional (see Example
\ref{exam_loc_fin_comp_gr}). For $\psi\in\widehat{G},\ \psi:G\to\dT$
put
\begin{gather}\label{eq_grading_from_action}
\cA_\psi=\set{a\in\cA|\ \alpha_g(a)=\psi(g)a,\ \mbox{for all}\ g\in
G}.
\end{gather}

If $\cA$ is a $\widehat{G}$-graded $*$-algebra, we define an action
of $\widehat{\widehat{G}}=G$ on $\cA$ as follows. For
$a=\sum_{\psi\in\widehat{G}}a_\psi,\ a_\psi\in \cA_\psi$ and $g\in
G,$ define a $*$-automorphism $\alpha_g$ by putting
\begin{gather}\label{eq_action_from_grading}
\alpha_g(a):=\sum_{\psi\in G}\psi(g)a_\psi.
\end{gather}

\begin{prop}
Equations (\ref{eq_grading_from_action}) and
(\ref{eq_action_from_grading}) give a one-to-one correspondence
between locally finite-dimensional actions of $G$ on $\cA$ and
$\widehat{G}$-gradings of $\cA.$
\end{prop}
\noindent\textbf{Proof.} Let $\alpha:G \to Aut(\cA)$ be locally
finite-dimensional action and let $\cA_\psi$ be defined by
(\ref{eq_grading_from_action}). We consider $\cA$ as $G$-module and
$\cA_\psi$ as unitary $G$-submodule. Take a finite-dimensional
$\alpha$-invariant linear subspace $V$ of $\cA$. Since $G$ is
compact, $V$ is unitarizable and hence spanned by its subspaces
$\cA_\psi.$ Since the action of $G$ is locally finite-dimensional,
$\cA$ is spanned by such subspaces $V$ and so by $\cA_\psi$,
$\psi\in\widehat{G}.$ It is easily checked that
$\cA=\oplus_{\psi\in\widehat{G}} ~\cA_\psi$ is a
$\widehat{G}$-grading of $\cA.$

Conversely, suppose $\cA$ is a $\widehat{G}$-graded $*$-algebra. It
is clear that (\ref{eq_action_from_grading}) defines an action of
$G$ on $\cA$. Each element $a\in \cA$ is of the form $a=\sum_{i=1}^k
a_{\psi_i}$, where $a_{\psi_i}\in\cA_{\psi_i}$ and the elements
$\psi_i\in \widehat{G}$ are pairwise distinct. The elements
$a_{\psi_i}$ span a finite-dimensional subspace of $\cA$ which is
obviously invariant under the action (\ref{eq_action_from_grading}).
Hence the action (\ref{eq_action_from_grading}) is locally
finite-dimensional. \hfill $\Box\ \edex$
\end{exam}

\noindent\textbf{Remark.} For the study of modules over a $G$-graded
ring $\cA=\oplus_{g\in G}\cA_g,$ it is usually assumed that for all
$g,h\in G$ the linear span of $\cA_g\cA_h$ is equal to $\cA_{gh}$,
see \cite{nastas},\cite{markus}. Likewise in \cite{fd} it is
supposed that this linear span is dense in $\cA_{gh}$.
We have not made such an assumption, because it is not satisfied in
most of our standard examples. For instance, if $\cA$ is the Weyl
algebra (Example \ref{exam_weyl_cone}), then we have $\cB=\dC[N]$,
$\cA_1=a\cB$ and $\cA_{-1}=a^*\cB=\cB a^*.$ Therefore, the linear
span of $\cA_{-1}\cdot\cA_1$ is equal to $N\cdot\dC[N]$ which is
different from $\cB.$

\section{Systems of imprimitivity}\label{sect_sys_impr}
Let $\cA=\oplus_{g \in G}~ \cA_g$ be a $G$-graded $*$-algebra. We
retain the notation of the previous section. Recall that for a
subgroup $H\subseteq G$, the left $G$-space of left $H$-cosets is
denoted by $G/H$.
\begin{defn}\label{defn_sys_impr}
Let $\pi$ be a $*$-representation of the $*$-algebra $\cA$ and let
$E$ be a mapping from the set $G/H$ to the set of projections of
the underlying Hilbert space $\cH_\pi$ such that
\begin{enumerate}
  \item[$(i)$] $E(t_1)E(t_2)=0$ for all $t_1,t_2\in G/H,\ t_1\neq t_2,$ and $\sum_{t\in G/H} E(t)=I,$
  \item[$(ii)$] $E(gt)\pi(a_g)\subseteq \pi(a_g)E(t)\ \mbox{for all}\ g\in G,\ t\in G/H,\ a_g\in\cA_g.$
\end{enumerate}
We call the pair $(\pi,E)$ a \textit{system of imprimitivity} for
the algebra $\cA$ over $G/H.$
\end{defn}
\mn
Let $(\pi,E)$ be a system of imprimitivity. Let $t \in G/H$ and set $\cD_t(\pi):=\Ran E(t)\cap\cD(\pi)$.
The conditions in Definition \ref{defn_sys_impr} imply that
$$E(t)\cD(\pi)\subseteq\cD(\pi),~~\pi(\cA_{g})\cD_t(\pi)\subseteq\cD_{gt}(\pi)~\mbox{for}~g\in G,~\mbox{and}~\cD(\pi)\subseteq\widetilde\oplus_{t\in G/H} \cD_t(\pi),$$
where $\widetilde\oplus$ denotes the direct Hilbert sum.

A system of imprimitivity $(\pi,E)$ is called \textit{non-degenerate} if for all $t\in G/H$ the subspace $\pi(\cA_{t})\cD_H(\pi)$ is dense in $\cD_t(\pi)$ with respect to the graph topology of $\pi.$ Otherwise, we say that $(\pi,E)$ is \textit{degenerate}.

\begin{lemma}\label{lemma_closure_impr_sys}
Let $H$ be a subgroup of $G$ and let $(\pi,E)$ be a system of
imprimitivity for the algebra $\cA$ over $G/H.$ Then the pair
$(\overline{\pi},E)$ is again a system of imprimitivity for $\cA$
over $G/H.$ Moreover, if $(\pi,E)$ is non-degenerate, then
$(\overline{\pi},E)$ is also non-degenerate.
\end{lemma}
\noindent\textbf{Proof.} From condition $(ii)$ we obtain
$\norm{\pi(a_g)E(t)\varphi}\leq\norm{\pi(a_g)\varphi}$ for
$a_g\in\cA_g$ and $\varphi \in \cD(\pi)$. This shows that $E(t)$ is
a continuous mapping of $\cD(\pi)$ with respect to the graph
topology of $\pi.$ Hence condition $(2)$ extend by continuity to the
closure $\overline{\pi}$ of $\pi.$ Obviously,
$(\overline{\pi},E)$ is non-degenerate if $(\pi,E)$ is. \hfill
$\Box$\mn

Systems of imprimitivity arise from induced representations in the
following way (see e.g. \cite{fd}, p.1248, for the case of finite
groups). Let $\rho$ be a non-zero inducible representation of the
algebra $\cA_H$ on a dense domain $\cD(\rho)$ of the Hilbert space
$\cH_\rho$ and let $\pi=\Ind_{\cA_H\uparrow\cA}\rho.$

Since $\cA=\bigoplus_{t\in G/H}\cA_t,$ we get
$$
\cA\otimes_{\cA_H}\cD(\rho)=\bigoplus_{t\in
G/H}\cA_t\otimes_{\cA_H}\cD(\rho).
$$
Recall that the representation space $\cH_\pi$ of $\pi$ is the
completion of the quotient space of the tensor product
$\cA\otimes_{\cA_H}\cD(\rho)$ by the kernel $\cK_\rho$ of the
sesquilinear form $\langle\cdot,\cdot\rangle_0$ defined by
(\ref{eq_welldef1}). Let $\cH_{t,0}$ denote the subspace of vectors
$\xi_t\in\cA_t\otimes_{\cA_H}\cD(\rho),\ t\in G/H,$ such that
$\langle \xi_t,\xi_t\rangle_0=0.$ Take $\eta=\sum_{t\in
G/H}\eta_t\in\cH_0$, where $\
\eta_t\in\cA_t\otimes_{\cA_H}\cD(\rho).$ Since
$\langle\eta_t,\eta_s\rangle_0=0$ for $t\neq s$ we get
$$
0=\langle\eta,\eta\rangle_0=\sum_{s,t\in
G/H}\langle\eta_s,\eta_t\rangle_0=\sum_{t\in
G/H}\langle\eta_t,\eta_t\rangle_0,
$$
that is, every $\eta_t$ belongs to $\cH_{t,0}.$ This implies that
$\cH_0=\bigoplus_{t\in G/H}\cH_{t,0} $ and hence
$$(\cA\otimes_{\cA_H}\cD(\rho))/\cH_0=\bigoplus_{t\in
G/H}(\cA_t\otimes_{\cA_H}\cD(\rho))/\cH_{t,0}.
$$
Note that for different left cosets $t\in G/H$ the subspaces
$(\cA_t\otimes_{\cA_H}\cD(\rho))/\cH_{t,0}$ are pairwise orthogonal.
For $t\in G/H,$ we denote by $E(t)$ the orthogonal projection from
$\cH_\pi$ onto the completion of the subspace
$(\cA_t\otimes_{\cA_H}\cH_\rho)/\cH_{t,0}.$

\begin{prop}\label{prop_impr_induced}
The pair $(\pi,E)$ constructed above is a non-degenerate system of imprimitivity for the algebra $\cA$ over $G/H.$
\end{prop}
\noindent\textbf{Proof.} Because of Lemma \ref{lemma_closure_impr_sys} it suffices to check the conditions in Definition \ref{defn_sys_impr} for the restriction of $\pi$ to its core $(\cA\otimes_{\cA_H}\cD(\rho))/\cH_0.$ One easily verifies condition $(i)$. We now show that condition $(ii)$ is satisfied. Since the vectors $[a_t\otimes v],\ a_t\in\cA_t,\ t\in G/H,\ v\in\cD(\rho),$ span a core for $\pi,$ it is enough to check $(ii)$ for vectors of this form. Let us fix elements $g\in G,\ a_g\in\cA_g,\ s,t\in G/H$ and $ v\in \cD(\rho).$ Then we have
$$
\pi(a_g)E(s)[a_t\otimes v]=\left\{
                                 \begin{array}{ll}
                                   [a_g a_t\otimes v], & \mbox{if}\ s=t; \\
                                   0, & \mbox{otherwise.}
                                 \end{array}
                               \right.
$$
Since the same result is obtained for $E(gs)\pi(a_g)[a_t\otimes
v]=E(gs)[a_g a_t\otimes v],$ $(ii)$ holds.

The equality $\pi(a_t)[\mathbf{1}_\cA\otimes v]=[a_t\otimes v]$
 implies that the span of $\pi(\cA_t)\cD_H(\pi)$ is equal to $\cD_t(\pi),$ so $(\pi,E)$
is non-degenerate. \hfill $\Box$\mn

We call the pair $(\pi,E)$ from Proposition \ref{prop_impr_induced}
the \textit{system of imprimitivity induced by $\rho.$}

\begin{thm}(First Imprimitivity Theorem)\label{thm_imprim_non_degenerate}
Let $\cA=\oplus_{g \in G}~ \cA_g$ be a $G$-graded $*$-algebra and
$H$ a subgroup of $G$. Suppose that $\pi$ is a closed
$*$-representation of $\cA$ and $(\pi,E)$ is a non-degenerate system
of imprimitivity for $\cA$ over $G/H$. Then there exists a unique,
up to unitary equivalence, closed $*$-representation $\rho$ of
$\cA_H$ such that
\begin{enumerate}
  \item[$(i)$] $\rho$ is inducible,
  \item[$(ii)$] $(\pi,E)$ is unitarily equivalent to the system of imprimitivity
induced by $\rho.$
\end{enumerate}
\end{thm}
\noindent\textbf{Proof.} By condition $(ii)$ in Definition
\ref{defn_sys_impr}, the projection $E(H)$ commutes with the
operators $\pi(a_H),\ a_H\in\cA_H.$ Hence the restriction of the
representation $\Res_{\cA_H}\pi$ to the subspace $\Ran E(H)$ is a
well-defined $*$-representation of the $*$-algebra $\cA_H$ denoted
by $\rho.$ The domain $\cD(\rho)$ is equal to $\Ran
E(H)\cap\cD(\pi)$ and the representation space $\cH_\rho$ is $\Ran
E(H).$

First we prove that $\rho$ is inducible. We have to show that the
 form $\langle\cdot,\cdot\rangle_0$ is nonnegative. Take a vector $\xi=\sum_r a_r\otimes
v_r\in\cA\otimes_{\cA_H}\cD(\rho)$, where $ v_r\in\cD(\rho),\
a_r\in\cA.$ Each $a_r$ can be presented as a finite sum
$a_r=\sum_{t\in G/H}a_{r,t}$, where $a_{r,t}\in\cA_t,\ t\in G/H.$
Then we have
\begin{gather}
\nonumber\langle\xi,\xi\rangle_0=\langle\sum_r a_r\otimes v_r,\sum_s a_s\otimes v_s\rangle_0=\sum_{r,s}\langle\rho(p_H(a_s^*a_r^{}))v_r,v_s\rangle=\\
\label{eq_pos_semi_def_ind_impr_sys}=\sum_{r,s}\langle\rho(\sum_{t\in G/H} a_{s,t}^*a_{r,t}^{})v_r,v_s\rangle=\sum_{t\in G/H}\sum_{r,s}\langle\rho(a_{s,t}^*a_{r,t}^{})v_r,v_s\rangle=\\
\nonumber=\sum_{t\in G/H}\sum_{r,s}\langle\pi(a_{s,t})v_r,\pi(a_{r,t})v_s\rangle=\sum_{t\in G/H}\langle\sum_r\pi(a_{r,t})v_r,\sum_s\pi(a_{s,t})v_s\rangle\geq 0.%
\end{gather}
This shows that $\rho$ is inducible.

Let $(\pi_1,E_1)$ denote the system of imprimitivity on the space $\cH_{\pi_1}$ induced by $\rho.$ We have to prove that $(\pi_1,E_1)$ is unitarily equivalent to $(\pi,E).$ Define a linear mapping $F_0:\cA\otimes\cD(\rho)\to\cD(\pi)$ by putting $F_0(a\otimes v)=\pi(a)v$, where $v\in\cD(\rho)\subseteq\cD(\pi),\ a\in\cA.$ It is clear that $F_0$ maps $\cA\otimes_{\cA_H}\cD(\rho)$ into $\cD(\pi).$ Recall that $\cK_\rho$ denotes the kernel of the sesqulinear form $\langle\cdot,\cdot\rangle_0.$ Reasoning in the same manner as in (\ref{eq_pos_semi_def_ind_impr_sys}) it follows that for any $\xi\in\cA\otimes_{\cA_H}\cD(\rho)$ we have $\langle\xi,\xi\rangle_0=\langle F_0(\xi),F_0(\xi)\rangle$. Therefore, the quotient mapping from $\cA\otimes_{\cA_H}\cD(\rho)/\cK_\rho$ to $\cH_\pi$ is a well-defined isometric linear mapping. We extend this mapping by continuity to an isometry $F:\cH_{\pi_1}\to\cH_\pi.$

We claim that $F$ intertwines the systems $(\pi,E)$ and $(\pi_1,E_1).$ Take $k\in G,\ a_k\in\cA_k,\ a_t\in\cA_t,\ t\in G/H,\ v\in\cD(\rho).$ Then we obtain
\begin{gather*}
F(\pi_1(a_k)([a_t\otimes v]))=F(a_ka_t\otimes v)=\pi(a_ka_t)v=\\
\pi(a_k)\pi(a_t)v=\pi(a_k)F([a_t\otimes v])
\end{gather*}
which means that $F$ intertwines $\pi$ and $\pi_1.$

For $v\in\cD(\rho),\ a_t\in\cA_t,\ t\in G/H$ condition $(ii)$ in Definition
\ref{defn_sys_impr} implies that $\pi(a_t)v\in\cD_t(\pi).$ The subspace $\cD_t(\pi_1),$ is spanned by the vectors $[a_t\otimes v],\ a_t\in\cA_t,\ v\in\cD(\rho),$ and we have $F([a_t\otimes v])=\pi(a_t)v\in\cD_t(\pi).$ Thus, $F(\cD_t(\pi_1))\subseteq\cD_t(\pi)$ and $F$ intertwines $E$ and $E_1.$

Since $(\pi,E)$ is non-degenerate, the vectors $F([a_t\otimes v])=\pi(a_t)v,\ a_t\in\cA_t,\ v\in\cD(\rho),$ span a dense linear subspace $\cD_t(\pi_1)$ of $\cD_t(\pi)$ in the graph topology of $\pi.$ In particular, we have $F(\Ran E_1(t))=\Ran E(t),$ so that $F$ is a unitary operator. Since the graph topology on $F(\cD_t(\pi_1))$ is the same as that of $\pi$ and $\pi_1$ is
closed by definition, we have $F(\cD_t(\pi_1))=\cD_t(\pi)$ for each $t\in G/H,$ which implies that $F(\cD(\pi_1))=\cD(\pi).$ That is, $\pi$ and $\pi_1$ are unitarily equivalent.

Let $\rho_1$ be an inducible closed $*$-representation of $\cA_H$ on the Hilbert space $\cH_{\rho_1}$ and let $(\pi_2,E_2)$ be the system of imprimitivity for $\cA$ over $G/H$ induced by $\rho_1.$ It follows from the previous considerations that $\rho_2:=\Res_{\cA_H}\pi_2\upharpoonright\Ran E_2(H)$ is well-defined $*$-representation of $\cA_H.$ One immediately verifies
that the canonical isomorphism $v\leftrightarrow[\mathbf{1}_\cA\otimes v]$ of $\cH_{\rho_1}$ and $\Ran E_2(H)$ defines a unitary equivalence of $\rho_1$ and $\rho_2.$ \hfill $\Box$\mn

Summarizing, we have shown that there is a one-to-one correspondence
between unitary equivalence classes of inducible representations of
$\cA_H$ and unitary equivalence classes of non-degenerate closed
systems of imprimitivity for $\cA$ over $G/H$. In particular, the
inducing representation $\rho$ is determined uniquely up to unitary
equivalence by the system of imprimitivity.\mn

The following example shows that the non-degeneracy assumption of
the system of imprimitivity is crucial in Theorem
\ref{thm_imprim_non_degenerate}.

\begin{exam} Let $\cA_q$ be the $*$-algebra $\dC\langle
a,a^*|aa^*-qa^*a=1\rangle,$ where $q>-1.$ Put $\lambda_0=0$ and
$\lambda_k=\sqrt{1+q+q^2+\dots+q^{k-1}},\ k\in\dN.$ Let $\cH$ be a Hilbert space with orthonormal base
$\set{e_k,\ k\in\dN_0}.$ There is a $\ast$-representations $\pi$ of $\cA_q$
on $\cD(\pi)=\rm{Lin} \set{e_k ; k \in \dN_0}$ such that
$$
\pi(a)e_k=\lambda_k e_{k-1},\pi(a^*)e_k=\lambda_{k+1}e_{k+1},\
\mbox{for}\ k\in\dN_0,
$$
where $e_{-1}:=0$. The representation $\pi$ is bounded if and only if $-1\leq q\leq 0.$ Note that in the
case $q=1$ the algebra $\cA_q$ is just the Weyl algebra and $\overline{\pi}$
is the Fock-Bargmann representation.

Let $E(n),\ n\in\dN,$ be the orthogonal projection onto
$\mathbb{C}{\cdot}e_{n-1}$ and put $E(n):=0$ for $n\leq 0.$ Then the
pair $(\pi,E)$ is a system of imprimitivity for $\cA$ over
$G=\mathbb{Z}.$ Since $E(0)=0$, it follows immediately from the
construction of the induced system of imprimitivity that $(\pi,E)$
is not induced by a $*$-representation of $\cB$.\hfill $\edex$
\end{exam}

We now define another construction of systems of imprimitivity. It
will also include the system of imprimitivity in the latter example.
Fix a system of imprimitivity $(\pi,E)$ for $\cA$ over $G/H$ and an
element $f\in G.$ Define a mapping $E^f$ from the set $G/fHf^{-1}$
into the set of projections on the space $\cH_\pi$ by
$E^f(k(fHf^{-1})):=E(kfH),\ k\in G.$

\begin{prop}
The pair $(\pi,E^f)$ constructed above is a well-defined system of
imprimitivity for $\cA$ over $G/fHf^{-1}.$
\end{prop}
\noindent\textbf{Proof.} Take $k_1(fHf^{-1}),k_2(fHf^{-1})\in
G/fHf^{-1}$, where $k_1,k_2\in G.$ The cosets $k_1(fHf^{-1})$ and
$k_2(fHf^{-1})$ are equal if and only if $k_2^{-1}k_1^{}\in
fHf^{-1}$ which is equivalent to $k_1fH=k_2fH.$ This implies that
$E^f$ is well-defined. It is straightforward to verify that
$(\pi,E^f)$ satisfies the two conditions in Definition
\ref{defn_sys_impr}. \hfill $\Box$

\begin{defn}
If $(\pi,E),\ f\in G,\ (\pi,E^f)$ are as above, we say that the
system $(\pi,E^f)$ is \textit{conjugated} to the system $(\pi,E)$
by the element $f\in G.$
\end{defn}

Our second Imprimitivity Theorem describes systems of imprimitivity
which are not necessarily non-degenerate. We prove it now for
bounded representations (cf. also the Imprimitivity Theorem in
\cite{fd}, p.1192). In Section \ref{sect_well_beh_sys_impr} we
formulate its analogue for well-behaved systems of imprimitivity
(Theorem \ref{thm_impr_for_well_beh_degen_system}).

The following definition and the subsequent lemma are used in the proof of Theorem \ref{thm_impr_for_bounded_degen_system} below.
\begin{defn}\label{defn_gener_sys_impr}
Let $(\pi,E)$ be a system of imprimitivity for $\cA$ over $G/H$ and let $fH\in G/H$. We say that $(\pi,E)$ is \textit{generated by the projection} $E(fH)$ if for every $gH\in G/H$ the linear subspace $\pi(\cA_{gHf^{-1}})(\cD_{fH}(\pi))$ is dense in $\cD_{gH}(\pi)$ with respect to the graph topology of $\pi.$
\end{defn}
\begin{lemma}\label{lemma_non_deg_iff_gen_by}
A system of imprimitivity $(\pi,E)$ is generated by the projection
$E(fH),\ fH\in G/H,\ f\in G,$ if and only if the conjugated system
of imprimitivity $(\pi,E^f)$ over $G/fHf^{-1}$ is non-degenerate.
\end{lemma}

The simple proof of Lemma \ref{lemma_non_deg_iff_gen_by} will be
omitted. The next theorem says that for bounded representations each
system of imprimitivity over $G/H$ can be obtained as a direct sum
of conjugated systems by elements of $G$.

\begin{thm}(Second Imprimitivity Theorem)\label{thm_impr_for_bounded_degen_system}
Let $\cA=\oplus_{g \in G}~ \cA_g$ a $G$-graded $*$-algebra, $H$ a
subgroup of $G$ and $(\pi,E)$ a system of imprimitivity for $\cA$
over $G/H.$ Suppose the $*$-representation $\pi$ acts by bounded
operators on $\cD(\pi)=\cH_\pi$. We fix one element $k_t\in G,\ t\in
G/H,$ in each left coset from $G/H.$ Then for every $t\in G/H$ there
exists a bounded $*$-representation $\rho_t$ of
$\cA_{k_t^{}Hk_t^{-1}}$ on a Hilbert space $\cH_t$ such that:
\begin{enumerate}
  \item[$(i)$] $\rho_t$ is inducible,
  \item[$(ii)$] $(\pi,E)$ is the direct sum of systems of imprimitivity $(\pi_t,E_t),\ t\in G/H,$ where
  $(\pi_t,E_t)$ is conjugated by the element $k_t$ to the system of imprimitivity induced by $\rho_t,\ t\in G/H.$
\end{enumerate}
\end{thm}
\noindent\textbf{Proof.} Let $(\pi_1,E_1)$ be an subsystem of
imprimitivity of $(\pi,E)$ over $G/H,$ that is, $\pi_1\subseteq\pi$
is a subrepresentation of $\pi$ on a Hilbert subspace
$\cH_1\subseteq\cH_\pi$ and for all $gH\in G/H$ we have $\Ran
E_1(gH)\subseteq\Ran E(gH).$ Since $\pi$ is a bounded
$*$-representation, there is a $*$-representation $\pi_2$ on
$\cH_2:=\cH_\pi \ominus\cH_1$ such that $\pi =\pi_1\oplus \pi_2$.
Put $E_2(gH):=E(gH)\ominus\Ran E_1(gH)$ for $gH\in G/H$. Then
$(\pi_2,E_2)$ is again a system of imprimitivity for $\cA$ over
$G/H.$ Indeed, condition $(i)$ in Definition \ref{defn_sys_impr} is
obvious and condition $(ii)$ follows immediately by subtracting the
equation $\pi_1(a_g)E_1(fH)=E_1(gfH)\pi_1(a_g)$ from
$\pi(a_g)E(fH)=E(gfH)\pi(a_g)$, where $g\in G,a_g\in\cA_g,\ fH\in
G/H.$ That is, we have shown that every subsystem of imprimitivity
has a complement.

Now we fix $fH\in G/H.$ Let $E_1(gH)$ denote the orthogonal
projection onto the closure of $\Ran\pi(\cA_{gHf^{-1}})E(fH)$ and
set $\cH_1:=\oplus_{t\in G/H}\Ran E_1(t).$ It is easily checked that
the family of projections $E_1(t),\ t\in G/H,$ satisfies condition
$(i)$ of Definition \ref{defn_sys_impr}. Let $g\in G, a_g\in\cA_g$
and $kH\in G/H.$ Then we have
\begin{gather*}
\pi(a_g)\Ran
E_1(kH)=\pi(a_g)\overline{\Ran\pi(\cA_{kHf^{-1}})E(fH)}\subseteq\overline{\Ran\pi(\cA_{gkHf^{-1}})E(fH)}=E_1(gkH),
\end{gather*}
which shows that the subspace $\cH_1$ is invariant under all
operators $\pi(a),\ a\in\cA.$ If we denote by $\pi_1$ the
restriction of $\pi$ to $\cH_1,$ then condition $(ii)$ in Definition
\ref{defn_sys_impr} holds for the pair $(\pi_1,E_1).$ Therefore,
$(\pi_1,E_1)$ is an subsystem of imprimitivity for $\cA$ over $G/H.$
The system $(\pi_1,E_1)$ is generated by $E_1(fH)=E(fH).$

Combining the considerations of the preceding paragraphs with Zorn's
lemma we conclude that there exist system of imprimitivity
$(\pi_t,E_t),\ t\in G/H,$ for $\cA$ over $G/H$ such that every
$(\pi_t,E_t)$ is generated by the projection $E_t(k_tH),\ t\in G/H,$
and $(\pi,E)$ is equal to the orthogonal direct sum of $(\pi_t,E_t),
t\in G/H.$

Lemma \ref{lemma_non_deg_iff_gen_by} together with Theorem
\ref{thm_imprim_non_degenerate} imply that each conjugated system
$(\pi_t,E_t^{k_t}),\ t\in G/H,$ is induced by some representation
$\rho_t$ of the $*$-algebra $\cA_{k_t^{}Hk_t^{-1}}.$ By the
construction of $\rho_t$ (see the proof of the Theorem
\ref{thm_imprim_non_degenerate}), $\rho_t$ it is a bounded
$*$-representation. \hfill $\Box$\mn

\noindent\textbf{Remark.} We do not know a generalization of Theorem
\ref{thm_impr_for_bounded_degen_system} for \textit{general
unbounded} representations. The main difficulty lies in the fact
that for a closed subrepresentation $\pi_1$ of closed
$*$-representation $\pi$ in general there is no representation
$\pi_2$ such that $\pi=\pi_1\oplus\pi_2.$

\section{A partial group action defined by the grading}\label{sect_grad_alg_com}
Throughout this section we assume that $\cA=\bigoplus_{g\in G}\cA_g$
is a $G$-graded unital $*$-algebra and that the
\textit{$*$-subalgebra $\cB:=\cA_e$ is commutative}. The canonical
conditional expectation of $\cA$ onto $\cB$ is denoted by $p.$

Let $\cBd$ be the set of all characters of $\cB$, that is, $\cBd$ is
the set of nontrivial $*$-homomorphisms $\chi:\cB \to \dC$. The set
of characters from $\cBd$ which are nonnegative on the cone
$\sum\cA^2\cap\cB$ is denoted by $\cBp.$

In addition we assume in this section that all characters $\chi\in\cBp$ satisfy the
following condition:
\begin{gather}\label{eq_cond_compatible}
\chi(c^*d)\chi(d^*c)=\chi(c^*c)\chi(d^*d)\ \mbox{for all}\
\chi\in\cBp,\ g\in G,\ \mbox{and}\ c,d\in\cA_g.
\end{gather}
Note that for $c,d\in\cA_g$ we have $c^*d,\ d^*c,\ c^*c,\ d^*d\in
\cA_{g^{-1}}\cdot\cA_g\subseteq\cA_e=\cB.$ Hence all expressions
in the equation (\ref{eq_cond_compatible}) are well-defined.

\begin{prop}\label{prop_crossed_prod_is_compat}
Let $\cA$ denote the crossed product algebra $A\times_{\alpha}G$
from Example \ref{exam_crossed_prod}. Assume that $A$ is
commutative, so that $\cB=A\otimes e$ is commutative. Then condition
(\ref{eq_cond_compatible}) is satisfied.
\end{prop}

Proposition \ref{prop_crossed_prod_is_compat} follows at once from
the more general

\begin{prop}\label{prop_full=>compatible} Assume that for every
$g\in G$ there exists an element $a_g\in \cA_g$ such that
$\cA_g=a_g\cB$ or $\cA_g=\cB a_g.$ Then condition
(\ref{eq_cond_compatible}) is satisfied.
\end{prop}
\noindent\textbf{Proof.} Fix a $g\in G.$ Assume that there exist an
element $a_g\in\cA_g$ such that $\cA_g=a_g\cB.$ Take $\chi\in\cBp$
and $c,d\in\cA_g.$ Then there exist $c_1,d_1\in\cB$ such that
$c=a_gc_1$ and $d=a_gd_1$. We now compute
$$\chi(c^*d)\chi(d^*c)=\chi(c_1^*)\chi(a_g^*a_g^{})\chi(d_1)\chi(d_1^*)\chi(a_g^*a_g^{})\chi(c_1)=\chi(c^*c)\chi(d^*d).$$
In the same way one proves (\ref{eq_cond_compatible}) in the case when
$\cA_g=\cB a_g,\ a_g\in\cA_g.$ \hfill $\Box$\mn

The main content of this section is the following partial action
of $G$ on the set $\cBp.$

\begin{defn}\label{defn_action}
Let $\chi\in\cBp$ and $g\in G.$ We say that $\chi^g$ is defined if
there exists an element $a_g\in \cA_g$ such that
$\chi(a_g^*a_g^{})\neq 0.$ In this case we set
\begin{gather}\label{eq_action_partial}
\chi^g(b):=\frac{\chi(a_g^*ba_g^{})}{\chi(a_g^*a_g^{})}\
\mbox{for}\ b\in\cB.
\end{gather}
For $g\in G$ we denote by $\cD_g$ the set of all
characters $\chi\in\cBp$ such that $\chi^g$ is defined.
\end{defn}

\noindent\textbf{Remarks.} 1. One could also define
$\chi^g$ as it was done in \cite{fd}. As noted in \cite{fd}, the
space $\cA_g,\ g\in G,$ has a natural structure of a $\cB$-rigged
$\cB{-}\cB$-bimodule, where $\cB$ acts by the
multiplication and the $\cB$-valued product is
$$[\cdot,\cdot]:\cA_g\times\cA_g\to\cB,\ [c,d]:=d^*c,\ c,d\in\cA_g.$$
Then $\chi^g$ is defined as the representation of $\cB$ induced from
$\chi$ via $\cA_g.$ Condition (\ref{eq_cond_compatible}) ensures
that $\chi^g$ is again a character.\mn

\noindent 2. Crossed-products defined by partial group actions on
$C^*$-algebras appeared in \cite{exel}. Our $G$-graded $*$-algebra
$\cA$ can be considered as another generalization of
crossed-product algebras. We shall not elaborate the details here.

\begin{prop}\label{prop_partial_action}
The map $\chi\mapsto\chi^g$ is a well-defined partial action of $G$
on the set $\cBp,$ that is:
\begin{enumerate}
\item[$(i)$] $\chi^g(b)$ in (\ref{eq_action_partial}) does not depend on
the choice of $a_g$ and we have $\chi^g \in\cBp,$

\item[$(ii)$] if $\chi^g$ and $(\chi^g)^h$ are defined, then $\chi^{hg}$
is defined and equal to $(\chi^g)^h,$

\item[$(iii)$] if $\chi^g$ is defined,
then $(\chi^g)^{g^{-1}}$ is defined and equal to $\chi,$

\item[$(iv)$] $\chi^e$ is defined and equal to $\chi.$
\end{enumerate}
\end{prop}
\noindent\textbf{Proof.} $(i)$: Let $\chi\in\cBp,\ g\in G,$ and $
c,d\in\cA_g$ such that $\chi(d^*d)\neq 0$ and $\chi(c^*c)\neq 0.$
Since $\cB$ is commutative, we have $bcd^*=cd^*b$ for $b\in\cB$.
Therefore we obtain
$$\chi(c^*bc)\chi(d^*d)=\chi(c^*bcd^*d)=\chi(c^*cd^*bd)=\chi(c^*c)\chi(d^*bd),
$$
so that
$$\frac{\chi(c^*bc)}{\chi(c^*c)}=\frac{\chi(d^*bd)}{\chi(d^*d)}.$$

We show that $\chi^g$ is again a character belonging to $\cBp.$ Let
$b_1,b_2\in\cB.$ Since $\cB$ is commutative, we have
$a_g^{}a_g^*b_1^{}=b_1^{}a_g^{}a_g^*.$ Hence we get
$$
\chi^g(b_1b_2)=\frac{\chi(a_g^*b_1^{}b_2^{}a_g^{})}{\chi(a_g^*a_g^{})}=
\frac{\chi(a_g^*a_g^{}a_g^*b_1^{}b_2^{}a_g^{})}{\chi(a_g^*a_g^{})\chi(a_g^*a_g^{})}=
\frac{\chi(a_g^*b_1^{}a_g^{}a_g^*b_2^{}a_g^{})}{\chi(a_g^*a_g^{})\chi(a_g^*a_g^{})}=
\chi^g(b_1)\chi^g(b_2).
$$

Next we prove the positivity of $\chi^g$. For take
$b\in\sum\cA^2.$ Since $\chi(\sum\cA^2)\geq 0$ and
$a_g^*ba_g\in\sum\cA^2$ we have $\chi^g(b)>0.$

$(ii)$: Let $\chi\in\cBp$ and $g,h\in G$ such that $(\chi^g)^h$ is
defined. Then there exists $a_g\in \cA_g$ such that
$\chi(a_g^*a_g)\neq 0.$ Since $(\chi^g)^h$ is defined, there exists
$a_h\in\cA_h$ such that
$$\chi^g(a_h^*a_h)=\frac{\chi(a_g^*a_h^*a_ha_g)}{\chi(a_g^*a_g)}\neq 0,$$
that is, $\chi((a_ha_g)^*a_ha_g)\neq 0.$ Since
$a_ha_g\in\cA_{hg},\ \chi^{hg}$ is well-defined. It is
straightforward to check that $(\chi^g)^h=\chi^{hg}.$

$(iii)$: Assume that $\chi^g$ is defined. Then there exists $a_g\in
\cA_g$ such that $\chi(a_g^*a_g)\neq 0.$ We have
$a_g^*\in\cA_{g^{-1}}$ and
$$
\chi^g(a_ga_g^*)=\frac{\chi(a_g^*a_ga_g^*a_g)}{\chi(a_g^*a_g)}=\chi(a_g^*a_g)\neq 0.
$$
Hence $(\chi^g)^{g^{-1}}$ is defined. One easily verifies that
$(\chi^g)^{g^{-1}}=\chi.$

$(iv)$ is trivial. \hfill $\Box$\mn

\noindent\textbf{Remark.} It follows from Proposition
\ref{prop_partial_action} that for each $g \in G$ the mapping $\chi\mapsto\chi^g$ defines a bijection $\alpha_g:\cD_g\to\cD_{g^{-1}}$
such that:
\begin{enumerate}
  \item[$(i)$] $\cD_e=\cBp$ and $\alpha_e$ is the identity mapping of $\cBp,$
  \item[$(ii)$] $\alpha_g(\cD_g\cap\cD_h)=\cD_{g^{-1}}\cap\cD_{hg^{-1}},$
  \item[$(iii)$] $\alpha_g(\alpha_h(x))=\alpha_{gh}(x),$ for $x\in\cD_g\cap\cD_{gh}.$
\end{enumerate}

In what follows, we shall use both notations $\alpha_g(\chi)$ and
$\chi^g$ for the partial action of $g\in G$ on $\chi\in\cBp$ and we
freely use the properties $(i)-(iii).$

It should be emphasized that up to now condition
(\ref{eq_cond_compatible}) has not been used for the partial action.
For the next proposition assumption (\ref{eq_cond_compatible}) is needed.
\begin{prop}\label{prop_defn_part_act_gener}
Let $a_g,c_g\in\cA_g,\ g\in G,$ and $\chi\in\cBp$ be such that
$\chi(a_g^*c_g^{})\neq 0.$ Then we have $\chi\in\cD_g$ and
\begin{gather}\label{eq_other_defn_part_act}
\chi^g(b)=\frac{\chi(a_g^*bc_g^{})}{\chi(a_g^*c_g^{})}\ \mbox{for
all}\ b\in\cB.
\end{gather}
\end{prop}
\noindent\textbf{Proof.} Since $\chi(a_g^*c_g^{})\neq 0$, we have
$\chi(c_g^*a_g^{})=\overline{\chi(a_g^*c_g^{})}\neq 0,$ so that
(\ref{eq_cond_compatible}) implies $\chi(a_g^*a_g^{})\neq 0,$ i.e.
$\chi\in\cD_g.$ Now (\ref{eq_other_defn_part_act}) follows from the
equality
$$
\chi(a_g^*ba_g^{})\chi(a_g^*c_g^{})=\chi(a_g^*a_g^{}a_g^*bc_g^{})=\chi(a_g^*a_g^{})\chi(a_g^*bc_g^{}).
$$
\hfill $\Box$\mn

Examples developed below show that in general $\chi^g$ is not
always defined, so that in general $\chi\mapsto\chi^g$ is not a
group action.

We introduce some more notation which will be kept till the
end of the paper.
For a fixed $\chi\in\cBp$ let
$$G_\chi=\set{g\in G|\chi^g\ \mbox{is defined}}.$$
We denote by $\Orb_\chi\subseteq\cBp$ the \textit{orbit} of the
$\chi,$ that is, $$\Orb_\chi=\set{\chi^g|\chi^g\ \mbox{is
defined}}.$$ Further, let $\St\chi\subseteq G_\chi$ denote the
\textit{stabilizer} of the element $\chi,$ that is,
$$\St\chi=\set{g\in G|\chi^g\ \mbox{is
defined and equal to}\ \chi}.
$$
A number of elementary properties of the partial action
of $G$ are collected in the following
\begin{prop}Let $\chi\in\cBp.$ Then we have:
\begin{enumerate}
  \item[$(i)$] $\St\chi$ is a subgroup of $G,$
  \item[$(ii)$] The union of sets $G_\psi,\ \psi\in\Orb\chi$ equipped with the multiplication derived from $G$ is a groupoid with identity element,
  \item[$(iii)$] if $\psi\in\cBp,$ then $\psi\in\Orb\chi$ if and only if $\Orb_\psi=\Orb_\chi,$
  \item[$(iv)$] if $\psi\in\Orb_\chi$, then $\St\chi$ and $\St\psi$ are
  conjugate subgroups of $G.$
\end{enumerate}
\end{prop}

Now we illustrate these concepts by a few examples.
\begin{exam} Let $A$ be a commutative $*$-algebra and
$\cA=A\times_{\alpha}G$ be the crossed-product algebra from Example \ref{exam_crossed_prod}. It was shown therein that $\sum\cA^2\cap\cB=\sum\cB^2.$ This implies that $\cBp=\cBd=\widehat{A}$ and the partial action defined by (\ref{eq_action_partial}) coincides with the usual group action of
$G$ on $\widehat{A}$ induced by the action of $G$ on $A.$ \hfill $\edex$
\end{exam}

\begin{exam}\label{exam_weyl_act} Let $\cA$ be the Weyl algebra. We retain
the notation from Examples \ref{exam_weyl_alg} and
\ref{exam_weyl_cone}. It follows from (\ref{eq_cone_in_weyl_alg})
that a character $\chi\in\cBd$ is non-negative on the cone
$\sum\cA^2\cap\cB$ if and only if $\chi(N)\in\dN_0.$ For
$k\in\dN_0,$ let $\chi_k$ denote the character of $\cBp$ defined by
$\chi_k(N)=k.$

Suppose that $n\in \dN_0$. Clearly, any element of the $\cA_n$ has
the form $a^np(N)$, where $ p\in\dC[N]$, and
$\chi^{}_k((a^np(N))^*a^np(N))\neq 0$ implies
$\chi^{}_k(a^{*n}a^n)\neq 0.$ So we obtain that
$$(\alpha_n(\chi_k))(N)=\frac{\chi^{}_k(a^{*n}Na^n)}{\chi^{}_k(a^{*n}a^n)}=
\frac{\chi_k(N(N-1)\dots(N-n+1)(N-n))}{\chi_k(N(N-1)\dots(N-n+1))}
$$
is defined if and only if $k\geq n$ and
$(\alpha_n(\chi_k))(N)=\chi_{k-n}(N).$

Analogously we conclude that
$$
(\alpha_{-n}(\chi_k))(N)=\frac{\chi^{}_k(a^nNa^{*n})}{\chi^{}_k(a^na^{*n})}=
\frac{\chi^{}_k((N+1)(N+2)\dots(N+n)^2)}{\chi^{}_k((N+1)(N+2)\dots(N+n))}
$$
is defined for all $n\in\dN$ and
$(\alpha_{-n}(\chi_k))(N)=\chi^{}_{k+n}(N),$ i.e.
$\alpha_{-n}(\chi_k)=\chi_{k+n}.$

The partial action is transitive, so $\cBp$ consists of a single
orbit. The stabilizer $\St{\chi_k}$ of each character $\chi_k$ is
trivial, the set $G_{\chi_k}$ is equal to
$\set{n\in\dZ|n\leq k}$. \hfill $\edex$
\end{exam}

The next proposition gives explicit formulas for representations
induced from characters. Recall that a character $\chi\in\cBp$ is a
one-dimensional $*$-representation of $\cB$ on the space $\dC$ and
the representation space $\cH_\pi$ of $\pi=\Ind\chi$ is spanned by
the vectors $[a\otimes 1]$, $a\in\cA$ (see Section
\ref{sect_rigged_induced}).

\begin{prop}\label{prop_space_ind_rep}
Let $\chi\in\cBp$ and $\pi=\Ind\chi$. Fix elements $a_g\in \cA_g,\ g\in
G_\chi,$ such that $\chi(a_g^*a_g^{})\neq 0,\ g\in G_\chi.$ Then we have:
\begin{enumerate}
  \item[$(i)$] The vectors
  $$e_g=\frac{[a_g\otimes 1]}{\sqrt{\chi(a_g^*a_g^{})}},\ g\in G_\chi,$$
  form an orthonormal base of the
  representation space $\cH_\pi$ of $\Ind\chi.$
  \item[$(ii)$] For $b_h\in\cA_h$ and $ h\in G$ we have
  $$\pi(b_h)e_g=\frac{\chi(a_{hg}^*b_h^{}a_g^{})}{\sqrt{\chi(a_{hg}^*a_{hg}^{})\chi(a_g^*a_g^{})}}e_{hg},\ \mbox{if}\ hg\in G_\chi$$
  and $\pi(b_h)e_g=0$ otherwise. In particular, if $b\in\cB$, then we have
  $$
  \pi(b)e_g=\frac{\chi(a_g^*ba_g^{})}{\chi(a_g^*a_g^{})}e_g=\chi^g(b)e_g.
  $$
\end{enumerate}
\end{prop}
\noindent\textbf{Proof.} First suppose that $b_g\in\cA_g$ and $g\notin G_\chi.$
Then $\norm{[b_g\otimes 1]}^2=\chi(b_g^*b_g^{})=0,$ so $\cH_\pi$ is
spanned by the vectors $[b_g\otimes 1]$, where $ b_g\in\cA_g$ and $ g\in
G_\chi.$

For $b_g\in\cA_g$ and $g\in G$ the equality
(\ref{eq_cond_compatible}) applied to $a_g$ and $b_g$ is
equivalent to the equation
$$|\langle[a_g\otimes 1],[b_g\otimes 1]\rangle|^2=\norm{[a_g\otimes
1]}^2\norm{[b_g\otimes 1]}^2,
$$
that is, we have equality in the Cauchy-Schwartz inequality. This
implies that $[a_g\otimes 1]=\lambda[b_g\otimes 1]$ for some
complex number $\lambda.$ Hence it follows that the elements
$[a_g\otimes 1], g\in G_\chi,$ span the space $\cH_\pi.$ Since
$\langle[a_g\otimes 1],[a_h\otimes
1]\rangle=\chi(p(a_h^*a_g))=\chi(0)=0$ for $g\neq h,$ the elements
$[a_g\otimes 1]$ are pairwise orthogonal. The square of the norm
of $[a_g\otimes 1]$ is equal to $\langle[a_g\otimes 1],[a_g\otimes
1]\rangle=\chi(a_g^*a_g^{}).$ Thus we have shown that the elements
$e_g,\ g\in G_\chi,$ form an orthonormal base of $\cH_\pi.$

Now let $b_h\in\cA_h,\ h\in H.$ If $hg\in G_\chi$ we have
\begin{gather*}
\pi(b_h)e_g=\frac{[b_ha_g\otimes
1]}{\sqrt{\chi(a_g^*a_g^{})}}=\frac{\lambda[a_{hg}\otimes
1]}{\sqrt{\chi(a_g^*a_g^{})}}=\lambda\frac{\sqrt{\chi(a_{hg}^*a_{hg}^{})}}{\sqrt{\chi(a_g^*a_g^{})}}e_{hg},
\end{gather*}
where $\lambda$ is equal to
$$\frac{\langle[b_ha_g\otimes 1],[a_{hg}\otimes 1]\rangle}{\langle[a_{hg}\otimes 1][a_{hg}\otimes 1]\rangle}=
\frac{\chi(a_{hg}^*ba_g^{})}{\chi(a_{hg}^*a_{hg})}.$$ This yields
the second statement of the theorem. \hfill $\Box$\mn

In Section \ref{sect_well_beh_sys_impr} we will derive a simple
criterion of the irreducibility of the induced representation by
showing that $\Ind\chi,\ \chi\in\cBp,$ is irreducible if and only if
the stabilizer group $\St\chi$ is trivial.

\section{Well-behaved representations}\label{sect_well_beh}
There is an essential difference between unbounded and bounded
representation theory of $*$-algebras in Hilbert space. The problem
of classifying \textit{all} or even \textit{all self-adjoint}
unbounded $*$-representations is not well-posed for arbitary
$*$-algebras. Let us explain this for the $*$-algebra $\dC[x_1,x_2]$
of polynomials in two variables. In \cite{s3} it was proved that for
any properly infinite von Neumann algebra $\cN$ on a separable
Hilbert space there exists a self-adjoint $*$-representation $\pi$
of $\dC[x_1,x_2]$ such that the operators $\overline{\pi(x_1)}$ and
$\overline{\pi(x_2)}$ are self-adjoint and their spectral
projections generate $\cN$. This result has been used in
{\cite{sam_tur} to show the representation theory of $\dC[x_1,x_2]$
is wild. Such a pathological behavior can be overcome if we restrict
to integrable representations. For the $*$-algebra $\dC[x_1,x_2]$ a
self-adjoint representation $\pi$ is integrable if and only the
operators $\overline{\pi(x_1)}$ and $\overline{\pi(x_2)}$ are
self-adjoint and their spectral projections commute. However, for
arbitrary $*$-algebras no method is known to single out such a class
of well-behaved representations. To define and classify well-behaved
representations of general $*$-algebras is a fundamental problem in
unbounded representation theory. One possible proposal was given in
\cite{s2}. In this section we develop a concept of well-behaved
representations for $G$-graded $*$-algebras $\cA$ with commutative
$*$-subalgebras $\cA_e$.
We begin with some simple technical facts.
\begin{lemma}\label{lemma_prelim}
Let $\pi$ be a $*$-representation of a $G$-graded $*$-algebra $\cA$ and $ \cB=\cA_e$. Then the graph topologies of $\pi$ and of $\Res_\cB\pi$ coincide. In particular, $\pi$ is closed if and only if $\Res_\cB\pi$ is closed.
\end{lemma}
\noindent\textbf{Proof.} Since $\cB$ is a $*$-subalgebra of
$\cA$, the graph topology of $\Res_\cB\pi$ is obviously weaker than
that of $\pi$. For $a_g\in\cA_g$ and $\varphi\in \cD(\pi)$, we have
$$
\norm{\pi(a_g)\varphi}=\langle\pi(a_g^*a_g)\varphi,\varphi\rangle^{1/2}\leq\norm{\pi(a_g^*a_g^{})\varphi}+\norm{\varphi}.
$$
Since $a_g^*a_g^{}\in\cB,$ the graph topology of
$\pi$ is weaker than the graph topology of $\Res_\cB\pi.$ Hence both
topologies coincide. Since closedness of a $*$-representation is equivalent to the completeness in the graph topology (see \cite{s1}, 8.1), it follows that $\pi$ is closed if and only if $\Res_\cB\pi$ is closed.
 \hfill $\Box$

\mn
\textit{Throughout the rest of this section we assume that $\cA= \oplus_{g \in
G}~ \cA_g$ is a $G$-graded $*$-algebra such that $\cA_e=\cB$ is
commutative and condition (\ref{eq_cond_compatible}) is
satisfied}.

We begin with some preliminaries. An element $b\in\cB$ can be viewed
as a function $f_b$ on the set $\cBp$, that is, $f_b(\chi)=\chi(b)$
for $ b\in\cB$ and $ \chi\in\cBp$. Let $\tau$ denote the weakest
topology on the set $\cBp$ for which all functions $f_b,\ b\in\cB,$
are continuous. This topology is generated by the sets
$f_b^{-1}((c,d)),\ -\infty\leq c\leq d\leq\infty.$
Clearly, the topology $\tau$ on $\cBp$ is Hausdorff. We assume in
addition that the topology $\tau$ on $\cBp$ is \textit{locally
compact.}

The topology $\tau$ on $\cBp$ defines a Borel structure which is
generated by all open sets.
Since the domain $\cD_g$ of the mapping $\alpha_g$ is the union of open sets
$f_{a_g^*a_g^{}}^{-1}((0,+\infty)),\ a_g\in\cA_g,$ the set $\cD_g$
is open and hence Borel.

\begin{lemma}\label{lemma_topol_cBp_ind_by_a_g*a_g}
Let $\tau_g,\ g\in G,$ be the weakest topology on $\cD_g$ for which
all functions $f_{a_g^*a_g^{}},\ a_g\in\cA_g,$ are continuous. Then
$\tau_g$ is induced from the topology $\tau$ on $\cBp.$
\end{lemma}
\noindent\textbf{Proof.} Let $\chi\in\cD_g.$ Since the topology
$\tau$ on $\cBp$ is locally compact, there is a compact neighborhood
$\Omega$ of $\chi.$ Since $\cD_g$ is open,
$\Omega_1=\Omega\cap\cD_g$ is again a neighborhood of $\chi.$ The
elements of $\cB$ separate the points of $\cBp$. The set
$\set{b^2|b=b^*,\ b\in\cB}$ generates $\cB,$ so it also separates
the points of $\cBp.$ It follows that the set $\set{a_g^*a_g^{},\
a_g\in\cA_g}$ separates the elements of $\cD_g.$ Since
$\Omega$ is compact, $\overline{\Omega}_1$ is also compact. Since
the functions $f_{a_g^*a_g^{}}$ are continuous on $\Omega_1$ and
vanish on the set $\overline{\Omega}_1\backslash\Omega_1,$ they
belong to the $C^*$-algebra $C_0(\Omega_1)$ of continuous functions
vanishing at infinity. By the Stone-Weierstra\ss\ theorem, the
functions $f_{a_g^*a_g^{}}$, where $ a_g\in\cA_g$, generate a
$*$-algebra which is dense in $C_0(\Omega_1).$ Hence the induced
topology of $\tau_g$ on $\Omega_1$ coincides with the induced
topology of $\tau.$ Since $\chi\in\cD_g$ is arbitrary, $\tau_g$ is
induced from the topology $\tau$ on $\cBp.$ \hfill $\Box$\mn

For $\Delta\subseteq\cBp$ and $g\in G,$ we define $\Delta^g$ by
$$\Delta^g=\set{\chi^g|\chi\in\cD_g\cap\Delta}.$$
By definition, $\emptyset_{}^g$ is $\emptyset.$ In particular, if
$\Delta\cap\cD_g=\emptyset,$ then $\Delta^g=\emptyset.$ We also
write $\alpha_g(\Delta)$ for $\Delta^g.$

\begin{lemma}\label{prop_alpha_g_maps_Borel_to_Borel}$ $
\begin{enumerate}
  \item[$(i)$] For any $g\in G,$ the mapping $\alpha_g$ is a homeomorphism of $\cD_g$ onto $\cD_{g^{-1}}.$
  \item[$(ii)$] If $\Delta\subseteq\cD_g$ is open (resp. Borel), then $\Delta^g$ is open (resp. Borel).
\end{enumerate}
\end{lemma}
\noindent\textbf{Proof.} $(i)$: By Proposition
\ref{prop_partial_action}, $\alpha_g$ is a bijection. The equality
$f_{a_g^*a_g}(\chi)=f_{a_g^{}a_g^*}(\chi^g),\ a_g\in\cA_g,$ implies
that for every open subset $X$ of $\dR$ the set
$(f_{a_g^*a_g}^{-1}(X))^g=f_{a_g^{}a_g^*}^{-1}(X)$ is open.
Therefore, by Lemma \ref{lemma_topol_cBp_ind_by_a_g*a_g},
$\alpha_{g^{-1}}$ is continuous. Replacing $g$ by $g^{-1}$ we
conclude that $\alpha_g$ is continuous. Since $\alpha_g$ and
$\alpha_{g^{-1}}$ are inverse to each other, $\alpha_g$ is a
homeomorphism.

$(ii)$: As noted above, $\cD_g$ is open. Therefore, if $\Delta$ is
open (resp. Borel), then $\Delta\cap\cD_g$ is open (resp. Borel).
Since $\alpha_g$ is a homeomorphism, $\Delta^g=(\Delta\cap\cD_g)^g$
is also open (resp. Borel). \hfill $\Box$\mn

After these preliminaries we are ready to give the main definition
of this section.
\begin{defn}\label{defn_well_beh}
A $*$-representation $\pi$ of $\cA$ is \textit{well-behaved} if the
following two conditions are satisfied:
\begin{enumerate}
\item[$(i)$] The restriction $\Res_\cB\pi$ of $\pi$ to $\cB$ is
integrable and there exists a spectral measure $E_\pi$ on the
locally compact space $\cBp[\tau]$ such that
\begin{gather}\label{pispec}
\overline{\pi(b)}=\int_\cBp f_b(\chi)d
E_\pi(\chi)\ \mbox{for }\ b\in\cB.
\end{gather}

\item[$(ii)$] For all $a_g\in\cA_g,\ g\in G,$ and all Borel
subsets $\Delta\subseteq\cBp$, we have
\begin{gather*}
E_\pi(\Delta^g)\pi(a_g)\subseteq\pi(a_g)E_\pi(\Delta).
\end{gather*}
\end{enumerate}
\end{defn}

If $(i)$ is fulfilled, we shall say that the spectral measure
$E_\pi$ is associated with $\pi.$

We give some equivalent forms of the conditions in Definition
\ref{defn_well_beh}. From Theorem \ref{thm_bhat} in the Appendix it
follows that condition $(i)$ is already fulfilled if $\Res_\cB\pi$
is integrable and $\cB$ is countably generated. The next proposition
contains a number of reformulations of condition $(ii)$.

\begin{prop}\label{prop_comm_rel_for_well_beh}
Let $\pi$ be a $*$-representation of $\cA$ satisfying condition
$(i)$ of Definition \ref{defn_well_beh}. Let $\cF_\pi$ denote the
set of Borel functions $f$ on $\cBp$ such that the operator $\int
fdE_\pi$ maps the domain $\cD(\pi)$ into itself. For $a_g\in\cA_g,\
g\in G,$ let $U_gC_g$ be the polar decomposition of
$\overline{\pi(a_g)}$. Then the following statements are equivalent:

$(i):$ Condition $(ii)$ of Definition \ref{defn_well_beh} is
fulfilled.

$(ii):$ For all $a_g\in\cA_g,\ g\in G,$ and all Borel sets
$\Delta\subseteq\cBp$ we have $U_gE_\pi(\Delta)=E_\pi(\Delta^g)U_g$.

$(iii):$ For any $E$-measurable function $f$ on $\cBp$ and $a_g\in\cA_g,\ g\in G,$ we have
\begin{gather*}\label{eq_comm_rel_pi(a_g)int f}
U_g\int f(\chi)dE_\pi(\chi)\subseteq
\int_{\cD_{g^{-1}}}f(\alpha_{g^{-1}}(\chi))dE_\pi(\chi)U_g.
\end{gather*}

$(iv):$ For any $f\in \cF_\pi$, $a_g\in\cA_g,\ g\in G,$ and $\varphi
\in \cD(\pi)$, we have
\begin{gather*}\label{eq_comm_rel_pi(a_g)int f}
\pi(a_g)\int f(\chi)dE_\pi(\chi) \varphi =
\int_{\cD_{g^{-1}}}f(\alpha_{g^{-1}}(\chi))dE_\pi(\chi)\pi(a_g)\varphi.
\end{gather*}
\end{prop}
\noindent\textbf{Proof.} $(i)\Rightarrow(ii):$ Fix
$\Delta\subseteq\cBp$. Since $\Res_\cB\pi$ is integrable,
$\overline{\pi(a_g^*a_g)}$ is self-adjoint. But
$\pi(a_g)^*\overline{\pi(a_g)}$ is self-adjoint extension of
$\overline{\pi(a_g^*a_g)}$, so that
$C_g^2=\pi(a_g)^*\overline{\pi(a_g)}=\overline{\pi(a_g^*a_g)}$.
Since $\overline{{\pi(a_g^*a_g)}}$ commutes with the projections
$E_\pi(\cdot),$ $C_g^2$ and hence $C_g$ commute with $E_\pi(\cdot).$
Thus we get $U_gE(\Delta^g)C_g\subseteq U_gC_g E(\Delta)=
\overline{\pi(a_g)}E(\Delta).$ From Definition \ref{defn_well_beh},
$(i)$ it follows that the kernel of $C_g^2=\overline{\pi(a_g^*a_g)}$
is equal to $\Ran E_\pi(f_{a_g^*a_g^{}}^{-1}(0))$. By the properties
of the polar decomposition, this kernel equals to $\ker U_g=\ker
C_g.$ If $v\in\ker C_g,$ then $E(\Delta^g)U_gv=0$ and, since
$P_0:=E_\pi(f_{a_g^*a_g^{}}^{-1}(0))$ commutes with $E_\pi(\cdot)$,
we get $U_gE(\Delta)v=U_gE(\Delta)P_0v=U_gP_0E(\Delta)v=0.$ Thus the
bounded operators $U_gE(\Delta)$ and $E(\Delta^g)U_g$ coincide on
the dense set $\Ran C_g {+}\ker C_g,$ so they coincide everywhere.

$(ii)\Rightarrow(iii):$ From $(ii)$ we get $(iii)$ for
characteristic functions, then for simple functions and by a limit
procedure for arbitrary measurable functions $f \in \cF_\pi$.

$(iii)\Rightarrow(iv):$ This follows from the relation
$\pi(a_g)\varphi=U_gC_g\varphi$ combined with the fact that $C_g$
and the first integral commute on vectors $\varphi \in \cD(\pi)$.

$(iv)\Rightarrow(i):$ Since $\pi$ is integrable, $\pi$ is closed and so is $\Res_\cB \pi$ by Lemma \ref{lemma_prelim}. Therefore, $\cD(\pi)=\cap_{b\in \cB} \cD(\overline{\pi(b)})$. By (\ref{pispec}) the latter implies that $E(\Delta)$ leaves the domain $\cD(\pi)$ invariant. Hence the characteristic function of $\Delta $ belongs to $\cF_\pi$ and (i) follows from (iv) applied to this characteristic function. \hfill $\Box$\mn

Many notions on unbounded operators are derived from appropriate
reformulations of the corresponding notions on bounded operators.
The next proposition says that bounded $*$-representations satisfy
the two conditions in Definition \ref{defn_well_beh}. This
observation was in fact the starting point for our definition
of well-behaved representations.

\begin{prop}\label{prop_well_behav_bounded}
If $\pi$ is a bounded $*$-representation of the $*$-algebra $\cA$ such that $\cD(\pi)=\cH_\pi$,
then $\pi$ is well-behaved.
\end{prop}
\noindent\textbf{Proof.} Since the representation $\pi$ is bounded,
the closure of $\pi(\cB)$ in the operator norm is a commutative
$C^*$-algebra. Hence condition $(i)$ follows from Theorem 12.22 in
\cite{rud}.

Fix $g\in G,\ a_g,b_g\in\cA_g.$ From assumption
(\ref{eq_cond_compatible}) we obtain that
$f_{a_g^*a_g^{}b_g^*b_g^{}}(\chi)=f_{a_g^*b_g^{}b_g^*a_g^{}}(\chi)$
on $\cBp$. Therefore, by condition $(i)$ we have
$\pi(a_g^*a_g^{}b_g^*b_g^{})=\pi(a_g^*b_g^{}b_g^*a_g^{})$ which can
be rewritten in the form
\begin{gather}\label{eq_aux}
\pi(a_g^*)\pi(a_g^{}b_g^*b_g^{})=\pi(a_g^*)\pi(b_g^{}b_g^*a_g^{}).
\end{gather}
Since $\pi(b_g^{}b_g^*)$ commutes with $\pi(a_g^*a_g^{}),$ it also
commutes with the projection onto the range of $\pi(a_g^{}).$ This
implies that $\pi(b_g^{}b_g^*)(\Ran(\pi(a_g^{})))$ is contained in
$\overline{\Ran(\pi(a_g^{}))},$ so the range of the operator
$\pi(b_g^{}b_g^*a_g)$ is contained in
$\overline{\Ran(\pi(a_g^{}))}.$ The range of the operator
$\pi(a_g^{}b_g^*b_g^{})$ is evidently contained in
$\Ran\pi(a_g^{}).$ From the relation
$\overline{\Ran(\pi(a_g^{}))}=\ker(\pi(a_g^*))^\perp$ it follows
that $\pi(a_g^*)$ restricted to $\overline{\Ran(\pi(a_g^{}))}$ is
injective. Therefore, from (\ref{eq_aux}) we get
$\pi(a_g^{}b_g^*b_g^{})=\pi(b_g^{}b_g^*a_g^{})$ and so
\begin{gather*}
\pi(a_g^{})\pi(b_g^*b_g^{})=\pi(b_g^{}b_g^*)\pi(a_g^{})
\end{gather*}
for all $b_g \in \cA_g$. Now we use a standard approximation
procedure. The preceding relation yields
\begin{gather*}
\pi(a_g^{})p_n(\pi(b_g^*b_g^{}))=p_n(\pi(b_g^{}b_g^*))\pi(a_g^{})
\end{gather*}
for all polynomials $p_n\in\dC[t]$ which implies that
\begin{gather*}
\pi(a_g^{})E_{\pi(b_g^*b_g^{})}(X)=E_{\pi(b_g^{}b_g^*)}(X)\pi(a_g^{}),
\end{gather*}
where $E_{\pi(\cdot)}$ denotes the spectral measure of the
self-adjoint operator $\pi(\cdot)$ and $X$ is a Borel subset of
$\dR$. The spectral measure $E_\pi$ on the space $\cBp$ associated
with $\pi$ is releated to the spectral measure of the operator
$\pi(b_h^*b_h^{}),\ b_h^{}\in\cA_h,\ h\in G,$ by the equation
$$
E_{\pi(b_h^*b_h^{})}(X)=E_\pi(f_{b_h^*b_h^{}}^{-1}(X)),
$$
where $f_{b_h^*b_h^{}}$ is the function on $\cBp$ defined by the
element $b_h^*b_h^{}\in\cB.$ From the equality
$$
\alpha_h(f_{b_h^*b_h^{}}^{-1}(X))=f_{b_h^{}b_h^*}^{-1}(X)
$$
we obtain
\begin{gather}\label{eq_aux1}
\pi(a_g^{})E_\pi(\Delta)=E_\pi(\Delta^g)\pi(a_g^{}),%
\end{gather}
where $g\in G,a_g\in\cA_g$, $\Delta=f_{c_g^*c_g}^{-1}(X)$, and $X$
is a Borel subset $\dR.$ Since (\ref{eq_aux1}) is valid for such
sets $\Delta$, it holds for the all sets from the $\sigma$-algebra
generated by the sets $\Delta$ as well. From Lemma
\ref{lemma_topol_cBp_ind_by_a_g*a_g} we conclude that
(\ref{eq_aux1}) holds for all Borel sets $\Delta\subseteq\cD_g.$

In particular, equation (\ref{eq_aux1}) is true for $\Delta=\cD_g,$
so also for $\Delta=\cBp\backslash\cD_g.$ Therefore we have
$\pi(a_g)E_\pi(\cBp\backslash\cD_g)=0$ which implies that
$\pi(a_g)E_\pi(\Delta_0)=0$ for all Borel subsets
$\Delta_0\subseteq\cBp\backslash\cD_g.$ Since
$E_\pi(\alpha_g(\Delta_0))=E_\pi(\emptyset)=0$, (\ref{eq_aux1}) is
valid for all Borel sets $\Delta_0$ of $\cBp\backslash\cD_g.$ Hence
condition $(ii)$ of Definition \ref{defn_well_beh} is satisfied.
\hfill $\Box$\mn

In the rest of this section we derive some basic properties of
well-behaved representations.
\begin{prop}\label{prop_self_adj_subrep_is_well_beh}
Let $\pi$ be a well-behaved representation of $\cA.$ Then any
self-adjoint subrepresentation $\pi_0\subseteq\pi$ is well-behaved.
\end{prop}
\noindent\textbf{Proof.} Since $\pi$ is well-behaved, it is
self-adjoint. By Corollary 8.3.13 in \cite{s1}, there exists a
representation $\pi_1$ of $\cA$ such that $\pi=\pi_0\oplus\pi_1.$
Since $\Res_\cB\pi$ is integrable, $\Res_\cB\pi_0$ is integrable by
Proposition 9.1.17 $(i)$ in \cite{s1}. Let $P\in\pi(\cA)'$ denotes
the projection on the representation space $\cH_{\pi_0}$ of $\pi_0.$
Then $PE_\pi(\cdot)\upharpoonright\cH_{\pi_0}$ is a spectral measure
$E_{\pi_0}(\cdot)$ associated with $\pi_0$. Let $a_g\in\cA_g,\ g\in
G,$ and let $\Delta$ be a Borel subset in $\cBp$ such that
$\Delta^g$ is defined. Suppose that $\varphi \in \cD(\pi_0)$. Using
Definition \ref{defn_well_beh}, $(ii)$ for $\pi$ we obtain
$$
E_{\pi_0}(\Delta^g)\pi_0(a_g)\varphi
=PE_{\pi}(\Delta^g)\pi(a_g)\varphi=P\pi(a_g)E_{\pi}(\Delta)\varphi=
\pi_0(a_g)E_{\pi_0}(\Delta)\varphi,
$$
that is, $E_{\pi_0}(\Delta^g)
\pi_0(a_g)\subseteq\pi_0(a_g)E_{\pi_0}(\Delta)$, so condition $(ii)$
of Definition \ref{defn_well_beh} holds for $\pi_0$. Hence $\pi_0$
is well-behaved. \hfill $\Box$

\begin{lemma}
As above, $H$ denotes a subgroup of $G$. Let $\rho$ be a well-behaved inducible representation of $\cA_H$,
$E_\rho$ a spectral measure on $\cBp$ associated with $\rho$ and
$\pi$ the induced representation $\Ind_{\cA_H\uparrow\cA}\rho$.
Suppose that $b\in\cB$ and $g\in G.$ Then the domain of the operator
$\int_{\cD_g}f_b(\alpha_g(\chi))dE_\rho(\chi)$ contains $\cD(\rho)$
and for arbitrary
$a_g\in\cA_g$ and $ v\in\cD(\rho)$ we have
\begin{gather}\label{eq_ResInd_rho_action}
\pi(b)[a_g\otimes{}v]=[ba_g\otimes{}v]=[a_g\otimes\left(\int_{\cD_g}f_b(\alpha_g(\chi))dE_\rho(\chi)\right)v].
\end{gather}
\end{lemma}
\noindent\textbf{Proof.} Let $[c_g\otimes w]\in\cH_\pi$, where
$c_g\in\cA_g,\ w\in\cD(\rho).$ Then we have
\begin{gather*}
\langle \pi(b)[a_g\otimes{}v],[c_g \otimes w] \rangle =
\langle[ba_g\otimes v],[c_g\otimes
w]\rangle=\langle\rho(c_g^*ba_g^{})v,w\rangle=\int_{\cBp}
f_{c_g^*ba_g^{}}(\chi)d \langle E_\rho(\chi)v,w\rangle.
\end{gather*}
From Proposition \ref{prop_defn_part_act_gener} we obtain the
equalities
$f_{c_g^*ba_g^{}}(\chi)=f_{b}(\alpha_g(\chi))f_{c_g^*a_g^{}}(\chi)$
for $\chi\in\cD_g$ and $f_{c_g^*ba_g^{}}(\chi)=0$ for
$\chi\in\cBp\backslash\cD_g$, so the preceding is equal to
\begin{gather*}
\int_{\cD_g}f_b(\alpha_g(\chi))f_{c_g^*a_g^{}}(\chi)d\langle
E_\rho(\chi)v,w\rangle=\langle\left(\int_{\cD_g}f_b(\alpha_g(\chi))f_{c_g^*a_g^{}}(\chi)dE_\rho(\chi)\right)v,w\rangle.
\end{gather*}
Since $v$ belongs to the domains of
$\int_{\cD_g}f_b(\alpha_g(\chi))f_{c_g^*a_g^{}}(\chi)dE_\rho(\chi)$
and $\int_{\cD_g}f_{c_g^*a_g^{}}(\chi)dE_\rho(\chi),$ the
multiplicativity property of the spectral integral (see e.g.
\cite{rud}, 13.24) implies that $v$ belongs to the domain of
$\int_{\cD_g}f_b(\alpha_g(\chi))dE_\rho(\chi)$ and we can proceed
\begin{gather*}
\langle\pi(b)[a_g\otimes{}v],[c_g \otimes w]\rangle=
\langle\left(\int_{\cD_g}f_{c_g^*a_g^{}}(\chi)dE_\rho(\chi)\right)\left(\int_{\cD_g}f_b(\alpha_g(\chi))dE_\rho(\chi)\right)v,w\rangle\\
=\langle\rho(c_g^*a_g^{})\left(\int_{\cD_g}f_b(\alpha_g(\chi))dE_\rho(\chi)\right)v,w\rangle=
\langle [a_g\otimes \left(\int_{\cD_g} f_b(\alpha_g(\chi))dE_\rho(\chi)\right) v],[c_g\otimes w]\rangle.%
\end{gather*}
Since the linear span of vectors $[c_g\otimes w]$, where
$c_g\in\cA_g$ and $w\in\cD(\rho)$, is dense in the closed subspace
to which $[ba_g\otimes v]$ and $[a_g\otimes \left(\int
f_b(\alpha_g(\chi))dE_\rho(\chi)\right) v]$ belong, the assertion
follows. \hfill $\Box$

\begin{prop}\label{prop_ind_of_cyc_well_beh_is_well_beh}
Assume that $\cB$ is countably generated. If $\rho$ is a well-behaved inducible cyclic representation of the $*$-algebra $\cA_H,$ then the induced representation $\pi=\Ind_{\cA_H\uparrow\cA}(\rho)$ is a well-behaved representation of the $*$-algebra $\cA.$
\end{prop}
\noindent\textbf{Proof.} Let $E_\rho$ be a spectral measure on $\cBp$ associated with $\rho.$ It follows from the Theorem \ref{thm_bhat}, $(ii)$ that $E_\rho$ is supported on $\cBp.$ We first show that $\Res_\cB\pi$ is defined by a spectral measure, i.e. condition $(i)$ in Definition \ref{defn_well_beh} holds for some spectral measure $E_\pi$ on $\cBp.$

Let $a_g\in\cA_g,\ g\in G,\ w\in\cD(\rho),$ and let $\Delta$ be a Borel subset of $\cBp.$ We define a linear operator $E_\pi(\Delta)$ on the tensor product $\cA\otimes\cD(\rho)$ by putting $E_\pi(\Delta)(a_g\otimes w):=a_g\otimes E_\rho(\Delta^{g^{-1}})w.$ Note that the vector $E_\rho(\Delta^{g^{-1}})w$ belongs to $\cD(\rho).$ Let $h\in H$ and $a_h\in\cA_h.$ Using Proposition \ref{prop_comm_rel_for_well_beh} $(i)$ we get that
\begin{gather*}
E_\pi(\Delta)(a_ga_h\otimes{}w-a_g\otimes\rho(a_h)w)=a_ga_h\otimes{}E_\rho(\Delta^{h^{-1}g^{-1}})w-a_g\otimes{}E_\rho(\Delta^{g^{-1}})\rho(a_h)w=\\
=a_ga_h\otimes{}E_\rho(\Delta^{h^{-1}g^{-1}})w-a_g\otimes{}\rho(a_h)E_\rho(\Delta^{h^{-1}g^{-1}})w,
\end{gather*}
belongs to the kernel of the quotient mapping $\cA\otimes\cD(\rho)\to\cA\otimes_{\cA_H}\cD(\rho),$ so $E_\pi(\Delta)$ defines a linear operator on $\cA\otimes_{\cA_H}\cD(\rho)$ which we denote again by $E_\pi(\Delta).$

Let $v\in\cD(\rho)$ be a cyclic vector for $\rho.$ Take $a\otimes
v\in\cA\otimes_{\cA_H}\cD(\rho).$ We write $a$ as a finite sum
$\sum_{i,k} a_{ik},\ a_{ik}\in\cA_{g_{ik}},$ where $g_{ik}\in G$ are
pairwise distinct and $g_{ik}^{-1}g_{jm}\in H$ if and only if $k=m.$
Then we have $\langle a_{ik}\otimes v,a_{jm}\otimes v\rangle_0=0$
for $k\neq m$ and remembering that $\rho$ is well-behaved we get
\begin{gather}
\nonumber\langle E_\pi(\Delta)(a\otimes v),E_\pi(\Delta)(a\otimes{}v)\rangle_0=%
\langle\sum_{i,k} a_{ik}\otimes E_\rho(\Delta^{g_{ik}^{-1}})v,\sum_{i,k}a_{ik}\otimes E_\rho(\Delta^{g_{ik}^{-1}})v\rangle_0=\\
\nonumber=\sum_{k}\langle\sum_i a_{ik}\otimes E_\rho(\Delta^{g_{ik}^{-1}})v,\sum_i a_{ik}\otimes E_\rho(\Delta^{g_{ik}^{-1}})v\rangle_0%
=\sum_k\sum_{i,j}\langle\rho(a_{kj}^*a_{ki}^{})E_\rho(\Delta^{g_{ik}^{-1}})v,E_\rho(\Delta^{g_{jk}^{-1}})v\rangle=\nonumber\\
\label{eq_aux6}=\sum_k\sum_{i,j}\langle E_\rho(\Delta^{g_{jk}^{-1}})\rho(a_{kj}^*a_{ki}^{})v,E_\rho(\Delta^{g_{jk}^{-1}})v\rangle%
=\sum_k\sum_{i,j}\langle\rho(a_{kj}^*a_{ki}^{})v,E_\rho(\Delta^{g_{jk}^{-1}})v\rangle=\\
\nonumber=\langle a\otimes v,E_\pi(\Delta)(a\otimes{}v)\rangle_0%
\end{gather}
Assume that $a\otimes v\in\cK_\rho,$ that is, $\langle a\otimes
v,a\otimes v\rangle_0=0.$ The preceding calculation implies that
$E_\pi(\Delta)(a\otimes v)\in\cK_\rho,$ so $E_\pi(\Delta)$ is a
well-defined linear operator on the linear span of vectors
$[a\otimes v]\in\cD(\pi)$ defined by
\begin{gather}\label{eq_action_ind_spec_measure}
E_\pi(\Delta)[a_g\otimes v]:=[a_g\otimes E_\rho(\Delta^{g^{-1}})v].
\end{gather}

Since $v$ is cyclic, the set of vectors $[a\otimes v]$ is dense in
$\cH_\pi$ by Lemma \ref{lemma_cyc_vec_of_ind_rep}. It follows from
(\ref{eq_aux6}) that $E_\pi(\Delta)$ is bounded and can be
extended by continuity to $\cH_\pi.$ From now on we consider
$E_\pi(\Delta)$ on the subspace $\cH_\pi.$

It can be easily seen that $E_\pi(\Delta)^2=E_\pi(\Delta).$ We
prove that $E_\pi(\Delta)$ is self-adjoint. For this it suffices to
show that
\begin{gather}\label{eq_self_adj_of_ind_measure}
\langle E_\pi(\Delta)[a_{g_1}\otimes{}v],[a_{g_2}\otimes{}v]\rangle=
\langle [a_{g_1}\otimes{}v],E_\pi(\Delta)[a_{g_2}\otimes{}v]\rangle
\end{gather}
for $a_{g_1}\in\cA_{g_1},\ a_{g_2}\in\cA_{g_2},\ g_1,g_2\in G.$
First we consider the case when $g_1H\neq g_2H.$ Then we get
\begin{gather*}
\langle E_\pi(\Delta)[a_{g_1}\otimes v],[a_{g_2}\otimes v]\rangle=
\langle [a_{g_1}\otimes E_\rho(\Delta^{g_1^{-1}})v],[a_{g_2}\otimes v]\rangle=\\
=\langle\rho(p_H(a_{g_2}^*a_{g_1}^{}))E_\rho(\Delta^{g_1^{-1}})v,v\rangle=0,
\end{gather*}
since $p_H(a_{g_2}^*a_{g_1}^{})=0.$ Analogously, $\langle
[a_{g_1}\otimes v],E_\pi(\Delta)[a_{g_2}\otimes v]\rangle=0,$ so
that (\ref{eq_self_adj_of_ind_measure}) holds in this case.
Now suppose that $g_1H=g_2H.$ Then we have
\begin{gather*}
\langle E_\pi(\Delta)[a_{g_1}\otimes v],[a_{g_2}\otimes v]\rangle= %
\langle[a_{g_1}\otimes E_\rho(\Delta^{g_1^{-1}})v],[a_{g_2}\otimes v]\rangle=\langle\rho(a_{g_2}^*a_{g_1}^{})E_\rho(\Delta^{g_1^{-1}})v,v\rangle.
\end{gather*}
Since $\rho$ is well-behaved and
$a_{g_2}^*a_{g_1}^{}\in\cA_{g_2^{-1}g_1^{}}$, the preceding equals to
\begin{gather*}
=\langle E_\rho(\Delta^{g_2^{-1}})\rho(a_{g_2}^*a_{g_1}^{})v,v\rangle %
=\langle \rho(a_{g_2}^*a_{g_1}^{})v,E_\rho(\Delta^{g_2^{-1}})v\rangle %
=\langle [a_{g_1}\otimes v],E_\pi(\Delta)[a_{g_2}\otimes v]\rangle.%
\end{gather*}
Thus, $E_\pi(\Delta)$ is self-adjoint.

Take $a_g\in\cA_g,$ a Borel set $\Delta\subseteq\cBp$ and
$a_k\in\cA_k.$ Then we get
\begin{gather}\label{eq_aux5}
\pi(a_g)E_\pi(\Delta)[a_k\otimes v]=\pi(a_g)[a_k\otimes E_\rho(\Delta^{k^{-1}})v]=[a_ga_k\otimes E_\rho(\Delta^{k^{-1}})v]=\\
\nonumber[a_ga_k\otimes E_\rho((\Delta^g)^{(gk)^{-1}})v]=E_\pi(\Delta^g)[a_ga_k\otimes v]=E_\pi(\Delta^g)\pi(a_g)[a_k\otimes v].%
\end{gather}

Next we prove that $E_\pi(\Delta)\cD(\pi)\subseteq\cD(\pi).$
Take $d_g\in\cA_g,\ g\in G.$ Using (\ref{eq_aux5}) we
obtain
\begin{gather*}
\norm{E_\pi(\Delta)[a\otimes v]}_{d_g}^2=\norm{\pi(d_g)E_\pi(\Delta)[a\otimes v]}^2=%
\langle\pi(d_g)E_\pi(\Delta)[a\otimes v],\pi(d_g)E_\pi(\Delta)[a\otimes v]\rangle=\\
=\langle E_\pi(\Delta^g)\pi(d_g)[a\otimes v],E_\pi(\Delta^g)\pi(d_g)[a\otimes v]\rangle=%
\langle\pi(d_g)[a\otimes v],E_\pi(\Delta^g)\pi(d_g)[a\otimes v]\rangle=\\
=\langle\pi(d_g),\pi(d_g)E_\pi(\Delta)[a\otimes v]\rangle\leq%
\norm{[a\otimes v]}_{d_g}\cdot\norm{E_\pi(\Delta)[a\otimes v]}_{d_g},%
\end{gather*}
and hence $\norm{E_\pi(\Delta)[a\otimes v]}_{d_g}\leq\norm{[a\otimes
v]}_{d_g}.$ By Lemma \ref{lemma_cyc_vec_of_ind_rep}, the set of vectors $[a\otimes v]$ is a core
for $\pi$ . Therefore, the preceding shows that $E_\pi(\Delta)$ is
continuous in the graph topology of $\pi$. This in turn implies that
$E_\pi(\Delta)\cD(\pi)\subseteq\cD(\pi)$.

Now we prove that $E_\pi(\cdot)$ defines a spectral measure on
$\cBp.$ For $a_g\in\cA_g$ we have
\begin{gather*}
\langle E_\pi(\cBp)[a_g\otimes v],[a_g\otimes v]\rangle=\langle[a_g\otimes E_\rho(\cD_g)v],[a_g\otimes v]\rangle=\\
=\langle\rho(a_g^*a_g^{})E_\rho(\cD_g)v,v\rangle=\langle\rho(a_g^*a_g^{})v,v\rangle=\langle[a_g\otimes v],[a_g\otimes v]\rangle
\end{gather*}
which shows that $E_\pi(\cBp)=I$. The countable additivity
$E_\pi(\cdot)$ follows at once from the countable additivity of
$E_\rho(\cdot).$

Next we show that $\Res_\cB\pi$ is an integrable representation
associated with spectral measure $E_\pi.$ It suffices to prove that
\begin{gather}\label{eq_aux3}
\langle b[a_{g_1}\otimes v],[a_{g_2}\otimes v]\rangle=\int
f_b(\chi)d\langle E_\pi(\chi)[a_{g_1}\otimes v],[a_{g_2}\otimes
v]\rangle.
\end{gather}
for all $[a_{g_1}\otimes v],[a_{g_2}\otimes v]\in\cH_\pi$. In the
case $g_1H\neq g_2H$ one easily checks that the both sides of
(\ref{eq_aux3}) are equal to zero. In the case $g_1H=g_2H$ we use
(\ref{eq_ResInd_rho_action}) and compute
\begin{gather*}
\langle\pi(b)[a_{g_1}\otimes v],[a_{g_2}\otimes
v]\rangle=\langle[a_{g_1}\otimes\int_{\cD_{g_1}}f_b(\alpha_{g_1}(\chi))dE_\rho(\chi)\otimes
v],[a_{g_2}\otimes v]\rangle\\
=\langle\rho(a_{g_2}^*a_{g_1}^{})\int_{\cD_{g_1}}f_b(\alpha_{g_1}(\chi))dE_\rho(\chi)v,v\rangle.
\end{gather*}
Applying Proposition \ref{prop_comm_rel_for_well_beh} $(iv)$ we continue
\begin{gather*}
=\langle\int_{\cD_{g_2}}f_b(\alpha_{g_2}(\chi))dE_\rho(\chi)\rho(a_{g_2}^*a_{g_1}^{})v,v\rangle=
\int_{\cD_{g_2}}f_b(\alpha_{g_2}(\chi))d\langle
E_\rho(\chi)\rho(a_{g_2}^*a_{g_1}^{})v,v\rangle\\
=\int_{\cD_{g_2}}f_b(\alpha_{g_2}(\chi))d\langle\rho(a_{g_2}^*a_{g_1}^{})E_\rho(\alpha_{g_1^{-1}g_2^{}}(\chi))v,v\rangle=
\int_{\cD_{g_2^{-1}}}f_b(\chi)d\langle\rho(a_{g_2}^*a_{g_1}^{})E_\rho(\alpha_{g_1^{-1}}(\chi))v,v\rangle\\
=\int_{\cD_{g_2^{-1}}}f_b(\chi)d\langle [a_{g_1}^{}\otimes E_\rho(\alpha_{g_1^{-1}}(\chi))v],[a_{g_2}\otimes v]\rangle=%
\int_{\cD_{g_2^{-1}}}f_b(\chi)d\langle E_\pi(\chi)[a_{g_1}^{}\otimes v],[a_{g_2}\otimes v]\rangle\\
=\int_{\cD_{g_2^{-1}}}f_b(\chi)d\langle [a_{g_1}^{}\otimes v],[a_{g_2}\otimes E_\rho(\alpha_{g_2^{-1}}(\chi))v]\rangle=%
\int_{\cBp}f_b(\chi)d\langle E_\pi(\chi)[a_{g_1}^{}\otimes v],[a_{g_2}\otimes v]\rangle.%
\end{gather*}

It follows from (\ref{eq_aux5}) that the equality in the Definition
\ref{defn_well_beh}, $(ii)$ holds on the span of vectors $[a\otimes
v]\in\cD(\pi)$ which is a core of $\pi$ by Lemma
\ref{lemma_cyc_vec_of_ind_rep}. Since $\pi(a_g)$ and $E_\pi(\Delta)$
are continuous in the graph topology of $\pi,$ condition $(ii)$ in
Definition \ref{defn_well_beh} holds for $\pi.$ This completes the
proof. \hfill $\Box$ \mn

In what follows, we want to induce from arbitrary
well-behaved representations of subalgebras $\cA_H.$ For this reason
we shall need the decomposition of well-behaved representations into
direct sums of cyclic well-behaved representations. This aim will be
achieved by Proposition \ref{prop_well_beh_is_sum_of_cyclic} below.
First we develop some more preliminaries.
\begin{lemma}\label{lemma_U_liegt_im_Doppelkomm}
Suppose that $\pi$ is a well-behaved representation of $\cA$. Let
$a_g\in\cA_g$ and let $UC$ be the polar decomposition of
$\overline{\pi(a_g)}.$ Then $U$ belongs to $\pi(\cA)''.$
\end{lemma}
\noindent\textbf{Proof.} Let $T\in\pi(\cA)'.$ As noted already in
the proof of Proposition \ref{prop_comm_rel_for_well_beh}, we have
$C^2=\overline{\pi(a_g^*a_g^{})}.$ Since $T$ commutes with
$\pi(a_g^*a_g^{})$, it commutes with $C^2$ and therefore with $C.$

Take $\varphi\in\cD(C).$ Then we obtain
$TU(C\varphi)=T\overline{\pi(a_g)}\varphi=
\overline{\pi(a_g)}T\varphi=UCT\varphi=UT(C\varphi).$ Now let
$\psi\in\ker C=\ker U=\ker\overline{\pi(a_g)}.$ Then we have
$\overline{\pi(a_g)}T\psi=T\overline{\pi(a_g)}\psi=0,$ i.e. $T\ker
U\subseteq \ker U,$ so that $UT\psi=0=TU\psi.$ Therefore, $T$ and
$U$ commute on the linear dense subspace $\ker C+\Ran C$. Since $T$
and $U$ are bounded, they commute on $\cH_\pi$. This shows that
$U\in\pi(\cA)''.$ \hfill $\Box$\mn

\begin{lemma}\label{lemma_aux1}
If $\pi$ is a well-behaved representation of $\cA$, then we have:
\begin{enumerate}
  \item[$(i)$] $\overline{\pi(a_g^*)}=\pi(a_g)^*$ for $a_g\in \cA_g.$
  \item[$(ii)$] $\overline{\pi(a_ga_k)}=\overline{\overline{\pi(a_g)}\cdot\overline{\pi(a_k)}}$%
for $a_g\in \cA_g$ and $a_k\in\cA_k.$
\end{enumerate}
\end{lemma}
\noindent\textbf{Proof.} $(i):$ It is clear that
$\pi(a_g^{}a_g^*)\subseteq\pi(a_g^*)^*\overline{\pi(a_g^*)}$.
Since $\pi$ is well-behaved, $\Res_\cB \pi$ is integrable, so
$\pi(a_g^{}a_g^*)$ is essentially self-adjoint (\cite{s1}, 9.1.2).
Hence it follows that
$\overline{\pi(a_g^{}a_g^*)}=\pi(a_g^*)^*\overline{\pi(a_g^*)}=|\overline{\pi(a_g^*)}|^2.$
By the same reasoning we obtain
$\overline{\pi(a_g^{}a_g^*)}=\overline{\pi(a_g)}\pi(a_g)^*=|\pi(a_g)^*|^2.$
Combining these relations with the fact that $\cD(T)=\cD(|T|)$ for
a closed operator $T$ we get
$$
\cD(\overline{\pi(a_g^*)})=\cD(|\overline{\pi(a_g^*)}|)=
\cD((\overline{\pi(a_g^{}a_g^*)})^{1/2})=\cD(|\pi(a_g)^*|)=\cD(\pi(a_g)^*).
$$
Since $\overline{\pi(a_g^*)}\subseteq\pi(a_g)^*$, the preceding implies that
$\overline{\pi(a_g^*)}=\pi(a_g)^*.$

$(ii):$ Clearly,
$\pi(a_k^*a_g^*a_g^{}a_k^{})\subseteq\left(\overline{\pi(a_g)}
\cdot\overline{\pi(a_k)}\right)^*\overline{\overline{\pi(a_g)}\cdot\overline{\pi(a_k)}}.$
Since $a_k^*a_g^*a_g^{}a_k^{}\in \cB$, the operator
$\overline{\pi(a_k^*a_g^*a_g^{}a_k^{})}$ is self-adjoint, so we have
the equality
$\overline{\pi(a_k^*a_g^*a_g^{}a_k^{})}=\left(\overline{\pi(a_g)}
\cdot\overline{\pi(a_k)}\right)^*\overline{\overline{\pi(a_g)}\cdot\overline{\pi(a_k)}}$
which yields $\cD((\overline{\pi(a_k^*a_g^*a_g^{}a_k^{})})^{1/2})=
\cD(\overline{\overline{\pi(a_g)}\cdot\overline{\pi(a_k)}})$. As
shown in the proof of $(i)$ we also have that
$\cD(\overline{\pi(a_ga_k)})=\cD((\overline{\pi(a_k^*a_g^*a_g^{}a_k^{})})^{1/2}$.
Combining these two equalities with the obvious inclusion
$\overline{\pi(a_ga_k)}\subseteq\overline{\overline{\pi(a_g)}\cdot\overline{\pi(a_k)}}$,
the assertion follows. \hfill $\Box$

\begin{lemma}\label{lemma_commut_polar_decompos}
Let $\pi$ be a well-behaved $*$-representation of $\cA$. We denote
by $\cU_\pi$ the set of all partial isometries in the polar
decompositions of elements $\overline{\pi(a_g)}$, where $ a_g\in\cA_g,\ g\in
G.$ Then
$$
\gA_0=\set{\sum_{i=1}^n \lambda_iU_iE_\pi(\Delta_i):\
\lambda_i\in\dC,\ U_i\in\cU_\pi,\ \Delta_i\subseteq\cBp,\ \Delta_i\
\mbox{is a Borel set}}
$$
is a dense $*$-subalgebra of $\pi(\cA)''$ in the strong operator
topology.
\end{lemma}
\noindent\textbf{Proof.} Since $\cU_\pi\subseteq\pi(\cA)''$ by Lemma
\ref{lemma_U_liegt_im_Doppelkomm} and the spectral projections
$E_\pi(\cdot)$ belong to $\pi(\cB)''\subseteq\pi(\cA)'',$ we
conclude that $\gA_0\subseteq\pi(\cA)''.$

We prove that $\gA_0$ is a $*$-algebra. Take $a_g\in\cA_g$ and let
$U_g|\overline{\pi(a_g)}|$ be the polar decomposition of the closed
operator $\overline{\pi(a_g)}.$ By Lemma \ref{lemma_aux1}, $(i)$ we
have $\overline{\pi(a_g^*)}=\pi(a_g)^*$. It is well-known (see e.g.
\cite{kato}, p. 421), that $U_g^*|\overline{\pi(a_g^*)}|$ is the
polar decomposition of the adjoint operator
$\overline{\pi(a_g^*)}=\pi(a_g)^*$ of $\overline{\pi(a_g)}$.
Therefore, $U_g^*\in\gA_0$ which proves that $\gA_0$ is
$*$-invariant.

Take another element $a_k\in\cA_k,\ k\in G$ and let $U_kC_k$ be the
polar decomposition of $\overline{\pi(a_k)}.$ Then using Lemma
\ref{lemma_aux1} and Proposition \ref{prop_comm_rel_for_well_beh}
$(iii)$ we get
\begin{gather}\label{eq_polar_decompos_pi(a_ga_k)}
\overline{\pi(a_ga_k)}\supseteq{}U_gC_gU_kC_k\supseteq{}U_gU_k\int_{\cD_k}f_{a_g^*a_g^{}}(\alpha_k(\chi))dE_\pi(\chi)\cdot{}C_k.
\end{gather}

From the properties of the polar decomposition and
the equality $\overline{\pi(a_g^*a_g^{})}=\int f_{a_g^*a_g^{}}dE_\pi$ we conclude
that $U_g^*U_g^{}=E_\pi(f_{a_g^{*}a_g^{}}^{-1}(0,+\infty))$. Similarly, $U_k^*U_k^{}=E_\pi(f_{a_k^{*}a_k^{}}^{-1}(0,+\infty)).$
Using Proposition \ref{prop_comm_rel_for_well_beh} $(ii)$ it follows that
\begin{gather}\nonumber
(U_gU_k)^*U_gU_k=U_k^*E_\pi(f_{a_g^{*}a_g^{}}^{-1}(0,+\infty))U_k=U_k^*U_kE_\pi(\alpha_{k^{-1}}(\cD_{k^{-1}}\cap f_{a_g^{*}a_g^{}}^{-1}(0,+\infty)))=\\
\label{eq_(U_gU_k)^*U_gU_k}=E_\pi(f_{a_k^{*}a_k^{}}^{-1}(0,+\infty))E_\pi(\alpha_{k^{-1}}(\cD_{k^{-1}}\cap f_{a_g^{*}a_g^{}}^{-1}(0,+\infty)))%
\end{gather}
is a projection. Hence $U_gU_k$ is a partial isometry. We denote by
$S_{gk}$ the closure of the operator
$\int_{\cD_k}f_{a_g^*a_g^{}}(\alpha_k(\chi))dE_\pi(\chi)\cdot C_k.$
From (\ref{eq_(U_gU_k)^*U_gU_k}) and the properties of the partial
action we conclude that the kernels of $U_gU_k$ and $S_{gk}$ are
equal. Since $S_{gk}$ is positive and its domain $\cD(S_{gk})$
contains $\cD(\pi),$ it follows from
(\ref{eq_polar_decompos_pi(a_ga_k)}) that the polar decomposition of
$\overline{\pi(a_ga_k)}$ is $U_gU_kS_{gk}$. Hence $U_gU_k$ belongs
to $\cU_\pi.$ By Proposition \ref{prop_comm_rel_for_well_beh}
$(ii),$ $\gA_0$ is closed under multiplication. That is, $\gA_0$ is
a unital $*$-algebra.

Since any $T\in\gA_0'$ commutes with $\cU_\pi$ and with the spectral
projections $E_\pi(\cdot),$ we have $T\in\pi(\cA)'.$ That is,
$\gA_0'\subseteq\pi(\cA)'$ and so $\gA_0''\supseteq\pi(\cA)''$ which
implies that $\gA_0''=\pi(\cA)''.$ Hence $\gA_0$ is dense in
$\pi(\cA)''$ in the strong operator topology. \hfill $\Box$

\begin{prop}\label{prop_pi_cyclic_iff_pi''_cyclic}
Suppose that $\pi$ is a well-behaved representation of
algebra $\cA$ such that the graph topology of $\pi$ is metrizable. Then $\pi$ is cyclic
if and only if the von Neumann algebra $\pi(\cA)''$ is cyclic.
\end{prop}
\noindent\textbf{Proof.} Suppose that $\varphi_0\in\cH_\pi$ is a
cyclic vector for $\pi.$ Let $\psi\in\cD(\pi)$ and $\varepsilon>0.$
Then there exists an element $a\in\cA$ such that
$\norm{\pi(a)\varphi_0-\psi}<\varepsilon.$ Clearly, $a$ is a finite
sum $a_1+a_2+\dots+a_k,$ where each $a_i$ belong to some vector
space $\cA_g,\ g\in G.$ Let $\overline{\pi(a_i)}=U_iC_i$ be the
polar decomposition of $\overline{\pi(a_i)}$. Since the operators
$U_i$ (by Lemma \ref{lemma_commut_polar_decompos}) and the spectral
projections $E_{C_i}(\cdot)$ of $C_i$ belong to $\pi(\cA)'',$ the
operators
$$
A_{i,r}:=U_i\int_{-r}^r\lambda dE_{C_i}(\lambda),\ r\in\dN,
$$
are in the von Neumann algebra $\pi(\cA)''.$ We choose $r\in\dN$ such that
$\norm{(A_{i,r}-\pi(a_i))\varphi_0}<\varepsilon/k$, $i=1,\dots,k$, and put
$A_r:=A_{1,r}+\dots+A_{k,r}.$ Then we have
\begin{gather*}
\norm{A_r\varphi_0-\psi}\leq\norm{(A_r-\pi(a))\varphi_0}+\norm{\pi(a)\varphi_0-\psi}\leq
\sum_{i=1}^k\norm{(A_{i,r}-\pi(a_i))\varphi}+\norm{\pi(a)\varphi_0-\psi}<2\varepsilon.
\end{gather*}
Since $A_r\in\pi(\cA)'',$ this shows that $\varphi_0$ is cyclic for $\pi(\cA)''.$

Conversely, suppose that $\varphi_0$ is a cyclic vector for the von
Neumann algebra $\pi(\cA)''.$ Let $P_0$ be the orthogonal projection
onto the closure of $\pi(\cB)''\varphi_0.$ Obviously,
$P_0\in\pi(\cB)'.$ Since $\Res_\cB\pi$ is self-adjoint by Definition
\ref{defn_well_beh}, $P_0\cH_\pi$ reduces $\Res_\cB\pi$ to a
self-adjoint subrepresentation $\rho$ (\cite{s1}, 8.3.11) which is
also integrable (\cite{s1}, 9.1.17). The graph topology of $\pi$ is
metrizable by assumption, so are the graph topologies of
$\Res_\cB\pi$ and $\rho$ by Lemma \ref{lemma_prelim}, $(i)$.
Therefore, a theorem of R.T. Powers (\cite{pow}, see \cite{s1},
9.2.1) applies and states that $\rho$ is cyclic, that is, there
exists a vector $\psi_0\in\cD(\rho)$ such that $\rho(\cB)\psi_0$ is
dense in $\cD(\rho)$ in the graph topology. In particular
$\overline{\rho(\cB)\psi_0}=P_0\cH_\pi=\overline{\pi(\cB)''\varphi_0}.$
Hence $\psi_0$ is also cyclic for the commutative von Neumann
algebra $\rho(\cB)''=P_0\pi(\cB)''P_0.$ Our aim is to show that
$\psi_0$ is cyclic for $\pi,$ that is, $\pi(\cA)\psi_0$ is dense in
$\cD(\pi)$ in the graph topology of $\pi.$

We first show that the subspace $\cH_0:=\pi(\cA)\psi_0$ is dense in
$\cH_\pi.$ Let $\gA_0$ be as in Lemma
\ref{lemma_commut_polar_decompos}. Since $\gA_0$ is dense in
$\pi(\cA)''$ in the strong operator topology, the vector $\varphi_0$
is also cyclic for $\gA_0.$ Let $U_g\in\cU_\pi$ and $a_g\in\cA_g,\
g\in G,$ be such that the polar decomposition of
$\overline{\pi(a_g)}$ is $U_gC_g.$ It suffices to show, that for any
Borel $\Delta_0\subseteq\cBp$ and $\varepsilon>0$ there exists
$b_1\in\cB$ such that
\begin{gather}\label{eq_aux_for_cyclic_1}
\norm{U_gE_\pi(\Delta_0)\varphi_0-\pi(a_gb_1)\psi_0}<\varepsilon.
\end{gather}
Let $b_0$ be such that
$\norm{\rho(b_0)\psi_0-E_\pi(\Delta_0)\varphi_0}<\varepsilon/3.$
Denote by $E_{C_g}$ the spectral measure on $\dR_+$ associated with
$C_g.$ Since $U_gE_{C_g}([0,+\infty))=U_gE_{C_g}((0,+\infty)),$ we
can choose $n$ such that
\begin{gather}\label{eq_aux_for_cyclic_2}
\norm{U_g(E_{C_g}([0,1/n])+E_{C_g}([n,+\infty)))\rho(b_0)\psi_0}<\varepsilon/3.
\end{gather}
Further, let $f$ be the function on $\dR$ defined by $f(x)=1/x$ if $
x\in(1/n,n)$ and $f(x)=0$ otherwise. Then the bounded operator
$f(C_g)$ is quasi-inverse to $C_g,$ that is, we have
$$Id_{\cH_\pi}=C_g f(C_g)+E_{C_g}([0,1/n])+E_{C_g}([n,+\infty)).$$

Since $\psi_0$ is strongly cyclic and $\overline{\pi(a_g^*a_g^{})}=C_g^2,$
there exists $b_1\in\cB$ such that
\begin{gather}\label{eq_aux_for_cyclic_3}
\norm{(1+C_g^2)(f(C_g)\rho(b_0)-\rho(b_1))\psi_0}<\varepsilon/3.
\end{gather}

Using (\ref{eq_aux_for_cyclic_2}) and (\ref{eq_aux_for_cyclic_3}) we
derive
\begin{gather*}
\norm{U_gE_\pi(\Delta_0)\varphi_0-\pi(a_gb_1)\psi_0}\leq\norm{U_g(E_\pi(\Delta_0)\varphi_0-\rho(b_0)\psi_0)}+\norm{U_g(\rho(b_0)-C_g\rho(b_1))\psi_0}\\
\leq\norm{U_g}\varepsilon/3+\norm{U_g\left(E_{C_g}([0,1/n])+E_{C_g}([n,+\infty))\right)\rho(b_0)\psi_0}+
\norm{U_g(C_gf(C_g)\rho(b_0)-C_g\rho(b_1))\psi_0}\\
\leq\varepsilon/3+\varepsilon/3+\norm{U_gC_g(1+C_g^2)^{-1}}\cdot\norm{(1+C_g^2)(f(C_g)\rho(b_0)-\rho(b_1))\psi_0}<\varepsilon.
\end{gather*}
Thus we have shown that $\cH_0$ is dense in $\cH_\pi.$

Let $\cD_0$ denote the closure of $\pi(\cA)\psi_0$ in the graph
topology of $\pi.$ We show that the representation
$\pi_0:=\pi\upharpoonright\cD_0$ of $\cA$ is self-adjoint. Since
$\rho$ is a restriction of $\Res_\cB\pi$, it is inducible. Let
$\cH_1$ denote the representation space of $\Ind\rho.$ Define a
linear operator $T:\cA\otimes\cD(\rho)\to\cD_0\subseteq\cD(\pi)$ by
$T(a\otimes\psi_0):=\pi(a)\psi_0.$ One easily checks that $T$ gives
rise to a unitary operator $\widetilde{T}$ of $\cH_1$ onto $\cH_0$
such that $ \widetilde{T}[a\otimes \psi_0]=\pi(a)\psi_0$ and that
$\widetilde{T}$ defines a unitary equivalence of representations
$\Ind\rho$ and $\pi_0.$ Since $\rho$ is cyclic and well-behaved,
$\Ind\rho$ is well-behaved by Proposition
\ref{prop_ind_of_cyc_well_beh_is_well_beh} and hence self-adjoint by
Lemma \ref{lemma_prelim}. Therefore, $\pi_0$ is self-adjoint. Since
$\cD(\pi_0)=\cD_0$ is dense in $\cH_\pi$ as shown in the preceding
paragraph, the $*$-representation $\pi$ of $\cA$ is an extension of
the self-adjoint representation $\pi_0$ acting on the same Hilbert
space $\cH_0$. By Corollary 8.3.12 in \cite{s1} this implies that
$\cD_0=\cD(\pi)$, that is, $\psi_0$ is a cyclic vector for $\pi$.
\hfill $\Box$

\begin{prop}\label{prop_well_beh_is_sum_of_cyclic}
Let $\pi$ be a well-behaved representation of $\cA$ on the Hilbert
space $\cH_\pi$ such that the graph topology of $\pi$ is metrizable.
Then $\pi$ can be decomposed into a direct orthogonal sum of cyclic
well-behaved representations.
\end{prop}
\noindent\textbf{Proof.} The identity representation of the von
Neumann algebra $\pi(\cA)''$ can be decomposed into a direct sum of
cyclic representations, i.e. there exists a decomposition
$\cH_\pi=\oplus_{i\in I}\cH_i$ such that the orthogonal projections
$P_i$ onto $\cH_i$ belong to $\pi(\cA)'$ and each von Neumann
algebra $P_i\pi(\cA)''$ is cyclic on $\cH_i$. By Proposition 8.3.11 in
\cite{s1} each representation $\pi_i:=\pi\upharpoonright
P_i\cD(\pi)$ is self-adjoint. It is straightforward to check that
$\pi=\oplus_{i\in I}\pi_i.$ Since $\pi$ is well-behaved, it follows
from Proposition \ref{prop_self_adj_subrep_is_well_beh} that
$\pi_i$, $i\in I$, is well-behaved. By Proposition
\ref{prop_pi_cyclic_iff_pi''_cyclic}, each representation $\pi_i$ is
cyclic. \hfill $\Box$\mn

Proposition \ref{prop_well_beh_is_sum_of_cyclic} combined with
Lemmas \ref{lemma_nonneg_sesq_form} and \ref{lemma_prelim} implies the
following
\begin{prop}
Let $H$ be a subgroup of $G$ and let $\rho$ be a well-behaved
representation of $\cA_H$ with metrizable graph topology. Then
$\rho$ is inducible to a $*$-representation of $\cA$ if and only
$\rho$ is $\cC$-positive, where $\cC:=\sum\cA^2\cap\cA_H$.
\end{prop}

\section{Well-behaved systems of imprimitivity}\label{sect_well_beh_sys_impr}%

In this section we shall prove an analogue of the Imprimitivity Theorem for well-behaved representations. A crucial step for this is to show that representations induced from well-behaved ones are again well-behaved. In the view of Proposition \ref{prop_ind_of_cyc_well_beh_is_well_beh} we assume for this section that $\cB$ is countably generated. We retain the notation from the previous section. Throughout $H$ denotes a subgroup of the group $G$.

\begin{defn}
A system of imprimitivity $(\pi,E)$ for $\cA$ over $G/H$ is called \textit{well-behaved} if
\begin{enumerate}
  \item[$(i)$] $\pi$ is a well-behaved representation of $\cA$,
  \item[$(ii)$] the projections $E$ and $E_\pi$ commute, that is,
$E(t)E_\pi(\Delta)=E_\pi(\Delta)E(t)$ for all $t\in G/H$ and all Borel subsets $\Delta$ of $\cBp$.
\end{enumerate}
\end{defn}
From Propositions \ref{prop_ind_of_cyc_well_beh_is_well_beh} and \ref{prop_well_beh_is_sum_of_cyclic} we obtain the following result.

\begin{prop}\label{prop_ind_of_well_beh_is_well_beh}
If $\rho$ is a well-behaved inducible representation of the
$*$-algebra $\cA_H$ with metrizable graph topology, then the induced
representation $\pi=\Ind_{\cA_H\uparrow\cA}(\rho)$ is a well-behaved
representation of the $*$-algebra $\cA.$
\end{prop}

The next proposition is an analogue of Proposition
\ref{prop_impr_induced}.

\begin{prop}
If $\rho$ is a well-behaved inducible $*$-representation of $\cA_H,$
then the system of imprimitivity induced by $\rho$ is non-degenerate
and well-behaved.
\end{prop}
\noindent\textbf{Proof.} Let $(\pi,E)$ be the system of
imprimitivity induced by $\rho$ and let $E_\pi(\cdot)$ be a
spectral measure associated with $\pi.$ It follows from
Proposition \ref{prop_impr_induced} that $(\pi,E)$ is
non-degenerate. By Proposition
\ref{prop_ind_of_well_beh_is_well_beh} the representation $\pi$ is
well-behaved. From the construction of $E(\cdot)$ (see Section
\ref{sect_grad_alg}) and relation (\ref{eq_action_ind_spec_measure})
it follows easily that $E(\cdot)$ and $E_\pi(\cdot)$ commute. \hfill
$\Box$

\begin{thm}(Imprimitivity Theorem for well-behaved representations)
Let $H$ be a subgroup of $G$ and let $(\pi,E)$ be a non-degenerate
well-behaved system of imprimitivity for $\cA$ over $G/H.$ Then
there exists a unique, up to unitary equivalence, inducible
well-behaved representation $\rho$ of $\cA_H$ such that $(\pi,E)$ is
unitarily equivalent to the system of imprimitivity induced by
$\rho.$
\end{thm}
\noindent\textbf{Proof.} Define $\rho$ as in the proof of the
Theorem \ref{thm_imprim_non_degenerate}. By Theorem
\ref{thm_imprim_non_degenerate} we only need to prove that $\rho$ is
well-behaved. Recall that the representation space $\cH_\rho$ is
defined as $\Ran E(H)$ and the domain $\cD(\rho)$ of $\rho$ is
$\cD(\pi)\cap\Ran E(H).$ For a Borel set $\Delta\subseteq\cBp$ put
$E_\rho(\Delta):=E_\pi(\Delta)E(H).$ Since $E_\pi(\cdot)$ commutes
with $E(\cdot),\ E_\rho$ is a well-defined spectral measure on
$\cBp$ whose values are projections in the Hilbert space $\Ran
E(H)=\cH_\rho.$ One easily checks that $\Res_\cB\rho$ is integrable
and defined by $E_\rho(\cdot).$

Let $a_h\in\cA_h,\ h\in H,\ v\in\cD(\rho),$ and let $\Delta\subseteq\cBp$
be a Borel set. Since $\pi(a_h)v=E(H)\pi(a_h)v$, we compute
$$
\rho(a_h)E_\rho(\Delta)v=\pi(a_h)E_\pi(\Delta)v=E_\pi(\Delta^h)\pi(a_h)v=E_\rho(\Delta^h)\rho(a_h)v.
$$
Hence $\rho$ is well-behaved. \hfill $\Box$

For the sake of completeness we formulate an analogue of Theorem
\ref{thm_impr_for_bounded_degen_system} for well-behaved
representations. Using the fact that well-behaved subrepresentations
have complements, the proof is similar to that of Theorem
\ref{thm_impr_for_bounded_degen_system}.
\begin{thm}\label{thm_impr_for_well_beh_degen_system}
Let $H$ be a subgroup of $G$ and let $(\pi,E)$ be a well-behaved
system of imprimitivity for $\cA$ over $G/H.$ Fix one element
$k_t\in G,\ t\in G/H,$ in each left coset from $G/H.$ Then for every
$t\in G/H$ there exists a well-behaved $*$-representation $\rho_t$
of $\cA_{k_t^{}Hk_t^{-1}}$ on a Hilbert space $\cH_t$ such that:
\begin{enumerate}
  \item[$(i)$] $\rho_t$ is inducible,
  \item[$(ii)$] $(\pi,E)$ is the direct sum of systems of imprimitivity $(\pi_t,E_t),\ t\in G/H,$ where
  $(\pi_t,E_t)$ is conjugated by the element $k_t$ to the system of imprimitivity induced by $\rho_t,\ t\in G/H.$
\end{enumerate}
\end{thm}

\begin{defn}
Let $\pi$ be a well-behaved representation of $\cA.$ We say that
\textit{$\pi$ is associated with an orbit $\Orb\chi,$} where
$\chi\in\cBp,$ if the spectral measure $E_\pi$ associated with $\pi$
is supported on the set $\Orb\chi.$
\end{defn}


The next theorem is a central result of the Mackey analysis (cf.
\cite{fd}, p. 1251 and p. 1284).

\begin{thm}\label{thm_main}
Assume that the group $G$ is countable. Let $\chi\in\cBp$ be a
character and let $H=\St\chi$ be its stabilizer group. Then the map
\begin{gather}\label{eq_map_inducing_A_H_to_A}
\rho\mapsto\Ind_{\cA_H\uparrow\cA}(\rho)=\pi
\end{gather}
is a bijection from the set of unitary equivalence classes of
inducible representations $\rho$ of $\cA_H$ for which
\begin{gather}\label{eq_cond_rep_A_H}
\Res_\cB\rho\ \mbox{corresponds to a multiple of the character}\ \chi
\end{gather}
onto the set of unitary classes of well-behaved representations
$\pi$ of $\cA$ associated with $\Orb\chi$. A $*$-representation
$\rho$ satisfying (\ref{eq_cond_rep_A_H}) is bounded and inducible.
Moreover, the von Neumann algebras $\rho(\cA_H)'$ and $\pi(\cA)'$
are isomorphic. In particular, $\pi$ is irreducible if and only if
$\rho$ is irreducible.
\end{thm}
\noindent\textbf{Proof.} Let $\pi$ be a well-behaved representation
of $\cA$ associated with $\Orb\chi,\ \chi\in\cBp.$ Since $G$ is
countable, the orbit $\Orb\chi$ is also countable. Therefore the
spectral measure $E_\pi$ is discrete. From the definition of $E_\pi$
it follows that $E_\pi(\set{\psi}),\ \psi\in\Orb\chi,$ is the
eigenspace of each operator $\pi(b),\ b\in\cB,$ corresponding to the
eigenvalue $\psi(b).$ Hence for all $\psi\in\Orb\chi$ the range
$\Ran E_\pi(\set{\psi})$ is contained in the domain of $\Res_\cB\pi$
which is equal to $\cD(\pi).$

Since $H$ is the stabilizer of $\chi,$ the projections
$E_\pi(\set{\chi}^{g_1})$ and $E_\pi(\set{\chi}^{g_2})$ are equal if
$g_1H=g_2H$ and for all $v\in\cD(\pi)$ we have
$$
\pi(a_g)E_\pi(\set{\chi}^k)v=E_\pi(\set{\chi}^{gk})\pi(a_g)v.
$$
(Note that if $\chi\in\cD_g,$ then $E_\pi(\set{\chi}^g)$ is equal to
$E_\pi(\set{\alpha_g(\chi)}),$ otherwise it is zero projection.)
Therefore, we can define a system of imprimitivity $E$ for $\cA$
over $G/H$ by putting $E(gH):=E_\pi(\set{\chi}^g).$

We show that
$(\pi,E)$ is non-degenerate. Let $g\in G$ be such that $\chi\in\cD_g$ and
let $e_{\chi^g}\in\Ran E(gH)$ be a non-zero vector. Since
$\chi^g\in\cD_{g^{-1}},$ there exists $a_{g^{-1}}\in\cA_{g^{-1}}$
such that $\chi^g(a_{g^{-1}}^*a_{g^{-1}}^{})>0.$ Since $e_{\chi^g}$
belongs to $\Ran E(gH)$ and $a_{g^{-1}}\in\cA_{g^{-1}},$ the vector
$\pi(a_{g^{-1}})e_{\chi^g}$ belongs to $\Ran E(H).$ Set
$e_\chi=(\chi^g(a_{g^{-1}}^*a_{g^{-1}}^{}))^{-1}\pi(a_{g^{-1}})e_{\chi^g}.$
Then, since $a_{g^{-1}}^*\in\cA_g$ and $e_{\chi^g}\in\Ran
E_\pi(\set{\chi^g})$, we obtain
$$
\pi(a_{g^{-1}}^*)e_\chi=(\chi^g(a_{g^{-1}}^*a_{g^{-1}}^{}))^{-1}\pi(a_{g^{-1}}^*a_{g^{-1}}^{})e_{\chi^g}=e_{\chi^g}.
$$
Thus, we have shown that the set $\set{\pi(a_g)e_\chi|a_g\in\cA_g,\
e_\chi\in\Ran E(H)}$ is equal to $\Ran E(gH),$ that is, $(\pi,E)$ is
non-degenerate. Since $E(H)$ is equal to $E_\pi(\set{\chi}),$
condition (\ref{eq_cond_rep_A_H}) is satisfied.

Conversely, let $\rho$ be a $*$-representation of $\cA_H$ satisfying
condition (\ref{eq_cond_rep_A_H}). Since $\rho(a_h^*a_h^{}),\
a_h\in\cA_h,\ h\in H,$ is a multiple of the identity, $\rho(a_h)$ is
bounded. Therefore each $\rho(a),\ a\in\cA$, is bounded, in
particular $\cD(\rho)=\cH_\rho.$ We will show later (see
Proposition \ref{prop_rep_A_H_is_inducible}) that every
representation $\rho$ satisfying (\ref{eq_cond_rep_A_H}) is positive
on the cone $\sum\cA^2.$ Since $\rho$ is bounded, it is a direct sum
of cyclic representations and hence inducible by Lemma
\ref{lemma_nonneg_sesq_form}. Proposition
\ref{prop_ind_of_cyc_well_beh_is_well_beh} together with Lemma
\ref{lemma_ind_of_sum} imply that
$\pi=\Ind_{\cA_H\uparrow\cA}\rho$ is well-behaved. Let $E_\pi$ be
the spectral measure associated with $\pi.$ The equality
(\ref{eq_action_ind_spec_measure}) implies that $E_\pi$ is supported
on $\Orb\chi$ which means that $\pi$ is associated with $\Orb\chi.$

It was shown in the proof of the Theorem
\ref{thm_imprim_non_degenerate} that the map
$$
\pi\mapsto\Res_{\cA_H}\pi\upharpoonright\Ran E(H)
$$
is the inverse of the map (\ref{eq_map_inducing_A_H_to_A}). Thus, we have proved
that the mapping (\ref{eq_map_inducing_A_H_to_A}) is indeed a
bijection.

Now we prove that $\rho(\cA_H)'=\pi(\cA)'.$ Let $T\in\rho(\cA_H)'.$ Define linear operator $\widetilde{T}$ on $\cA\otimes\cH_\rho$ by putting
\begin{gather}\label{eq_defn_oper_in_comm_of_ind_rep}
\widetilde{T}(a\otimes v)=a\otimes Tv,\ a\in\cA, v\in\cH_\rho.
\end{gather}
Let $c_H\in\cA_H.$ Then for arbitrary $a\in\cA$ and $ v\in\cH_\rho$ we have
$$
\widetilde{T}(ac_H\otimes v-a\otimes c_Hv)=ac_H\otimes Tv-a\otimes
Tc_Hv=ac_H\otimes Tv-a\otimes c_HTv,
$$
so $\widetilde{T}$ defines a linear operator on
$\cA\otimes_{\cA_H}\cH_\rho$ which is also denoted by
$\widetilde{T}.$

Let $a\in\cA,\ v\in\cH_\rho.$ We denote by $\norm{\cdot}_0$ the
seminorm $\langle\cdot,\cdot\rangle_0^{1/2}$. Since
$\rho$ is inducible, $S:=\rho(p_H(a^*a))$ is a positive operator on
$\cH_\rho$ commuting with $T.$ Hence $T$ commutes with $S^{1/2}$ and
we get
\begin{gather*}
\norm{\widetilde{T}(a\otimes v)}_0^2=\langle\widetilde{T}(a\otimes v),\widetilde{T}(a\otimes v)\rangle_0=
\langle\rho(p_H(a^*a))Tv,Tv\rangle=\langle S^{1/2}Tv,S^{1/2}Tv\rangle \\
=\langle TS^{1/2}v,TS^{1/2}v\rangle\leq\norm{T}^2\langle S^{1/2}v,S^{1/2}v\rangle=
\norm{T}^2 \langle\rho(p_H(a^*a))v,v\rangle=\norm{T}^2\norm{a\otimes v}_0^2.%
\end{gather*}

Let $\rho$ be a direct sum of cyclic representations $\rho_i$ with
cyclic vectors $v_i,\ i\in I.$ Take $\xi=\sum a_k\otimes
v_k\in\cA\otimes_{\cA_H}\cH_\rho,$ where $a_k\in\cA$ and $v_k$ are
distinct, hence pairwise orthogonal, cyclic vectors. Then the
vectors $a_k\otimes v_k$ are pairwise orthogonal with respect to
$\langle\cdot,\cdot\rangle_0$. Using the preceding inequality and
the latter fact we derive
\begin{gather*}
\norm{\widetilde{T}\xi}^2= \norm{\widetilde{T}(\sum_k a_k\otimes v_k)}_0^2
\leq(\sum_k\norm{T}\norm{a_k\otimes v_k}_0)^2
=\norm{T}^2\sum_k\langle a_k\otimes v_k,a_k\otimes
v_k\rangle_0\\ =\norm{T}^2\langle\sum_k a_k\otimes v_k,\sum_k a_k\otimes
v_k\rangle_0=\norm{T}^2\norm{\sum_k a_k\otimes v_k}_0^2 =\norm{T}^2 \norm{\xi}^2.
\end{gather*}
This shows that $\widetilde{T}$ gives rise to a bounded operator on
$\cH_\pi,$ which we denoted again by $\widetilde{T}.$ It is
straightforward to check that $\widetilde{T}$ commutes with all
operators $\pi(a),\ a\in\cA,$ and that the map
$\beta:T\mapsto\widetilde{T}$ is a $*$-homomorphism from
$\rho(\cA_H)'$ into $\pi(\cA)'$.

If $\widetilde{T}=0$, then in particular $\langle Tv,Tv\rangle
=\norm{\widetilde{T}(1\otimes v)}^2=0$ for all $v\in\cD(\rho)$ which
implies that $T=0$. That is, $\beta$ is injective.

We prove that $\beta$ is surjective. Let $S$ be an operator from
$\pi(\cA)'.$ Then $S\in\pi(\cB)'.$ Since the restrictions of
$\Res_\cB\pi$ to $\Ran E(gH)=\Ran E_\pi(\set{\chi}^g)$ are disjoint
representations for distinct cosets $gH\in G/H,$ $S$ commutes with
all operators $E(gH).$ In particular, $S_1:=S\upharpoonright\Ran
E(H)$ is a bounded operator on the Hilbert space $\Ran E(H)$ which
commutes with all operators $\pi(a)\upharpoonright\Ran E(H),$ where
$ a\in\cA_H.$ By the canonical isomorphism of $\cH_\rho$ and $\Ran
E(H)$, $S_1$ is a bounded operator on $\cH_\rho$. By construction we
have $S_1\in\rho(\cA_H)'.$ One easily verifies that $\beta(S_1)$ is
equal to $S$. This shows that $\beta$ is surjective. Summarizing the
preceding, we have proved that the mapping $\beta$ is an isomorphism
of von Neumann algebras $\rho(\cA_H)'$ and $\pi(\cA)'.$ \hfill
$\Box$\mn

\noindent\textbf{Remark.} Suppose that $\rho$ is an inducible
well-behaved representation of $\cA_H$. If condition
(\ref{eq_cond_rep_A_H}) does not hold, then the mapping
$\beta:T\mapsto\widetilde{T}$ of $\rho(\cA_H)'$ into $\pi(\cA)'$ is
not surjective in general.

\mn
We now derive an important corollary from the previous theorem.
\begin{prop}\label{prop_irr_iff_St_triv}
Let $\chi\in\cBp.$ Then the induced representation $\pi=\Ind\chi$
is irreducible if and only if its stabilizer group $\St\chi$ is
trivial.
\end{prop}
\noindent\textbf{Proof.} If the stabilizer $\St\chi$ is trivial,
then $\pi$ is irreducible by Theorem \ref{thm_main}.

Assume that the stabilizer group is not trivial. Then there exists
$h\in H{=}\St\chi$ such that $h\neq e.$ We choose an element $a_h\in
\cA_h$ such that $\chi(a_h^*a_h)=1.$ Using similar arguments as in
the proof of the Theorem \ref{thm_main}, one shows that there is a
linear operator $T_h$ on the $\cH_\pi$ defined by
$$
T_h([a_g\otimes 1])=[a_ga_h\otimes 1],\ a_g\in\cA_g,\ g\in G.
$$
For vectors $[a_1\otimes 1],[a_2\otimes 1]\in\cH_\pi$, where $
a_i\in\cA_{g_i},\ g_i\in G,i=1,2,$ we have
$$
\langle T_h[a_1\otimes 1],T_h[a_2\otimes 1]\rangle=
\langle[a_1a_h\otimes 1],[a_2a_h\otimes 1]\rangle=
\chi(p(a_h^*a_2^*a_1^{}a_h^{})).
$$
If $g_1\neq g_2.$, the latter is equal to
$0=\langle[a_1\otimes{}1],[a_2\otimes{}1]\rangle$. If $g_1=g_2$, then $a_2^*a_1\in\cB$ and hence
$$
\chi(p(a_h^*a_2^*a_1^{}a_h^{}))=\chi(a_h^*a_2^*a_1^{}a_h^{})=
\chi(a_2^*a_1^{})=\langle[a_1\otimes{}1],[a_2\otimes{}1]\rangle.
$$
This shows that $T_h$ is unitary. Since $T_h$ acts as a weighted
shift (see Proposition \ref{prop_space_ind_rep}), it is not a scalar
multiple of the identity. One easily verifies that $T_h$ commutes
with all representation operators. Since the commutant of $\pi$
contains a non-trivial unitary, $\pi$ is not irreducible. \hfill
$\Box$\mn

We now classify all representations of $\cA_H$ satisfying condition
(\ref{eq_cond_rep_A_H}). The result is the same as in the case when
$\cA$ is the group algebra $\dC[G]$ and $\cB$ is the group algebra
$\dC[N]$ of a commutative normal subgroup (see \cite{kir} and
\cite{fd}, pp. 1252-1258). That is, we establish a correspondence
between $*$-representations $\rho$ of $\cA_H$ satisfying
(\ref{eq_cond_rep_A_H}) and unitary projective representations of
$H.$

Let $\chi\in\cBp$ and let $H$ be the stabilizer group of $\chi.$
Take a representation $\rho$ satisfying (\ref{eq_cond_rep_A_H}).
Since $\chi^h$ is defined for all $h\in H,$ we can find elements
$a_h$ in each $\cA_h,\ h\in H,$ such that
$\chi(a_h^{}a_h^*)=\chi^h(a_h^{}a_h^*)=\chi(a_h^*a_h^{})\neq 0.$
From (\ref{eq_cond_rep_A_H}) it follows that for $h \in H$ the
operator
\begin{gather}\label{eq_defn_coc_rep}
\zeta(h):=\chi(a_h^*a_h^{})^{-1/2}\rho(a_h)
\end{gather}
is unitary and for any $b_h\in\cA_h$ the operator
$\rho(b_h^*a_h^{})$ is a scalar multiple of the identity, so
$\rho(a_h)$ differs from $\rho(b_h)$ by a scalar. Thus, the
operators $\zeta(h)$ define a unitary projective representation of
$H.$ Hence (see \cite{kir}) there exists a 2-cocycle $\tau:H\times
H\to\dT$
such that
\begin{gather}\label{eq_cocyc_rep}
\zeta(hk)=\tau(h,k)\zeta(h)\zeta(k),~h,k \in H.
\end{gather}
For $k\in H$ we have the equality
$\rho(a_k)^{-1}=\chi(a_k^*a_k^{})^{-1}\rho(a_k^*),$ in particular,
$\chi(a_k^*a_k^{})=\chi(a_k^{}a_k^*).$ Using this we calculate
\begin{gather*}
\zeta(hk)=\chi(a_{hk}^*a_{hk}^{})^{-1/2}\rho(a_{hk})=\chi(a_{hk}^*a_{hk}^{})^{-1/2}\rho(a_ha_k)\rho(a_ha_k)^{-1}\rho(a_{hk})=\\
=\chi(a_{hk}^*a_{hk}^{})^{-1/2}\chi(a_{h}^*a_{h}^{})^{1/2}\zeta(h)\chi(a_{k}^*a_{k}^{})^{1/2}\zeta(k)
\chi(a_h^*a_h^{})^{-1}\rho(a_h^*)\chi(a_k^*a_k^{})^{-1}\rho(a_k^*)\rho(a_{hk})=\\
=\chi(a_{hk}^*a_{hk}^{})^{-1/2}\chi(a_{h}^*a_{h}^{})^{-1/2}\chi(a_{k}^*a_{k}^{})^{-1/2}\chi(a_h^*a_k^*a_{hk})\zeta(h)\zeta(k).
\end{gather*}
Thus we have
\begin{gather}\label{eq_cocycle_from_chi}
\tau(h,k)=\chi(a_{hk}^*a_{hk}^{})^{-1/2}\chi(a_{h}^*a_{h}^{})^{-1/2}
\chi(a_{k}^*a_{k}^{})^{-1/2}\chi(a_h^*a_k^*a_{hk}),~h,k \in H.
\end{gather}
The mapping $\zeta$ satisfying (\ref{eq_cocyc_rep}) will be called
$\tau$-representation. Let $t$ be the element of the cohomology
group $Z^2(H,\dT)$ of $H$ with values in $\dT$ defined by the
cocycle $\tau.$ Analogously to the group case we call $t$ the
\textit{Mackey obstruction of} $\chi.$

Conversely, having a cocycle $\tau$ of the form (\ref{eq_cocycle_from_chi}) and
a $\tau$-representation $\zeta$ of $H$ it is straightforward to verify that
(\ref{eq_defn_coc_rep}) defines a $*$-representation $\rho$ of $\cA_H$
satisfying (\ref{eq_cond_rep_A_H}).

The proof of the following proposition is similar to the group case
(see \cite{fd}, pp. 1252-1258).

\begin{prop}\label{prop_obstr}
The Mackey obstruction $t$ of $\chi$ is trivial if and only if
$\chi$ can be extended to a character $\widetilde{\chi}$ of the
algebra $\cA_H.$ Equation (\ref{eq_defn_coc_rep}) defines a
one-to-one correspondence between unitary equivalence classes of
$\tau$-representations $\zeta$ of $H$ and unitary equivalence
classes of $*$-representations $\rho$ of $\cA_H$ satisfying
(\ref{eq_cond_rep_A_H}). Moreover, $\rho$ is irreducible if and only
if $\zeta$ is irreducible.
\end{prop}

We now show that condition (\ref{eq_cond_rep_A_H}) implies
$\sum\cA^2$-positivity.
\begin{prop}\label{prop_rep_A_H_is_inducible}
Let $\chi\in\cBp$ and let $H$ be its stabilizer. If $\rho$ is a
$*$-representation of $\cA_H$ satisfying condition
(\ref{eq_cond_rep_A_H}), then $\rho$ is nonnegative on the cone
$\sum\cA^2\cap\cA_H.$
\end{prop}
\noindent\textbf{Proof.} It suffices to show that for any $a\in\cA,\
\rho(p_H(a^*a))$ is a positive operator. It is enough to consider
the case when $a$ belongs to $\cA_{gH}$ for some $gH\in G/H,$ i.e.
$a=\sum_{h\in H}a_{gh},\ a_{gh}\in\cA_{gh}.$ Using that $H$ is the
stabilizer group of $\chi$, we get
\begin{gather*}
\chi(a_{gh}^*a_{gk}^{}a_{gk}^*a_{gh}^{})=\chi^{gh}(a_{gk}^{}a_{gk}^*)\chi(a_{gh}^*a_{gh}^{})=
\chi^{gk}(a_{gk}^{}a_{gk}^*)\chi(a_{gh}^*a_{gh}^{})=\chi(a_{gk}^*a_{gk}^{})\chi(a_{gh}^*a_{gh}^{}).
\end{gather*}
Using (\ref{eq_defn_coc_rep}) and the latter equality we calculate
\begin{gather*}
\rho(p_H(a^*a))=\rho(a^*a)=\sum_{k,h\in H}\rho(a_{gk}^*a_{gh}^{})=
\sum_{k,h\in H}\chi(a_{gh}^*a_{gk}^{}a_{gk}^*a_{gh}^{})^{1/2}\zeta(k^{-1}h)=\\
=\sum_{k,h\in{}H}\chi(a_{gk}^{}a_{gk}^*)^{1/2}\chi(a_{gh}^*a_{gh}^{})^{1/2}\zeta(k)^*\zeta(h)=
\left(\sum_{h\in H}\chi(a_{gh}^*a_{gh}^{})^{1/2}\zeta(h)\right)^*\sum_{h\in H}\chi(a_{gh}^*a_{gh}^{})^{1/2}\zeta(h),%
\end{gather*}
which implies that $\rho(p_H(a^*a))$ is positive. \hfill$\Box$\mn

Next we want to associate well-behaved irreducible representations
with orbits. Under some technical assumption this aim will be
achieved by Proposition \ref{prop_irrorbit} below. For this some
preparations are necessary.


\begin{defn}
A Borel subset $\Delta$ of $\cBp$ is called \textit{invariant}
under the partial action of $G$ if $\Delta^g\subseteq\Delta$ for
every $g\in G.$ A spectral measure $E$ on $\cBp$ is called
\textit{ergodic} under the partial action of $G$ on $\cBp$ if for
every invariant Borel subset $\Delta$ of $\cBp$ either $E(\Delta)$
or $E(\cBp\backslash\Delta)$ is zero.
\end{defn}

\begin{lemma}\label{prop_irrep_ergodic}
Let $\pi$ be a well-behaved irreducible representation of the
$*$-algebra $\cA$ and let $E_\pi$ be an associated spectral
measure. Then $E_\pi$ is ergodic.
\end{lemma}
\noindent\textbf{Proof.} Let $\Delta$ be a Borel subset of $\cBp$
which invariant under the partial action of $G.$ From Proposition
\ref{prop_comm_rel_for_well_beh}$(i),$ it follows that
$E_\pi(\Delta)$ is a projection commuting with $\pi(\cA_g)$ for all
$g\in G$ and hence with $\pi(\cA).$ Since $\pi$ is irreducible,
$E_\pi(\Delta)$ is trivial, i.e. $E_\pi(\Delta)=0$ or
$E_\pi(\Delta)=I$. \hfill $\Box$\mn

The following concepts are taken from the paper \cite{eff}.

We shall say that a measurable space $(Y,\gB)$ is \textit{countably
separated} if there exists a countable subfamily $\gB_0$ of $\gB$
such that for any two points in $Y$ there exists a member of $\gB_0$
containing one point but not the other. A measurable subset
$\Gamma\subseteq Y$ is said to be \textit{countably separated} if
$(\Gamma,\gB_\Gamma)$ is countably separated, where $\gB_\Gamma$ is
the induced Borel structure.

A subset $\Gamma\subseteq\cBp$ is called a \textit{section of the
partial action} of $G$ on $\cBp$ if it contains precisely one point
from each orbit. Recall that a (spectral) measure is called an
\textit{atom} if it attains only two values. An atom is called
\textit{trivial} if it is supported at a single point.

The proof of the following simple lemma is borrowed from the proof
of Theorem 2.6 in \cite{eff}.

\begin{lemma}\label{lemma_atom_on_coun_sep_is_triv}
Let $E$ be a spectral measure on a countably separated measurable
space $(X,\gB)$. If $E$ is an atom, then it is trivial.
\end{lemma}
\noindent\textbf{Proof.} Let $\{B_k; k \in \dN\}$ be a countable
family of Borel subsets of $X$ which separates the points of $X$ and
is closed under taking complements. Let $B_{k_n},\ n\in\dN,$ be
those sets with $E(B_{k_n})=I$ and put $B=\cap_{n\in\dN}B_{k_n}.$
Then we have $ E(B_{k_1}\cap\dots\cap B_{k_n})=E(B_{k_1})\dots
E(B_{k_n})=I$ which implies that $E(B)=I$ and $B\neq\emptyset.$

Assume to the contrary that there exist distinct points $p$ and $q$
in $B$. Then there exists $j\in \dN$ such that $p\in B_j$ and
$q\notin B_j.$ Due to the latter relation, we have
$B_j\notin\set{B_{i_n}}$ and $X\backslash B_j\notin\set{B_{i_n}}$
which implies that $E(B_j)$ and $E(X\backslash B_j)$ are zero. Hence
$E(X)=0$ which is a contradiction. \hfill $\Box$\mn


\begin{prop}\label{prop_irrorbit}
Let $G$ be a countable group. Suppose that the partial action of $G$ on
$\cBp$ possesses a measurable countably separated section $\Gamma.$
Then every ergodic spectral measure $E$ on $\cBp$ is supported on
a single orbit. In particular, each irreducible well-behaved
representation of $\cA$ is associated with an orbit.
\end{prop}
\noindent\textbf{Proof.} We first show that the spectral measure $E$
restricted to $\Gamma$ is either zero or an atom. Suppose that $E$
restricted to $\Gamma$ is non-zero. Assume to the contrary that $E$
restricted to $\Gamma$ is not an atom. Then $\Gamma$ is a disjoint
union of two Borel sets $\Gamma_1$ and $\Gamma_2$ such that
$E(\Gamma_1)\neq 0$ and $E(\Gamma_2)\neq 0$. By Proposition
\ref{prop_alpha_g_maps_Borel_to_Borel}, the sets
$\Omega_i=\cup_{g\in G}\Gamma_i^g,\ i=1,2,$ are Borel. The
properties of the partial action imply that the sets $\Omega_i$ are
invariant and both projections $E(\Omega_i)$ are non-zero which is a
contradiction. Thus, $E$ restricted to $\Gamma$ is an atom.

Since $\Gamma$ is countably separated, Proposition
\ref{prop_alpha_g_maps_Borel_to_Borel} implies that all
$\Gamma^g,\ g\in G,$ are countably separated. Since $\cBp$ is the
union of sets $\Gamma^g,$ it follows from Lemma
\ref{lemma_atom_on_coun_sep_is_triv} that there exist points
$\chi_k\in\Gamma^k,\ k\in I\subseteq G,$ such that
$E({\chi_k})\neq 0$ for all $k \in I$ and $E$ is supported on the
(at most countable) set $\set{\chi_k}_{k\in I}.$ Since the set
$\Orb\chi_k$ is invariant and $E(\Orb\chi_k)\neq 0$ for all $k$,
the ergodicity of $E$ implies that all $\chi_k$ belong to a single
orbit. \hfill $\Box$
%

\section{Example: Enveloping algebras of some complex Lie
algebras}\label{sect_env_alg}

In this section we illustrate the concepts of the previous sections
on three examples: enveloping algebras $\cU(su(2)),\ \cU(su(1,1))$
and $\cU(Vir),$ where $Vir$ denotes the Virasoro algebra
\cite{chari},\cite{fqs}. Is is easily checked that in these cases condition (\ref{eq_cond_compatible}) is satisfied and the space $\cBp$ is locally compact, so the theory developed in the preceding sections applies.

First let $\gog$ be one of the real Lie algebras $su(2)$ or
$su(1,1)$ and let $\gog_{\dC}$ be its complexification. Then
$\gog_\dC=sl(2,\dC)$ has a vector space basis $\set{E,F,H}$ with
commutation relations
\begin{gather}\label{eq_cr_sl_2}
[H,E]=2E,\ [H,F]=-2F,\ [E,F]=H.
\end{gather}
From (\ref{eq_cr_sl_2}) it follows that in the complex universal
enveloping algebra $\cU(\gog)$ we have
\begin{gather}
E q(H)=q(H-2)E,\ F q(H)=q(H+2)F\label{eq_cr_Eq(H)_Fq(H)} \\ %
HE^n=E^n(H+2n),\ FE^n=E^{n-1}(EF-n(H+n-1)),\ n\in\dN,\label{eq_cr_HE^n_FE^n}\\ %
HF^n=F^n(H-2n),\ EF^n=F^{n-1}(FE+n(H-n+1)),\ n\in\dN. %
\end{gather}
for each polynomial $q\in\dC[x]$ and that the Casimir element
$$C:=2(EF+FE)+H^2=4FE+H(H+2)=4EF+H(H-2)
$$
belongs to the center of $\cU(\gog).$

The complex unital algebra $\cU(\gog)$ becomes a $*$-algebra with
involution determined by $x^*=-x$ for $x\in \gog.$ In terms of the
generators $\set{E,F,H}$ of the algebra $\cU(\gog)$ this means that
\begin{gather}
E^*=F,\ H^*=H\ \mbox{for}\ \gog=su(2),\label{eq_involution_su2}\\
E^*=-F,\ H^*=H\ \mbox{for}\ \gog=su(1,1)\label{eq_involution_su11}.
\end{gather}

Using the commutation relation (\ref{eq_cr_sl_2}) it follows by
induction that
$$
\cU(\gog)_0:=\Lin\set{E^lF^lH^k;~k,l\in\dN_0}=
\Lin\set{(EF)^lH^k;~k,l\in\dN_0}=\Lin\set{C^lH^k;~k,l\in\dN_0}.
$$
In particular, $\cB:=\cU(\gog)_0$ is a commutative unital
$*$-subalgebra of $\cA=\cU(\gog).$ For $n\in\dN_0,$ let
$$\cA_n=E^n\cB=\Lin\set{E^{n+l}F^lH^k; k,l\in\dN_0},
\cA_{-n}=F^n\cB=\Lin\set{E^lF^{n+l}H^k; k,l\in\dN_0}.$$ By the
Poincare-Birkhoff-Witt theorem, $\set{E^iF^jH^l;\ i,j,l\in\dN_0}$ is
a vector space basis of $\cU(\gog).$ From this fact and the
definitions (\ref{eq_involution_su2}) and (\ref{eq_involution_su11})
of the involution we derive that
\begin{gather}\label{eq_grading}
\cA=\bigoplus_{n\in\dZ}\cA_n
\end{gather}
is a $\dZ$-graded $*$-algebra. Let $p:\cA\to\cB$ be the canonical conditional expectation (see Proposition \ref{prop_gradexp}). In both cases $\gog=su(2)$ and $\gog=su(1,1)$ the conditional expectation $p$ is not strong, because we have $E^*E\in\sum\cA^2\cap \cB,$ but $E^*E\notin\sum\cB^2.$\mn

\noindent\textbf{Remarks 1.} The $\dZ$-graded $*$-algebra (\ref{eq_grading}) is the special case $\gog=sl(2,\dC)$ of Example
\ref{semienv}. In this case, $Q=\dZ$ and $\cB=\cU(\gog)_0$ is just the commutant of the element $H$ in the algebra $\cU(\gog)$. Note
that $sl(2,\dC)$ is the only simple Lie algebra $\gog$ for which $\cB=\cU(\gog)_0$ is
commutative.\\

\noindent\textbf{2.} For the real Lie algebra $\gog=sl(2,\dR)$ the involution of the enveloping algebra $\cU(\gog)$ is given by $E^*=E,\ F^*=F,\
H^*=-H$. In this case the decomposition (\ref{eq_grading}) remains valid and shows that $\cU(\gog)$ is a $\dZ$-graded algebra. But
since $(\cU(\gog)_n)^*=\cU(\gog)_n$ for $n\in \dZ,\ \cU(\gog)=\oplus_n \cU(\gog)_n$ is not a $\dZ$-graded $*$-algebra.

We derive three simple lemmas which will be needed below.
\begin{lemma}\label{lemma_chi_in_cBp_inequalities}
Let $\gog$ be one of the real Lie algebras $su(2)$ or $su(1,1).$ A character $\chi\in\cBd$ belongs to $\cBp$ if and only $
\chi(F^{*k}F^k)\geq 0$ and $\chi(E^{*k}E^k)\geq 0$ for all $k\in\dN$.
\end{lemma}
\noindent\textbf{Proof.} Recall that $\chi\in\cBp$ if and only if $\chi(b)\geq 0$ for all $b\in\cA^2\cap\cB.$ Hence the necessity of
the condition is obvious. We prove that it is also sufficient. By Corollary \ref{corconeb}, it suffices to show $\chi(a_n^*a_n^{})\geq
0$ for all homogeneous elements $a_n\in\cA_n,\ n\in \dZ.$

Let $n\in\dN_0$ and take $a_n\in\cA_n.$ By the definition of $\cA_n$ we have $a_n=E^nb$ for some $b\in\cB$. Since $\chi(E^{*n}E^n)\geq 0$
by assumption, $\chi(a_n^*a_n^{})=\chi(b^*E^{*n}E^nb)=\chi(E^{*n}E^n)\chi(b^*b)\geq 0$. Similarly, for $n<0$ the inequality $\chi(F^{*n}F^n)\geq 0$
implies that $\chi(a_n^*a_n^{})\geq 0$ for all $a_n\in\cA_n.$ \hfill $\Box$

\begin{lemma}\label{lemma_EnFn} For $n \in \dN$ we have
\begin{gather}\label{equ2}
E^nF^n=EF(EF+H{-}2)(EF+H{-}2+H{-}4){\cdots}(EF+H{-}2+{\cdots}+H{-}2(n-1)),\\
\nonumber
F^nE^n=(EF-H-(H{+}2)-{\dots}-(H{+}2(n-1))){\cdots}(EF-H-(H{+}2))(EF-H)\\
=FE(FE-(H{+}2)){\cdots}(FE-(H{+}2)-{\cdots}-(H{+}2(n{-}1)))\label{equ3}
\end{gather}
\end{lemma}
\noindent\textbf{Proof.} We prove the first equality (\ref{equ2}) by induction on $n.$ The two equalities concerning $F^nE^n$ are
verified in a similar manner. Using the commutation relation (\ref{eq_cr_sl_2}) we compute
\begin{gather*}
E^{n+1}F^{n+1}=E^n(FE+H)F^n=E^nFEF^n+(H-2n)E^nF^n=\\
=E^{n-1}(FE+H)EF^n+(H-2n)E^nF^n=\\
=E^{n-1}FE^2F^n+(H-2(n-1))E^nF^n+(H-2n)E^nF^n=\dots\\
\dots=(EF+H-2+\dots+(H-2n))E^nF^n.
\end{gather*}
Inserting the induction hypothesis (\ref{equ2}) for $n$ and remembering that all elements $E^kF^k$ and $H^l$ mutually commute, we obtain (\ref{equ2}) for $n+1$. \hfill $\Box$

\begin{lemma}\label{bpolynom}
$\cB\equiv\cU(\gog)_0=\dC[EF,H]=\dC[C,H].$
\end{lemma}
\noindent\textbf{Proof.} Since the elements $EF$ and $H$ of $\cU(\gog)$ commute, there is an algebra homomorphism
$\sigma:\dC[x_1,x_2]\to \cU(\gog)$ given by $\sigma(x_1)=EF$ and $\sigma(x_2)=H$. From the Poincare-Birkhoff-Witt theorem we derive
easily that $\sigma$ is injective which gives $\cU(\gog)_0=\dC[EF,H].$ Clearly, we have also $\dC[EF,H]=\dC[C,H].$
\hfill $\Box$\mn

Lemma \ref{bpolynom} implies that the map $\cBd\ni\chi\mapsto(\chi(C),\chi(H))\in\dR^2$ is bijective. Denote by $\chi_{st}\in\cBd,\ s,t\in\dR$ a character such that
\begin{gather}\label{eq_defn_chi_st}
\chi_{st}(C)=s,\ \chi_{st}(H)=t.
\end{gather}

Propositions \ref{prop_cBp_su_2} and \ref{prop_cBp_su_11} below describe the set of parameters $s,t\in\dR$ for which
$\chi_{st}\in\cBp$ in the cases $\gog=su(2)$ and $\gog=su(1,1),$ respectively.

\begin{prop}\label{prop_p_action_sl_2}
Let $\gog$ be one of the real Lie algebras $su(2)$ or $su(1,1)$. If a character $\chi_{st}$ belongs to $\cBp$ and if $\chi_{st}^n$ is defined
for $n\in\dZ$, then we have
\begin{gather}\label{eq_p_action_cBp_env_alg}
\chi_{st}^n=\chi_{s,t+2n}.
\end{gather}
\end{prop}
\noindent\textbf{Proof.} For $n=0$ the proof is trivial. Assume that
$n>0.$ In the case $n<0$ the proof is similar. Since $\chi_{st}^n$
is defined, $\chi_{st}(E^{*n}E^n)>0.$ We compute
$$
\chi_{st}^n(H)=\frac{\chi_{st}(F^nHE^n)}{\chi_{st}(F^nE^n)}
=\frac{\chi_{st}(F^nE^n(H+2n))}{\chi_{st}(F_{}^nE_{}^n)}=\chi_{st}(H+2n)=t+2n=\chi_{s,t+2n}(H).
$$
Since $C$ belongs to the center of $\cA$, we have
$\chi_{st}^n(C)=\chi_{st}^{}(C).$ By the definition of $\chi_{st}$ we
obtain (\ref{eq_p_action_cBp_env_alg}). \hfill $\Box$\mn

\subsection{The case $\gog=su(2)$} In this subsection we let
$\cA=\cU(su(2))$ and $\cB=\cA_0=\dC[EF,H]=\dC[C,H].$ The next
proposition describes the set $\cBp.$
\begin{prop}\label{prop_cBp_su_2}
A character $\chi_{st}$ defined by (\ref{eq_defn_chi_st}) belongs to
$\cBp$ if and only if $t\in\dZ$ and $s=(t+2n)(t+2n+2)$ for some
$n\in\dN_0$ such that $n+t\geq 0.$
\end{prop}
\noindent\textbf{Proof.} Since $E^{*n}=F^n,$ Lemmas
\ref{lemma_chi_in_cBp_inequalities} and \ref{lemma_EnFn} imply that
$\chi$ belongs to $\cBp$ if and only if the following inequalities
are fulfilled for arbitrary $k\in\dN$:
\begin{gather}
\chi(E^kF^k)\equiv\chi(EF)\chi(EF+H-2)\dots\chi(EF+H-2+\dots+H-2k)\geq 0,\label{eq_ineq_def_B+_su2_1}\\
\chi(F^kE^k)\equiv\chi(EF-H)\chi(EF-H-(H+2))\dots\chi(EF-H-\dots-(H+2k))\geq 0.\label{eq_ineq_def_B+_su2_2}%
\end{gather}

We claim that for every $\chi\in\cBp$ there exist $m,n\in\dN_0$ such
that
\begin{gather}\label{eq_equa_chi_mn}
\chi(EF+m(H-(m+1)))=0,\ \chi(EF-(n+1)(H+n))=0.
\end{gather}

Assume to the contrary that $\chi(EF+k(H-(k+1)))\neq 0$ for all
$k\in\dN_0.$ It follows from (\ref{eq_ineq_def_B+_su2_1}) that
$\chi$ is positive on all factors in (\ref{eq_ineq_def_B+_su2_1}),
that is,
$$\chi(EF+H-2+\dots+H-2k)=\chi(EF+k(H-(k+1)))=\chi(EF)+k(\chi(H)-(k+1)))>0
$$
for all $k\in\dN_0$ which is a contradiction. Hence
$\chi(EF+m(H-(m+1)))=0$ for some $m\in\dN.$ In the same way one
proves the second equality in (\ref{eq_equa_chi_mn}).

The solution of the system of equations (\ref{eq_equa_chi_mn}) is
\begin{gather}\label{eq_defn_chi_mn}
\chi(EF)=m(n+1),\ \chi(H)=m-n.
\end{gather}
It is easy to verify that for all $m,n\in\dN_0$ the characters $\chi$ defined by
(\ref{eq_defn_chi_mn}) satisfy both inequalities
(\ref{eq_ineq_def_B+_su2_1}) and (\ref{eq_ineq_def_B+_su2_2}).

Putting $t=m-n$ in (\ref{eq_defn_chi_mn}) we get
\begin{gather*}
\chi(C)=4\chi(EF)+\chi(H^2-2H)=4m(n+1)+(m-n)^2-2m+2n=\\
=(m+n)(m+n+2)=(t+2n)(t+2n+2),
\end{gather*}
i.e. $\chi=\chi_{st}$ where $t=m-n\in\dZ$ and $s=(t+2n)(t+2n+2).$
Clearly, we have $m,n\in\dN_0$ if and only if $t\in\dZ,\
n+t\geq 0.$ \hfill $\Box$\mn

We denote by $\psi_n,\ n{\in}\dN_0,$ the character
$\chi_{n(n+2),-n}\in\cBp$ and by $\Gamma$ the subset
$\set{\psi_n,\ n\in\dN_0}$ of $\cBp.$ By Propositions
\ref{prop_p_action_sl_2} and \ref{prop_cBp_su_2}, each orbit
under the partial action of $\dZ$ on $\cBp$ contains precisely one
of the characters from $\Gamma,$ i.e. $\Gamma$ is a section of the
partial action of $\dZ$ on $\cBp.$

\begin{prop}
The representations $\Ind\chi,\ \chi\in\Gamma,$ are pairwise
non-equivalent and irreducible. Each irreducible well-behaved
representation of $\cA$ is unitarily equivalent to $\Ind\chi$ for
some $\chi\in\Gamma.$ A $*$-representation $\pi$ of $\cA=\cU(su(2))$
is well-behaved (in the sense of Definition \ref{defn_well_beh}) if
and only if $\pi$ is integrable (that is, $\pi{=}dU$ for some
unitary representation $U$ of the Lie group $SU(2).)$
\end{prop}
\noindent\textbf{Proof.} Clearly, the bijection
$\chi_{st}\mapsto(s,t)$ of the space $\cBd$ onto $\dR^2$ (by Lemma
\ref{bpolynom}) is a homeomorphism. Hence Proposition
\ref{prop_cBp_su_2} implies that $\cBp$ is a discrete space. It
follows from the formulas for the partial action of $\dZ$ that
$\Gamma$ is a Borel section. By Proposition \ref{prop_irrorbit} all
irreducible well-behaved representations are associated with orbits.
Therefore, by Theorem \ref{thm_main} we have that $\Ind\chi,\
\chi\in\Gamma,$ are up to unitary equivalence all irreducible
well-behaved representations. It follows from Proposition
\ref{prop_cBp_su_2} that $\Orb\psi_n,\ n\in\dN_0$ consists of $n+1$
elements, and Proposition \ref{prop_space_ind_rep} implies that
$\Ind\psi_n,\ n\in\dN_0$ has dimension $n+1.$ The latter implies in
particular that each representation $\Ind\chi,\ \chi\in\Gamma$ is
integrable.

Let $\pi$ be a well-behaved representation of $\cA$ and let $E_\pi$
be the associated spectral measure on $\cBp.$ Denote by $\rho$ the
restriction of $\Res_\cB\pi$ to $\Ran(E_\pi(\Gamma)).$ It is easily
checked that $\pi$ is unitarily equivalent to $\Ind\rho.$ Since
$\cBp$ is discrete, $\rho$ is equivalent to a direct sum of
characters $\chi\in\Gamma$ (taken with multiplicities), so that
$\pi$ is equivalent to a direct sum of representations $\Ind\chi,\
\chi\in\Gamma.$ Because $\Ind\chi$ is integrable as shown in the
preceding paragraph, $\pi$ is integrable.

Conversely, if $\pi$ is an integrable representation, $\pi$ is a
direct sum of integrable irreducible representations $\pi_i$. Since
each representation $\pi_i$ is finite dimensional and hence
well-behaved by Proposition \ref{prop_well_behav_bounded}, $\pi$ is
well-behaved.\hfill$\Box$\mn

It is well-known that for each $n\in\dN_0$ the spin $\frac{n}{2}$
representation is the unique (up to unitary equivalence) irreducible
$(n{+}1)$-dimensional $*$-representation of $\cA{=}\cU(su(2))$.
Since the $*$-representation $\Ind\psi_n$ of $\cA$ is irreducible
and of dimension $n{+}1$, $\Ind\psi_n$ is equivalent to the spin
$\frac{n}{2}$ representation. We want to establish this equivalence
by explicit formulas.

Recall from Proposition \ref{prop_space_ind_rep}, $(i)$ that the
vectors
$$
\set{\frac{[E^k\otimes 1]}{\norm{[E^k\otimes 1]}},\ k=0,1\dots n}
$$
form an orthonormal base of the representation space of
$\Ind\psi_n.$ By definition of $\psi_n$ we have $\psi_n(H)=-n$ and
$\psi_n(EF)=\frac{1}4\psi_n(C-H^2+2H)=0.$ Using Lemma
(\ref{lemma_EnFn}) we compute
\begin{gather*}
\norm{[E^k\otimes
1]}^2=\psi_n(F^kE^k)=\psi_n((EF-H)(EF-2(H+1))\dots(EF-k(H+k-1)))=\\
=n(2(n-1))\dots(k(n-k+1))=\frac{k!\cdot n!}{(n-k)!},\ k=0,1\dots,n.
\end{gather*}
Putting $l=\frac{n}{2},\ \pi_l:=\Ind\psi_n$ and
$$
e_m:=\frac{[E^{l+m}\otimes 1]}{\norm{[E^{l+m}\otimes
1]}}=\sqrt{\frac{(l-m)!}{(2l)!(l+m)!}}\ [E^{l+m}\otimes 1],\
m=-l,l+1,\dots, l,
$$
we calculate
\begin{gather*}
\pi_l(E)e_m=\frac{[E^{l+m+1}\otimes 1]}{\norm{[E^{l+m}\otimes
1]}}=\frac{\norm{[E^{l+m+1}\otimes 1]}}{\norm{[E^{l+m}\otimes
1]}}e_{m+1}=\sqrt{\frac{(2l)!(l+m+1)!}{(l-m-1)!}}\sqrt{\frac{(l-m)!}{(2l)!(l+m)!}}e_{m+1}=\\
=\sqrt{(l-m)(l+m+1)}e_{m+1},\ m=-l,l+1,\dots, l.
\end{gather*}
In the same manner we derive
$$
\pi_l(F)e_m=\sqrt{(l-m+1)(l+m)}e_{m-1},\ \pi_l(H)e_m=2me_m,\
m=-l,l+1,\dots, l.
$$
These are the formulas for the actions of $E,F,H$ in the spin $l$
representation of $\cU(su(2)).$\mn

We now show that the representations $\pi_l$ can be also induced
from the $*$-subalgebra $\cC=\dC[H].$ Let $p_3=p_2\circ p_1,$ where
$p_1$ is the canonical conditional expectation $p_1:\cA\to\cB$ and
$p_2:\cB\to\cC$ is conditional expectation defined by
$p_2((EF)^k)=0, p_2(H^k)=H^k,\ k\in\dN.$ Using Lemma
\ref{lemma_EnFn} we obtain
\begin{gather*}
p_3(\sum\cA^2)=\sum\cC^2-H\sum\cC^2+H(H+(H+2))\sum\cC^2\\
-H(H+(H+2))(H+(H+2)+(H+4))\sum\cC^2+\dots=\\
=\sum\cC^2-H\sum\cC^2+H(H+1)\sum\cC^2-H(H+1)(H+2)\sum\cC^2+\dots+\\
+(-1)^kH(H+1)(H+2)\dots(H+k-1)\sum\cC^2+\dots.
\end{gather*}
Obviously, $p_3$ is a $(\sum\cA^2,p_3(\sum\cA^2))$-conditional
expectation. It is easy to check that $\sum\cA^2\cap\dC[H]=\sum\cC^2.$
Since $p_3(\sum\cA^2)$ is strictly larger than $\sum\cC^2,\ p_3$ is
not a conditional expectation according to Definition
\ref{defn_cond_exp}. In particular we have seen that the composition of
two conditional expectations is not
a conditional expectation in general.

It is clear from the preceding formulas that the set of characters
on $\dC[H]$ which are non-negative on the cone $p_3(\sum\cA^2)$ and
hence inducible via $p_3$ is the set $\set{\chi_k,\ k\in\dN_0}.$
Note that $\chi_k(H)=-k.$ It is not difficult to compute that the
corresponding induced representation $\Ind\chi_{2l},\
l\in\frac{1}2\dN_0,$ is unitarily equivalent to $\pi_l.$

\subsection{The case $\gog=su(1,1)$} In this subsection let $\cA=\cU(su(1,1))$ and $\cB=\cA_0=\dC[EF,H]=\dC[C,H].$

We denote by $\chi_{st}\in\cBd$ the characters determined by
(\ref{eq_defn_chi_st}). It is convenient to introduce the following
subsets of $\cBd:$
\begin{gather*}
X_{00}=\set{\chi_{00}},\\
X_{1k}=\set{\chi_{st}|2k\leq t<2k+2,\ -\infty<s<(t-2k)(t-2(k+1))},\ k\in\dZ,\\ %
X_{2k}=\set{\chi_{st}|2k<t<2k+2,\ s=(t-2k)(t-2(k+1))},\ k\in\dZ,\\ %
X_{3k}=\set{\chi_{st}|t\geq 2k+2,\ s=(t-2k)(t-2(k+1))},\ k\in\dN_0,\\ %
X_{4k}=\set{\chi_{st}|t\leq 2k,\ s=(t-2k)(t-2(k+1))},\ k\in\dZ\backslash\dN_0.%
\end{gather*}

The following two propositions describe the set $\cBp$ and the
partial action of $\dZ$ on it.
\begin{prop}\label{prop_cBp_su_11} The set $\cBp$ is equal to the disjoint union
$$X_{00}\cup\bigcup_{k\in\dZ}X_{1k}\cup\bigcup_{k\in\dZ}
X_{2k}\cup\bigcup_{k\in\dN_0}X_{3k}\cup\bigcup_{k\in\dZ\backslash\dN_0}X_{4k}.
$$
\end{prop}
\noindent\textbf{Proof.} The equality $E^{*n}=(-1)^nF^n$ and Lemmas
\ref{lemma_chi_in_cBp_inequalities} and \ref{lemma_EnFn} imply that
a character $\chi\in\cBd$ belongs to $\cBp$ if and only if the
following inequalities hold:
\begin{gather}
(-1)^{k}\chi(EF(EF+H-2)\dots(EF+H-2+H-4+\dots+H-2(k-1)))\geq 0,\ k\in\dN,\label{eq_ineq_def_B+_su11_1}\\
(-1)^k\chi((EF-H)(EF-H-(H+2))\cdot\dots\label{eq_ineq_def_B+_su11_2}\\
\dots(EF-H-(H+2)-\dots-(H+2(k-1))))\geq 0,\ k\in\dN.\nonumber
\end{gather}
Straightforward calculations show that the solutions of the latter
system of inequalities are precisely the characters belonging to one
of the above sets $X_{ij}.$ One easily verifies that the sets
$X_{ij}$ are pairwise disjoint for different $(i,j).$
\hfill$\Box$\mn

\begin{prop}\label{prop_action_Z_su11}
$ $
\begin{enumerate}
    \item[$(i)$] $\chi_{00}^n$ is defined only for $n=0.$
    \item[$(ii)$] For $\chi_{st}\in X_{1k}\cup X_{2k},\ k\in\dZ,$ the $\chi_{st}^n$ is
    defined for all $n\in\dZ.$
    \item[$(iii)$] For $\chi_{st}\in X_{3k},\ k\in\dN_0,$ the $\chi_{st}^n$ is
    defined for $n\geq-k.$
    \item[$(iv)$] For $\chi_{st}\in X_{4k},\ k\in\dZ,$ the $\chi_{st}^n$ is
    defined for $n\leq k-1.$
\end{enumerate}
\end{prop}
\noindent\textbf{Proof.} Follows directly from Propositions
\ref{prop_p_action_sl_2} and \ref{prop_cBp_su_11}. \hfill $\Box$\mn

Set
$$
\Gamma:=X_{00}\cup X_{10}\cup X_{20}\cup X_{30}\cup
X_{4,-1}\subseteq\cBp.
$$

It follows from the previous propositions that each orbit under
the partial action of $\dZ$ on $\cBp$ intersects $\Gamma$ exactly
in one point, i.e. $\Gamma$ is a section of the partial action. As
in the case of $su(2),$ the topology on $\cBp$ is induced from the
standard topology on $\dR^2.$ Hence $\Gamma$ is a countably
separated Borel section of the partial action of $\dZ$ on $\cBp.$

Explicit formulas for the representations $\Ind\chi,\
\chi\in\Gamma,$ can be derived in a similar manner as in case of
$su(2).$ We omit the details. In the standard terminology of
representation theory of Lie algebras we have:
\begin{itemize}
  \item[-] the representation $\Ind\chi,\ \chi\in X_{00},$ is the trivial representation,
  \item[-] the representations $\Ind\chi,\ \chi\in X_{10},$ form the
  principal unitary series,
  \item[-] the representations $\Ind\chi,\ \chi\in X_{20},$ form the
  supplementary unitary series,
  \item[-] the representations $\Ind\chi,\ \chi\in X_{30}\cup X_{40},$ form the
  discrete unitary series.
\end{itemize}

Using this description we obtain the following

\begin{prop}
The representations $\Ind\chi,\ \chi\in\Gamma,$ are pairwise
non-equivalent and irreducible. Each irreducible well-behaved
representation of $\cA$ is unitarily equivalent to $\Ind\chi$ for
precisely one $\chi\in\Gamma.$ A $*$-representation of
$\cA=\cU(su(1,1))$ is well-behaved (in the sense of Definition
\ref{defn_well_beh}) if and only it is of the form $dU$ for some
unitary representation $U$ of the universal covering group of the
Lie group $SU(1,1).$
\end{prop}

We close this subsection with the following

\noindent\textbf{Remark.} For a character $\chi\in\cBp$ the
following three statements are equivalent:
\begin{enumerate}
\item[$(i)$] $\chi$ belongs to one of the series $X_{1k}$ or $X_{2k},\
k\in\dZ,$ corresponding to the principal or supplementary unitary
series,
\item[$(ii)$] $\chi^k$ is defined for all $k\in\dZ,$
\item[$(iii)$] $\chi(C)<0,$ where $C$ is the Casimir element defined above.
\end{enumerate}

\subsection{Enveloping algebra of the Virasoro algebra.}
Recall that the Virasoro algebra is the complex Lie algebra $Vir$
with generators $L_n,\ n \in \dZ$, and $C$ and defining relations
\begin{gather}\label{eq_vira}
[L_n,L_m]=(m-n)L_{n+m}+\delta_{n,-m}(n^3-n)/12{\cdot}C\ \mbox{and}\ [L_n,C]=0\ \mbox{for}\ n,m\in\dZ.%
\end{gather}

In this subsection we show that the unitary representations with
finite-dimensional weight spaces of the Virasoro algebra can be
identified with the well-behaved representations with respect to a
canonical grading of a quotient algebra of its enveloping algebra.
For results on unitary representations of $Vir$ we refer to
\cite{chari} and references therein.

Let $\cW$ denote the enveloping algebra of $Vir,$ that is, $\cW$ is
the unital $*$-algebra with generators $L_n,\ n\in\dZ,$ and $C$ and
the same defining relations (\ref{eq_vira}). It is a $*$-algebra
with involution determined by $L_n^*=L_{-n}$ for $n\in\dZ$ and
$C^*=C$. Lemma \ref{lemma_factorgrad} implies that $\cW$ is
$\dZ$-graded such that $L_n\in\cW_n$ and $C\in\cW_0.$

The main result in \cite{chari} states that there are precisely two families
of irreducible unitary representations of $\cW$ with
finite-dimensional weight spaces. The first series consists of
highest (resp. lowest) weight representations, i.e.
representations generated by a vector $v$ such that:\\
$(i)$ $L_0v=av$ for some $a\in\dC,$ $(ii)$ $L_nv=0$ for all $n>0$
(resp. $n<0$), $(iii)$
$Cv=zv$ for some $z \in \dC$.\\
These representations are uniquely defined by the pair
$(a,z)\in\dC^2.$ The possible values of $(a,z)$ for the
representation to be unitary (that is, a $*$-representation in our
terminology) are the following ones (see \cite{fqs}):
\begin{gather}\label{eq_values_a_z}
a\geq 0,\ z\geq 1,\ \mbox{or}\ z_n=1-\frac{6}{n(n+1)},\ a_n^{(p,q)}=\frac{(np+q)^2-1}{4n(n+1)},%
\end{gather}
where the integers $n,p,q$ satisfy $n\geq 2$ and $0\leq p<q<n.$

The other series of unitary representations are defined on spaces of
$\lambda$-densities (see \cite{chari}). They can be described as
follows. Let $\set{w_k}_{k\in\dZ}$ be an orthonormal base of
$l^2(\dZ).$ Then the action of $\cW$ on $l^2(\dZ)$ is given
 by
\begin{gather}\label{eq_action_Vir_lambda_dens}
L_kw_n=(n+a+k\lambda)w_{n+k},\ Cw_n=0,\ k,n\in\dZ,\
\lambda\in\frac{1}2+i\dR,\ a\in\dR.
\end{gather}

Let $\cI$ denote the two-sided $*$-ideal of $\cW$ generated by
elements
\begin{gather*}
bd-db,\ b,d\in\cW_0\ \mbox{and}\ a_k^*c_k^{}c_k^*a_k^{}-a_k^*a_k^{}c_k^*c_k^{},\ a_k,c_k\in\cW_k,\ k\in\dZ.%
\end{gather*}

\begin{lemma}\label{prop_common_ideal_Vir}
$\cI$ is contained in the intersection of all kernels of
representations described above.
\end{lemma}
\noindent\textbf{Proof.} We prove the assertion for
$*$-representations defined by (\ref{eq_action_Vir_lambda_dens}).
For highest and lowest weight representations the proof is similar.

We fix a $*$-representation $\pi$ given by
(\ref{eq_action_Vir_lambda_dens}), $k\in\dZ$ and $a_k,c_k\in\cW_k.$
It follows from (\ref{eq_action_Vir_lambda_dens}) that
$\pi(a_k)w_m=\mu_m w_{m+k},\ \pi(c_k)w_m=\nu_m w_{m+k},\ m\in\dZ,$
for some $\mu_m,\ \nu_m\in\dC.$ This implies that
$$
\pi(a_k^*c_k^{}c_k^*a_k^{})w_m=\overline{\lambda_m}\nu_m\overline{\nu_m}\lambda_m\cdot
w_m=\pi(a_k^*a_k^{}c_k^*c_k^{})w_m,
$$
for all $m\in\dZ.$ Taking $b,d\in\cW_0$ the same reasoning shows that
$\pi(bd)w_m=\pi(db)w_m,\ m\in\dZ.$ Therefore $\cI$ is contained in
$\ker\pi.$ \hfill$\Box$\mn

In view of Lemma \ref{prop_common_ideal_Vir} we introduce the
$*$-algebra $\cA=\cW/\cI.$ Let $\iota:\cW\to\cA$ be the quotient
mapping and put $l_k:=\iota(L_k)$ for $ k\in\dZ$ and $c=\iota(C).$ Since
the generators of $\cI$ are homogeneous, Lemma
\ref{lemma_factorgrad} implies that $\cA$ is again a $\dZ$-graded
$*$-algebra such that $l_k\in\cA_k,\ k\in\dZ,$ and $c\in\cA_0.$ As
usual we denote by $\cB$ the subalgebra $\cA_0.$

Because of the PBW-theorem there are two "natural" bases of the vector space $\cW:$
\begin{gather*}
\mathbf{B}_1=\set{C^kL_{n_1}L_{n_2}\dots L_{n_r}|n_1\leq
n_2\leq\dots\leq n_r,\ k,r\in\dN_0,n_i\in\dZ},\\
\mathbf{B}_2=\set{C^kL_{n_1}L_{n_2}\dots L_{n_r}|n_1\geq
n_2\geq\dots\geq n_r,\ k,r\in\dN_0,n_i\in\dZ}.
\end{gather*}
Fix $i{=}1,2$. Since all elements in $\mathbf{B}_i$ are homogeneous, the
elements $C^kL_{n_1}L_{n_2}\dots L_{n_r}\in\mathbf{B}_i$, $ \sum_j
n_j=0,$ form a vector space base of the algebra $\cW_0$. To
define a character of $\cW_0,$ it is therefore sufficient to define it on these
elements $C^kL_{n_1}L_{n_2}\dots L_{n_r}\in\mathbf{B}_i$.

Let $\pi$ be an irreducible unitary highest weight representation of
$Vir$ with weight vector $v.$ It defines a $*$-representation of
$\cW$ denoted also by $\pi.$ One easily checks that the subspace
$\dC\cdot v$ is invariant under all operators $\pi(b),\ b\in\cW_0.$
Therefore it defines a character $\chi$ on $\cW_0$ given by
$\chi(L_{n_1}\dots L_{n_k})=0,\ \chi(L_0)=a,\ \chi(C)=z$, where
$n_1\leq\dots\leq n_k,\ \sum_r n_r>0,$ and $(a,z)$ is one of the
pairs defined by (\ref{eq_values_a_z}). By Lemma
\ref{prop_common_ideal_Vir}, $\chi$ annihilates the ideal $\cI$, so
 it gives a character on the quotient algebra $\cB=\iota(\cW_0)$
which we denote again by $\chi.$ It is defined by
\begin{gather}\label{eq_high_wei_char}
\chi(l_{n_1}\dots l_{n_k})=0,\ \chi(l_0)=a,\ \chi(c)=z,\
\mbox{where}\ n_1\leq\dots\leq n_k\neq 0,\ \sum_r n_r=0,
\end{gather}
where $(a,z)$ is given by (\ref{eq_values_a_z}).
The character $\chi$ obviously belongs to $\cBp.$

From the lowest weight representations we get another series
of characters $\chi{\in}\cBp$ determined by
\begin{gather}\label{eq_low_wei_char}
\chi(l_{n_1}\dots l_{n_k})=0,\ \chi(l_0)=a,\ \chi(c)=z,\
\mbox{where}\ n_1\geq\dots\geq n_k\neq 0,\ \sum_r n_r=0,
\end{gather}
where $(a,z)$ is as in (\ref{eq_values_a_z}).

Let $\pi$ be a representation given by
(\ref{eq_action_Vir_lambda_dens}). Considering the restriction of $\pi$
to the subspace $\dC\cdot w_0$ we obtain a series of
characters $\chi\in\cBp$ defined by
\begin{gather}\label{eq_lamb_dens_char}
\chi(l_{n_1}\dots{}l_{n_k})=\prod_{r=1}^k(a-\sum_{s=1}^rn_s+n_r\lambda),\ \chi(c)=0,%
\end{gather}
where $a\in\dR,\ \lambda\in\frac{1}2+i\dR.$

Let $\Gamma\subseteq\cBp$ denotes the union of all characters
defined by the equations (\ref{eq_high_wei_char}),
(\ref{eq_low_wei_char}) and (\ref{eq_lamb_dens_char}).
\begin{prop}
Each orbit under the partial action of $\dZ$ on $\cBp$ contains
precisely one character from $\Gamma.$ The stabilizer of each
character in $\cBp$ is trivial. For every $\chi\in\cBp$,
$\iota\circ\Ind\chi$ is a $*$-representation of $\cW$ with
finite-dimensional weight spaces. Every irreducible
$*$-representation of $\cW$ with finite-dimensional weight spaces is
unitarily equivalent to $\iota\circ\Ind\chi$ for precisely one
$\chi\in\Gamma.$
\end{prop}
\noindent\textbf{Proof.} A straightforward computation shows that
$$
[l_0,l_{n_1}l_{n_2}\dots
l_{n_r}]=(n_1+n_2+\dots+n_r)l_{n_1}l_{n_2}\dots l_{n_r},\
n_i\in\dZ,\ r\geq 1.
$$
Since every $a_n\in\cA_n$ is a linear combination of the elements
$l_{n_1}l_{n_2}\dots l_{n_r},\ n_1+n_2+\dots+n_r=n,$ it follows that
\begin{gather}\label{eq_[L_0_a_n]}
[l_0,a_n]=na_n,\ \mbox{for all}\ a_n\in\cA_n,\ n\in\dZ.
\end{gather}

Let $\chi\in\cBp$ and $ n\in\dZ$. Assume that $\chi^n$ is defined.
Then there exists an $a_n\in\cA_n$ such that $\chi(a_n^*a_n^{})>0.$
Using (\ref{eq_[L_0_a_n]}) we get
\begin{gather}\label{eq_chi^n(L_0)}
\chi^n(l_0)=\frac{\chi(a_n^*l_0a_n)}{\chi(a_n^*a_n)}=\frac{\chi(a_n^*a_n^{}l_0+na_n^*a_n)}{\chi(a_n^*a_n)}=\chi(l_0)+n.
\end{gather}

Let $\pi:=\Ind\chi.$ Since $\chi$ satisfies condition
(\ref{eq_cond_compatible}), we can choose an orthonormal base of
vectors $e_k$ of the representation space $\cH_\pi$ such that
$\pi(l_0)e_k=\lambda_ke_k,$ where
$\lambda_k=\chi^k(l_0)=\chi(l_0)+k.$ This implies that $\pi(l_0)$
acts as a semisimple operator and that all eigenspaces of $\pi(l_0)$
are finite dimensional. It is also clear that the stabilizer of
$\chi$ is trivial, so the representation $\pi$ is irreducible by
Proposition \ref{prop_irr_iff_St_triv}. Therefore, by Theorem 0.5 in
\cite{chari} the representation $\iota\circ\pi$ is unitarily
equivalent either to a highest or lowest weight representation or to
a representation defined by (\ref{eq_action_Vir_lambda_dens}).

On the other hand, one easily verifies that $\Ind\chi$ gives rise
via $\iota$ either to a highest or lowest weight representation or
to a representation defined by (\ref{eq_action_Vir_lambda_dens}).
This implies that $\cBp$ is equal to the union of all orbits
$\Orb\chi,$ where $\chi\in\Gamma.$ \hfill $\Box$\mn

\section{Example: Representations of dynamical systems}\label{sect_dyn_sys}
Let $f\in\dR[x]$ be a fixed polynomial. In this section we consider
the $*$-algebra
$$\cA=\dC\langle a,a^*|aa^*=f(a^*a)\rangle.$$
Representations of the relation $aa^*=f(a^*a)$ for a measurable
real-valued function $f$ have been studied in detail in \cite{osam}
by other means. From the very beginning this important example gave
us intuition for developing our theory.\mn
%

By Lemma \ref{lemma_factorgrad} the $*$-algebra $\cA$ is
$\dZ$-graded with grading determined by $a\in\cA_1$ and
$a^*\in\cA_{-1}.$ From the definition of $\cA$ it follows that every
element of $\cA$ is a linear combination of elements
$$a^m,\ m\geq 0;\ a^{*k},\ k>0;\ a^{*k_1}a^{m_1}\dots a^{*k_r}a^{m_r},\ r\geq 1,\ k_1>0,\ m_r>0.$$
This implies that $\cA_n$ is the linear span of elements
$$a^{*k_1}a^{m_1}\dots a^{*k_r}a^{m_r},\ r\geq 1,\ k_1\geq 0,\ m_r\geq 0,\
\sum m_j-\sum k_i=n.$$ From the defining relation $aa^*=f(a^*a)$ we
easily derive that
\begin{gather}\label{apol}
ap(a^*a)=p(f(a^*a))a \ , \ p(a^*a)a^*=a^*p(f(a^*a))\ {\rm for}\ p \in \dC[t].
\end{gather}

\begin{lemma}
The $*$-algebra $\cB$ is commutative and spanned by the Hermitian
elements
\begin{gather}\label{eq_aux2}
a^{*k_1}a^{m_1}\dots a^{*k_r}a^{m_r},\ r\geq 1,\ k_1>0,\ m_r>0,\
\sum k_i=\sum m_j.
\end{gather}
\end{lemma}
\noindent\textbf{Proof.} For $k\in\dN$, let $\cB_{k}$ be the
subalgebra of $\cB$ generated by words $w$ in $a^*$ and $a$ satisfying equation (\ref{eq_aux2}) and of length $|w|$ less or
equal to $2k$.

We first prove by induction on $k$ that the algebra $\cB_{k}$ is generated
by words $w,\ |w|\leq 2k,$ of the form $a^*Q$ for some word $Q.$ For
$k=1$ the assertion holds, since $\cB_1$ is generated by the element
$a^*a.$ Suppose that the assertion is valid for $k>1.$ Let $w\in\cB,\
|w|\leq 2k+2, k>1.$ If $w=a^*Q$ for some word $Q,$ then the induction proof is
complete. Let $w=a^ra^*P,\ r>0,$ for some word $P.$ Using
( \ref{apol}) we get
$$w=a^ra^*P=a^{r-1}f(a^*a)P=a^{r-2}f(f(a^*a))aP=\dots=f^{r}(a^*a)a^{r-1}P.
$$
The word $a^{r-1}P$ belongs to the algebra $\cB_{k-1}$ and the
element $f^{r}(a^*a)$ belongs to $\cB_1.$ It follows that
$w\in\cB_{k-1}$ and the induction hypothesis applies. This completes our first induction proof.

A second similar induction proof shows that $\cB_k,\ k\geq 1,$ is
generated by words $w,\ |w|\leq 2k,$ of the form $a^*Qa$ for some
word $Q.$

We now prove by induction on $k$ that $\cB$ is commutative. The algebra
$\cB_1$ is generated by the single element $a^*a,$ so it is
commutative. Suppose that $\cB_{k},\ k\geq 1,$ is commutative. Let
$w_1$ and $w_2$ be words of length between $2k$ and $2k+2.$ Then, it
is enough to consider the case when the words $w_i$ have the form
$a^*P_ia,\ i=1,2,$ for some words $P_i.$ Remembering that $aa^* \in
\cB_1\subseteq \cB_k$ and using the induction hypothesis we compute
$$w_1w_2=a^*P_1aa^*P_2a=a^*aa^*P_1P_2a=a^*aa^*P_2P_1a=a^*P_2aa^*P_1a=w_2w_1.$$
Thus, $\cB_{k+1}$ is commutative.\hfill $\Box$\mn


\noindent\textbf{Remark.} The algebra $\cB$ is in general rather
"large" when the polynomial $f$ is not linear. We shall see this
from the description of the set $\cBp\subseteq\cBd$ given below.

The following Proposition allows us to use the theory developed in
the Section \ref{sect_grad_alg_com}.
\begin{prop}\label{prop_ds_is_compatible}
The $\dZ$-grading of the algebra $\cA$ introduced above satisfies
 condition (\ref{eq_cond_compatible}).
\end{prop}
\noindent\textbf{Proof.} Using a simple induction argument one can prove the
equalities
\begin{gather}\label{eq_ds_compatible}
\cA_n=\cB a^n,\ \cA_{-n}=a^{*n}\cB,\ n\in\dN.
\end{gather}
Then Proposition \ref{prop_full=>compatible} completes the proof.
\hfill$\Box$\mn

We now describe the set $\cBp,$ the partial action of $\dZ$ on it and
the representations associated with orbits of this partial action.

Let
$\chi\in\cBp$ be fixed and let $\pi$ be the induced representation
$\Ind\chi$. Let $h_k$ denote the vector
$[a^k\otimes 1]\in\cH_\pi$ for all $k\in\dZ.$ We always put
$a^{-k}:=a^{*k}$ for $k\in\dN$ and $a^0:=\mathbf{1}_\cA.$

If $h_k=0$ for some $k>0,$ then for any $c_k\in\cA_k$ we have
$[c_k\otimes 1]=0.$ Indeed, by (\ref{eq_ds_compatible}) there
exists $b\in\cB$ such that $c_k=ba^k$ which implies $[c_k\otimes
1]=[ba^k\otimes 1]=\pi(b)[a^k\otimes 1]=0.$ Moreover, for all
$m>0$ we have $h_{k+m}=\pi(a^m)h_k=0.$

Analogously, if $h_{-k}=0$ for some $k>0,$ then for any
$c_{-k}\in\cA_{-k}$ we have $[c_{-k}\otimes 1]=0.$ Indeed, by
(\ref{eq_ds_compatible}) there exists $b\in\cB$ such that
$c_{-k}=a^{*k}b.$ It implies $[c_{-k}\otimes 1]=[a^{*k}b\otimes
1]=[a^{*k}\otimes\chi(b)]=\chi(b)[a^{*k}\otimes 1]=0.$ For all
$m>0$ we have $h_{-k-m}=\pi(a^{*m})h_{-k}=0.$

Summarizing the above considerations we conclude that there exist
$K,M\in\dN\cup\set{\pm\infty},\ K<0<M$ such that $h_k\neq 0$ if
and only if $K<k<M.$ All $h_k$ are pairwise orthogonal and
Proposition \ref{prop_space_ind_rep} implies that the vectors $h_k$ span
$\cH_\pi.$ Using Proposition \ref{prop_space_ind_rep} we also
conclude that $\pi(a)h_k=\mu_kh_{k+1}$ for some $\mu_k\in\dC.$ We
choose numbers $\nu_k\in\dC\backslash\set{0},\ k\in\dZ,\ \nu_0=1,$ such
that the vectors $e_k:=\nu_k h_k,\ k\in\dZ$ are of the norm $1$ if
$h_k\neq 0$ and
\begin{gather}\label{eq_ds_rep_l_2_Z}
\pi(a)e_k=\lambda_ke_{k+1},\ \pi(a^*)e_k=\lambda_{k-1}e_{k-1}\
\mbox{for some}\ \lambda_k\geq 0,\ k\in\dZ.
\end{gather}
Thus the vectors $e_k,\ K<k<M,$ form an orthonormal base of
$\cH_\pi.$ Furthermore, $\lambda_k>0$ for $K<k<M-1$ and relation
(\ref{eq_ds_rep_l_2_Z}) together with the defining relation
$aa^*=f(a^*a)$ imply $\lambda_{k-1}^2=f(\lambda_k^2)$ for all
$K<k<M.$ In the case when $K$ resp. $M$ is finite we have also
$f(\lambda_{K+1}^2)=\lambda_K^2=0,$ resp. $\lambda_{M-1}=0,\
f(0)=\lambda_{M-2}^2.$

For the fixed character $\chi\in\cBp$ we consider the possible cases
depending on $K$ and $M.$\mn

\noindent\textbf{1.} Let $K<0$ and $M>0$ be finite, so that
$\lambda_{k-1}^2=f(\lambda_k^2)$ for $K<k<M,\ f(\lambda_{K+1}^2)=0,\
f(0)=\lambda_{M-2}^2.$ Since $\chi(c_k^*c_k^{})=\norm{[c_k\otimes
1]}^2=0$ for all $c_k\in\cA_k,\ k\leq K,\ k\geq M,$ the character
$\chi^k$ is defined only for $K<k<M.$ It implies that the stabilizer
of $\chi$ is trivial. Thus $\pi$ is an irreducible
finite-dimensional representation. Using (\ref{eq_ds_rep_l_2_Z}) we
get
\begin{gather*}
\pi(a)e_k=\lambda_ke_{k+1},\ \mbox{for}\ K<k<M-1,\ \pi(a)e_{M-1}=0,\\
\pi(a^*)e_k=\lambda_{k-1}e_{k-1}\ \mbox{for}\ K+1<k<M,\ \pi(a^*)e_{K+1}=0.%
\end{gather*}

\noindent\textbf{2.} Let only $M>0$ be finite, so that
$\lambda_{k-1}^2=f(\lambda_k^2)$ for all $k<M$ and
$f(0)=\lambda_{M-2}^2.$ As in the previous case we have that the
stabilizer of $\chi$ is trivial. Thus $\pi$ is an irreducible
infinite-dimensional representation. By (\ref{eq_ds_rep_l_2_Z})
we have
\begin{gather*}
\pi(a)e_k=\lambda_ke_{k+1},\ \mbox{for}\ k<M-1,\ \pi(a)e_{M-1}=0,\\
\pi(a^*)e_k=\lambda_{k-1}e_{k-1}\ \mbox{for}\ k<M.%
\end{gather*}
According to the terminology of \cite{osam}, $\pi$ is the \textit{Fock
representation.}

\noindent\textbf{3.} Let only $K<0$ be finite, so that
$\lambda_{k-1}^2=f(\lambda_k^2)$ for $K<k,\ f(\lambda_{K+1}^2)=0.$
As in the case \textbf{1.} the stabilizer of $\chi$ is trivial. Thus
$\pi$ is an irreducible infinite-dimensional representation. From
(\ref{eq_ds_rep_l_2_Z}) we obtain
\begin{gather*}
\pi(a)e_k=\lambda_ke_{k+1},\ \mbox{for}\ K<k,\\
\pi(a^*)e_k=\lambda_{k-1}e_{k-1}\ \mbox{for}\ K+1<k,\ \pi(a^*)e_{K+1}=0.%
\end{gather*}
In the terminology of \cite{osam}, $\pi$ is called \textit{anti-Fock
representation.}

\noindent\textbf{4.} Let both $K$ and $M$ be infinite, so that
$\lambda_{k-1}^2=f(\lambda_k^2)$ for $k\in\dZ.$ Recall that a sequence
$\set{\lambda_k}_{k\in\dZ}$ is called periodic if there exists
$m\in\dN,$ such that $\lambda_k=\lambda_{k+m}$ for all $k\in\dZ.$
The smallest such $m$ is called period of the sequence
$\set{\lambda_k}_{k\in\dZ}.$ We consider two subcases.

\noindent\textbf{4.1.} Let $\set{\lambda_k^2}_{k\in\dZ}$ be not
periodic. Then, in particular all numbers $\lambda_k,\ k\in\dZ,$ are pairwise
different. From (\ref{eq_ds_rep_l_2_Z}) we have
$\pi(a^*a)e_k=\lambda_k^2e_k$ and Proposition
\ref{prop_space_ind_rep} $(ii)$ implies that
$\chi^k(a^*a)=\lambda_k^2.$ Since $\set{\lambda_k^2}_{k\in\dZ}$ is
not periodic, all characters $\chi^k,\ k\in\dZ,$ are different. Thus,
the stabilizer of $\chi$ is trivial and representation $\pi$ defined
by (\ref{eq_ds_rep_l_2_Z}) is irreducible.

\noindent\textbf{4.2.} Let $\set{\lambda_k^2}_{k\in\dZ}$ be periodic
with a period $m\in\dN.$ Repeating the arguments from the previous
case it follows that the stabilizer $H$ of $\chi$ is equal to
$m\dZ\subset\dZ.$ Let $\cH_{\pi,m}$ be the Hilbert subspace spanned
by the vectors $e_{rm},\ r\in\dZ.$ Let $p\in\dN$ and
$c_{pm}\in\cA_{pm}.$ Then (\ref{eq_ds_compatible}) implies that
$c_{pm}=b_1a^{pm}$ for some $b_1\in\cB.$ Using
(\ref{eq_ds_rep_l_2_Z}) and Proposition \ref{prop_space_ind_rep}
$(ii)$ we get
\begin{gather*}
\pi(c_{pm})e_{rm}=\chi^{rm}(b_1)(\lambda_0\lambda_1\dots\lambda_{m-1})^pe_{(r+p)m}=
\chi(b_1)(\lambda_0\lambda_1\dots\lambda_{m-1})^pe_{(r+p)m}.
\end{gather*}
Thus $\pi(c_{pm})$ acts as a scalar multiple of the bilateral
shift on $\cH_{\pi,m}.$ This implies that
\begin{gather}\label{eq_chi_on_A_H}
\widetilde{\chi}(b_1a^{pm}):=\chi(b_1)(\lambda_0\lambda_1\dots\lambda_{m-1})^p,\
p\in\dN,
\end{gather}
defines a character on the algebra $\cA_H.$ The restriction of
$\widetilde{\chi}$ to $\cB$ coincides with $\chi.$ Therefore, by
Proposition \ref{prop_obstr} the Mackey obstruction of $\chi$ is
trivial. We denote by $\zeta_z,\ z\in\dT,$ the character of the
group $H=m\dZ$ defined by $\zeta_z(m)=z.$ Then, using
(\ref{eq_defn_coc_rep}) and (\ref{eq_chi_on_A_H}), we see that all
representations $\rho_z,\ z\in\dT,$ of $\cA_H$ satisfy condition
(\ref{eq_cond_rep_A_H}). These representations are
one-dimensional, that is, they are characters. For $c_{pm}=ba^{pm},\ p\in\dN,\
b\in\cB,$ we have
\begin{gather*}
\rho_z(c_{pm})=\chi(c_{pm}^*c_{pm}^{})^{1/2}\zeta_z(pm)=
\widetilde{\chi}(c_{pm}^*)^{1/2}\widetilde{\chi}(c_{pm}^{})^{1/2}z^p=
\chi(b^*b)^{1/2}(\lambda_0\lambda_1\dots\lambda_{m-1}z)^p,
\end{gather*}
where $z\in\dT.$

We now compute the representations induced from $\rho_z,\
z\in\dT.$ Let $\pi_z$ denotes the induced representation
$\Ind_{\cA_H\uparrow\cA}\rho_z$ on the space $\cH_z.$ One easily
verifies that the vectors
$$f_k=\chi(a^{*k}a^k)^{-1/2}[a^k\otimes{}1],\ k=0,\dots,m-1,$$
form an orthogonal base of the space $\cH_z.$ We calculate the
action of $\pi(a)$ on the base vectors $f_k.$ Using Proposition
\ref{prop_space_ind_rep} $(ii)$ and formulas
(\ref{eq_ds_rep_l_2_Z}) we find that
$\chi(a^{*k}a^k)=\lambda_0^2\lambda_1^2\dots\lambda_{k-1}^2,\
k\in\dN.$ Take $r=0,\dots,m-2.$ Then we have
$$
\pi_z(a)f_r=\frac{\chi(a^{(r+1)*}a^{r+1})^{1/2}}{\chi(a^{r*}a^r)^{1/2}}f_{r+1}=\lambda_rf_{r+1}.
$$
For $f_{m-1}$ we get
\begin{gather*}
\pi_z(a)f_{m-1}=\chi(a^{*(m-1)}a^{m-1})^{-1/2}[a^m\otimes
1]=\chi(a^{*(m-1)}a^{m-1})^{-1/2}[\mathbf{1}_\cA\otimes\rho_z(a^m)]=\\
=\chi(a^{*(m-1)}a^{m-1})^{-1/2}\widetilde{\chi}(a^m)[\mathbf{1}_\cA\otimes
1]=z\lambda_{m-1}f_0.
\end{gather*}\mn

Now suppose we are given a sequence $\lambda_k>0,\ K<k<M-1,$ where
$-\infty\leq K<0<M\leq\infty.$ Suppose also that
$f(\lambda_{K+1}^2)=0$ resp. $f(0)=\lambda_{M-2}^2$ in the case when
$K$ resp. $M$ is finite. We call such a sequence \textit{nonnegative
orbit of the dynamical system $(f,[0,+\infty)).$} Then
(\ref{eq_ds_rep_l_2_Z}) defines a $*$-representation $\pi$ of $\cA$
and the restriction of $\Res_\cB\pi$ to $\dC{\cdot} e_0$ gives a
character $\chi\in\cBp.$ Let us describe this characters $\chi$ in
the case \textbf{4.} explicitly. Take an element
$a^{*k_1}a^{m_1}\dots a^{*k_r}a^{m_r}\in\cB,\ r\geq 1,\ k_1>0,\
m_r>0,\ \sum k_i=\sum m_j.$ Using formulas (\ref{eq_ds_rep_l_2_Z})
we obtain
\begin{gather*}
\chi(a^{*k_1}a^{m_1}\dots
a^{*k_r}a^{m_r})=\prod_{i=0}^{m_r-1}\lambda_i\prod_{i=1}^{k_r}\lambda_{m_r-i}\dots\prod_{i=1}^{k_1}\lambda_{m_r-k_r+m_{r-1}-\dots+m_1-i}.
\end{gather*}

We summarize the above discussion in the following
\begin{prop}
The equations (\ref{eq_ds_rep_l_2_Z}) give a one-to-one
correspondence between nonnegative orbits of the dynamical system
$(f,[0,+\infty))$ and orbits of the partial action of $\dZ$ on
$\cBp.$ A representation $\pi$ defined by (\ref{eq_ds_rep_l_2_Z}) is
reducible if and only if the sequence $\lambda_k$ is periodic and
$\lambda_k>0$ for all $k\in\dZ.$
\end{prop}

Finally, we consider the problem of associating irreducible
well-behaved representations of $\cA$ with orbits in $\cBp$ (cf.
also \cite{osam}).
\begin{prop}
Assume that the function $f$ is one-to-one and there exists a
measurable set $\Gamma\subseteq[0,+\infty)$ containing precisely one
point from each nonnegative orbit of the dynamical system
$(f,[0,+\infty)).$ Then every irreducible well-behaved
representation of $\cA$ is associated with an orbit in $\cBp.$
\end{prop}
\noindent\textbf{Sketch of proof.} Let $\pi$ be an irreducible
well-behaved representation of $\cA.$ Then $\pi(a^*a)$ is
essentially self-adjoint. Using Proposition 33 in \cite{osam} we
conclude that the spectral measure of $\overline{\pi(a^*a)}$ is
ergodic with respect to $f.$ Applying Proposition 34 in
\cite{osam} it follows that the spectral measure of
$\overline{\pi(a^*a)}$ is concentrated on a single orbit of the
dynamical system $(f,[0,+\infty)).$ \hfill$\Box$\mn

For the case, when $f$ is not bijective, we refer to Theorem 15 in
\cite{osam}.

\section{Further examples}\label{sect_further_exam}

In this section we mention and briefly discuss some other classes of examples,
where the theory developed in the previous sections can be applied.

\begin{exam} \textit{(Compact quantum group algebras)}
The simplest example is the quantum group $SU_q(2)$, $q\in \dR$. The corresponding $*$-algebra $\cA$ has two generators $a$ and $c$ and defining relations
\begin{gather}\label{eq_SUq}
ac=qca,\ ca^*=qa^*c,\ c^*c=cc^*,\ aa^*+q^2cc^*=1,\ a^*a+c^*c=1.
\end{gather}
Then $\cA$ is $\dZ$-graded such that $a\in\cA_1,a^*\in\cA_{-1},\ c\in\cA_0.$

Set $N:=a^*a.$ Then  the subalgebra $\cB=\cA_0$ is equal to $\dC[c,c^*,N].$ It follows from (\ref{eq_SUq}) that $\cB$ is commutative and $\cA_k=a^k\cB,\ k\in\dZ.$ Proposition \ref{prop_full=>compatible} implies that condition (\ref{eq_cond_compatible}) is satisfied and our theory applies. From the defining relations (\ref{eq_SUq}) it follows at once that every $*$-representation is bounded and hence well-behaved by Proposition \ref{prop_well_behav_bounded}.

Suppose that $q\in(-1,1),\ q\neq 0.$ In what follows many arguments are similar to the case of the Weyl algebra (see Examples \ref{exam_weyl_alg}, \ref{exam_weyl_cone} and \ref{exam_weyl_act}). The last two equations in (\ref{eq_SUq}) imply $aa^*-q^2a^*a=1-q^2.$ By induction on $k\in\dZ$ one proves the following formulas:
\begin{gather}
\label{eq_SUq_tech} aa^{*k}=a^{*(k-1)}(q^{2k}(a^*a)-q^{2k}+1),\ a^*a^{k}=\frac{1}{q^{2k}}a^{k-1}(aa^*+q^{2k}-1),\\
\label{eq_SUq_tech2}a^na^{*n}=\prod_{k=1}^n(1-q^{2k}+q^{2k}N),\ a^{*n}a^n=\frac{1}{q^{2n}}\prod_{k=0}^{n-1}(N+q^{2k}-1).
\end{gather}

\noindent From Corollary \ref{corconeb} and formula (\ref{eq_SUq_tech2}) we obtain
\begin{gather}\label{eq_SUq_cone}
\sum\cA^2\cap\cB=\sum\cB^2+N\sum\cB^2+\dots+N(N+q^2-1)\dots(N+q^{2k}-1)\sum\cB^2+\dots
\end{gather}

\noindent Equations (\ref{eq_SUq}),(\ref{eq_SUq_cone}) imply that the only characters $\chi\in\cBd$ which are positive on $\sum\cA^2\cap\cB$ are:
\begin{itemize}
  \item $\chi_{k,u},k\in\dN_0,\ u\in\dC,\ |u|=1,$ defined by $\chi_{k,u}(N)=1-q^{2k},\ \chi_{k,u}(c)=q^{k}u,$ and
  \item $\chi_\infty$ defined by $\chi_\infty(N)=1,\ \chi_\infty(c)=0.$
\end{itemize}

\noindent From (\ref{eq_SUq_tech2}) we derive the partial action of $\dZ$ on $\cBp.$ For $\chi_{k,u}$,  $\alpha_n(\chi_{k,u})$ is defined and then equal to $\chi_{k-n,u}$
 if and only if $n\leq k.$ For  $\chi_\infty$ we have $\alpha_n(\chi_\infty)=\chi_\infty$  for all $n\in\dZ.$ The set $\set{\chi_{0,u},\ |u|=1}\cup\set{\chi_\infty}$ is a section of the action, i.e. it contains exactly one point from each orbit. By Proposition \ref{prop_irrorbit} every irreducible representation is associated to some orbit.

 The stabilizers of $\chi_{0,u},\ |u|=1$, are trivial. Hence, by Theorem \ref{thm_main}, $\pi_u:=\Ind\chi_{0,u}$ is the only irreducible representation, up to unitary equivalence, associated with $\Orb\chi_{0,u}.$ From Proposition \ref{prop_space_ind_rep} we obtain explicit formulas for the actions  on some orthobase $(f_k,k\in \dN_0)$, where $f_{-1}:=0$:
$$
\pi_u(N)f_k=(1-q^{2k})f_{k-1},\ \pi_u(a^*)f_{k}=(1-q^{2k+2})^{1/2}f_{k+1},\
\pi_u(c)f_k=q^k u f_k,\ k\in\dN_0.
$$

The stabilizer of $\chi_\infty$ is  $\dZ$ and $\cA_\dZ=\cA.$ Let $\rho$ be as in Theorem \ref{thm_main}, that is, $\rho$ is an irreducible representation of $\cA$ such that $\Res_\cB\rho$ is a multiple of $\chi_\infty.$ Then $\rho(c)=\rho(c^*)=0$ and $\rho(a^*a)=\rho(aa^*)=\textbf{1}.$ Hence $\rho$ is one-dimensional and equal to $\rho_u,$ where $\rho_u(c)=0,\ \rho_u(a)=u,\ u\in\dC,\ |u|=1.$ Since $\Ind\rho\simeq\rho,$ every irreducible representation associated with $\set{\chi_\infty}$ equals to some $\rho_u,\ |u|=1.$

\hfill $\edex$
\end{exam}

\begin{exam} \textit{(Quantum disk algebra.)} Suppose that $0\leq\mu\leq 1$, $ 0\leq q\leq
1$, and $(\mu,q)\neq(0,1)$. The two-parameter unit quantum disk
$*$-algebra $\cA$ has generators $a$ and $a^*$ and the defining
relation
$$qaa^*-a^*a=q-1+\mu(1-aa^*)(1-a^*a).
$$
Then $\cA$ is $\dZ$-graded such that $a\in\cA_1$ and
$a^*\in\cA_{-1}.$ As in the case of the dynamical systems in the
previous section one shows that $\cB=\cA_0$ is commutative and
condition (\ref{eq_cond_compatible}) is satisfied. There is a
one-to-one correspondence between orbits in $\cBp$ and orbits of the
dynamical system $(f,[0,+\infty))$ where
$$f(\lambda)=\frac{(q+\mu)\lambda+1-q-\mu}{\mu\lambda+1-\mu}.$$
For a more detailed analysis of this $*$-algebra see \cite{klimek}
and \cite{osam}, p.101.\hfill $\edex$
\end{exam}

\begin{exam} \textit{(Podles' quantum spheres.)} Let $q\in (0,\infty)$.
For $r\in[0,\infty)$, $\cO(\mathrm{S}^2_{qr})$ is the unital
$*$-algebra with generators $A\!=\!A^*, B,B^*$ and defining
relations (see \cite{podles} or \cite{ks}, 4.5)
\begin{gather*} AB\!=\!q^{-2} BA,\
AB^*\!=\!q^2 B^* A,\ B^* B\!=\!A-A^2+r,\ BB^* \!=\!q^2 A-q^4
A^2+r.
\end{gather*}
For $r=\infty$, the defining relations of
$\cO(\mathrm{S}^2_{q,\infty})$ are
\begin{gather*}\label{eq_rel_podles_sph}
AB=q^{-2}BA,\ AB^*=q^2B^*A,\ B^*B=-A^2+1,\ BB^*=-q^4 A^2+1.
\end{gather*}
In both cases $\cA=\cO(\mathrm{S}^2_{qr})$ is $\dZ$-graded such that
$B\in\cA_1$, $B^*\in\cA_{-1}$ and $A\in\cA_0.$ One can check that
$\cB=\cA_0$ is commutative and condition (\ref{eq_cond_compatible})
is fulfilled. It follows immediately from the defining relations
that all $*$-representations of $\cA$ are bounded.\hfill $\edex$
\end{exam}

\begin{exam} \textit{(Deformations of CAR algebra)} Let $q\in(0,1)$ be fixed.
The twisted canonical anti-commutation relations (briefly, TCAR)
$*$-algebra $\cA=\cA_q$ is generated by elements $a_i,a_i^*,\
i=1,\dots,d,$ with defining relations (see \cite{pusz})
\begin{gather*}
a_i^*a_i=1-a_i a_i^*-(1-q^2)\sum_{j<i}a_j a_j^*,\ i=1,\ldots,d,\ \\
a_i^*a_j=-q a_j a_i^*,\ a_j a_i =-q a_i a_j,\ i<j, a_i^2=0,\ i=1,\ldots,d.%
\end{gather*}

For $q=1$ we get the "usual" CAR algebra. For all $q\in(0,1]$, $\cA$ is
$(\dZ/2\dZ)^d$-graded such that $a_k^{},\ a_k^*\in\cA_{g_k}$, where
$g_1,\dots,g_d$ are generators of $(\dZ/2\dZ)^d$, the subalgebra
$\cB=\cA_0$ is commutative and condition (\ref{eq_cond_compatible})
is satisfied.

The Wick analogue of TCAR (denoted as WTCAR) was studied in
\cite{jsw,prosk,pst}. The WTCAR $*$-algebra $\cA$ is obtained from
TCAR by omitting the relations between $a_i$ and $a_j.$ Hence $\cA$
is $\dZ^d$-graded such that $a_k\in\cA_{g_k}$ where $g_1,\dots,g_d$
are generators of $\dZ^d.$ In this case the $*$-subalgebra
$\cB=\cA_0$ is not commutative. However, it was shown in
\cite{jsw,prosk} that in any irreducible representation of WTCAR the
relations
$$
a_ja_i=-q a_ia_j,\ i<j,\ a_i^2=0,\ i=1,\dots,d-1,
$$
hold automatically. Then our theory applies to the quotient of WTCAR
$*$-algebra by the latter relations.\hfill $\edex$
\end{exam}

\begin{exam} \textit{(Quantum algebras $\cU_q(su(2))$ and $\cU_q(su(1,1))$)}
For $q\in\dR$, $q^2\neq1$, the $q$-deformed enveloping algebra
$\cU_q(sl(2))$ is the complex unital (associative) algebra with
generators $E,F,K,K^{-1}$ and defining relations
\begin{gather*}
KK^{-1}=K^{-1}K=1,\ KEK^{-1}=q^2E,\ KFK^{-1}=q^{-2}F,\
[E,F]=\frac{K-K^{-1}}{q-q^{-1}}.
\end{gather*}
The involutions defining the $*$-algebras $\cU_q(su(2))$ and
$\cU_q(su(1,1))$ are given by the formulas
\begin{gather*}
E^*=F,\ F^*=E,\ K^*=K,\ K^{-1*}=K^{-1},\\
E^*=-F,\ F^*=-E,\ K^*=K,\ K^{-1*}=K^{-1},
\end{gather*}
respectively. Let $\cA$ be one of the $*$-algebras $\cU_q(su(2))$ or
$\cU_q(su(1,1)).$ Then $\cA$ is $\dZ$-graded with grading determined
by $E\in\cA_1,\ F\in\cA_{-1},$ and $ K,K^{-1}\in\cA_0,$ the
$*$-subalgebra $\cB=\cA_0$ is commutative, and condition
(\ref{eq_cond_compatible}) is valid. The Mackey analysis for $\cA$
is similar to that of $\cU(su(2))$ and $\cU(su(1,1))$.

The algebra $\cU_q(sl(2))$ was introduced in \cite{kulish}, see e.g.
\cite{ks}, 3.1. Representations of $\cU_q(su(2))$ and
$\cU_q(su(1,1))$ have been investigated in \cite{vaksman} and
\cite{burban}, respectively.\hfill $\edex$
\end{exam}

\begin{exam} \textit{(CAR algebras).} Let $\cA$ be the direct
limit of matrix $*$-algebras $M_{2^k}(\dC),\ k\in\dN,$ where the
embedding $M_{2^k}(\dC)\hookrightarrow M_{2^{k+1}}(\dC)$ is given by
the canonical injection $M_{2^k}(\dC)\otimes I_2\hookrightarrow
M_{2^{k+1}}(\dC)$. Here $I_2\in M_2(\dC)$ is the identity matrix.
The representation theory of $\cA$ was studied in \cite{garding},
see also \cite{sam_book} and \cite{kadison}.

Each matrix algebra $M_{n}(\dC)$ has a natural $\dZ$-grading such
that each matrix unit $e_{ij}$ belongs to the $(i{-}j)$-component.
Since the embeddings $M_{2^k}(\dC)\hookrightarrow M_{2^{k+1}}(\dC)$
respect this grading, $\cA$ is also $\dZ$-graded. One checks that
 condition (\ref{eq_cond_compatible}) is valid for
$M_{2^k}(\dC)$ which implies that the $\dZ$-grading on $\cA$ also
satisfies (\ref{eq_cond_compatible}). The $*$-subalgebra $\cB=\cA_0$
is the direct limit of commutative algebras $\dC^{2^k}$. It can be
considered as a (dense) $*$-subalgebra of the $*$-algebra of all
continuous functions on the Cantor set. The conditional expectation
defined by the $\dZ$-grading is strong, so $\cBp$ coincides with
$\cBd$ which is equal to the Cantor set. All representations of
$\cA$ are bounded. The partial action of $\dZ$ on $\cBp$ has trivial
stabilizers. All irreducible representations associated with orbits
in $\cBp$ are direct limits of representations. In this case the
assumptions of Proposition \ref{prop_irrorbit} are not satisfied and
there exist irreducible representations of $\cA$ arising from
ergodic measures under the partial action of $\dZ$ on $\cBp$ which
are not supported on single orbits.\hfill $\edex$
\end{exam}

\section*{Appendix}
The main result of this Appendix (Theorem \ref{thm_bhat}) is related
to condition $(i)$ of Definition \ref{defn_well_beh}, but it is also
of interest in itself. Its proof is based on the spectral theorem
for countable families of commuting self-adjoint operators, see
\cite{sam_book}, Theorem 1. We equip
$\dR^{\infty}=\dR\times\dR\times\dots$ with the product topology and
denote by $\mathbf{B}(\dR^{\infty})$ the Borel structure on
$\dR^\infty$ induced by this topology.
\begin{thm}\label{thm_spec_countable}
For each family $A_k,\ k\in\dN,$ of strongly commuting self-adjoint
operators there exists a unique resolution of the identity $E$ on the
Borel space $(\dR^\infty,\mathbf{B}(\dR^\infty))$ such that
\begin{gather*}
A_k=\int\lambda_k dE(\lambda_1,\lambda_2,\dots)\ \mbox{for all}\
k\in\dN.
\end{gather*}
\end{thm}

In the notation of Theorem \ref{thm_spec_countable}, the
\textit{joint spectrum} of the family $A_k,\ k\in\dN$, is the
intersection of all closed subsets $X$ of $\dR^\infty$ such that
$E(X)=E(\dR^\infty).$

Let $\cB$ be a commutative unital $*$-algebra. As in Section
\ref{sect_well_beh}, we equip the set $\cBd$ of all characters of
$\cB$ with the weakest topology for which all functions $f_b,\
b\in\cB,$ are continuous, where $f_b$ is defined by
$f_b(\chi)=\chi(b)$ for $\chi\in\cBd$. Clearly, if $\cB$ is
generated by elements $b_n,n\in\dN,$ then this topology coincides
with the weakest topology for which all functions $f_{b_n},\
n\in\dN,$ are continuous.

\begin{thm}\label{thm_bhat}
Suppose that $\cB$ is a countably generated commutative unital
$*$-algebra. We equip $\cBd$ with the Borel structure induced by the
weak topology. Let $\cC$ be a quadratic module of $\cB$ and let
$\cBp$ denote the set of all characters $\chi\in\cBd$ which are
nonnegative on $\cC.$ If $\pi$ is an integrable representation of
$\cB,$ then:
\begin{enumerate}
  \item[$(i)$] There exists a unique spectral measure $E_\pi$ on
  $\cBd$ such that
  \begin{gather*}
    \overline{\pi(b)}=\int f_b(\lambda)~d
    E_\pi(\lambda)\ \mbox{for all}\
   b \in \cB.
  \end{gather*}
  \item[$(ii)$] Assume in addition that $\langle\pi(c)\varphi,\varphi\rangle\geq
  0$ for all $c\in\cC$ and $\varphi\in\cD(\pi).$ Then the spectral measure
  $E_\pi$ is supported on $\cBp$ which
  is a closed subset of $\cBd.$
\end{enumerate}
\end{thm}
\noindent\textbf{Proof.} $(i)$ First we fix a sequence of
self-adjoint generators $b_k, k \in \dN,$ of the $*$-algebra $\cB$
and consider $\cBd$ as a subset of $\dR^\infty$ by identifying
$$
\cBd\ni\chi\longleftrightarrow(\chi(b_1),\chi(b_2),\chi(b_3),\dots)\in\dR^\infty.
$$

We prove that $\cBd$ is closed in $\dR^\infty,$ hence Borel. Let
$\chi_n=(\chi_n(b_1),\chi_n(b_2),\dots)\in\cBd,\ n\in\dN$ be a
sequence of characters converging to $\chi\in\dR^\infty$ in the
product topology. We claim that there is a character $\chi$ on $\cB$
such that $\chi(b_k):=\lim_{n\to\infty}\chi_n(b_k)$. Indeed, let
$m\in\dN$ and $p\in\dC[t_1,\dots,t_m]$ be a polynomial such that
$p(b_1,\dots,b_m)=0.$ Since
$$
p(\chi_n(b_1),\dots,\chi_n(b_m))=\chi_n(p(b_1,\dots,b_m))=0,
$$
we conclude that
$$
p(\chi(b_1),\dots,\chi(b_m))=p(\lim_{n\to\infty}\chi_n(b_1),\dots,\lim_{n\to\infty}\chi_n(b_m))=0
$$
for all $n\in\dN.$
Therefore $\chi\in\dR^\infty$ defines a character on $\cB,$ i.e.
$\chi\in\cBd.$

A sequence $\chi_n\in\cBd$ converges to $\chi\in\cBd$ if and only if
$\chi_n(b_k)=f_{b_k}(\chi_n)$ converges to $\chi(b_k)=f_{b_k}(\chi)$
for every fixed $k$ as $n\to\infty.$ Since the elements $b_k,k\in\dN,$ generate
$\cB,$ it follows that the topology on $\cBd$ induced from
$\dR^\infty$ coincides with the weak topology. In particular, the
Borel structure on $\cBd$ coincides with the one induced from
$\dR^\infty.$

Since $\pi$ is integrable, the operators $\overline{\pi(b_k)},\
k\in\dN,$ are self-adjoint and pairwise strongly commuting
(\cite{s1}, Corollary 9.1.14). Therefore, by Theorem
\ref{thm_spec_countable} there exist a spectral measure $E_\pi$ on
the set $\dR^\infty$ such that
\begin{gather*}
\overline{\pi(b_k)}=\int\lambda_k dE_\pi(\lambda_1,\lambda_2,\dots).
\end{gather*}
for all $ k\in\dN$. For every polynomial $p\in\dR[t_1,\dots,t_m]$ the operator
$p(\pi(b_1),\dots,\pi(b_m))$ is essentially self-adjoint and from basic
properties of spectral integrals we obtain
\begin{gather}\label{eq_spec_int_prop}
\overline{p(\pi(b_1),\dots,\pi(b_m))}=\int p(\lambda_1,\dots,\lambda_m)dE_\pi(\lambda_1,\lambda_2,\dots).%
\end{gather}

Next we show that the spectral measure $E_\pi$ is supported on $\cBd\subseteq\dR^\infty,$
or equivalently, that the joint spectrum
$\sigma(\overline{\pi(b_1)},\overline{\pi(b_2)},\dots)$ of the
family $\overline{\pi(b_k)},\ k\in\dN,$ is contained in $\cBd.$ Let
$x=(x_1,x_2,\dots)\in\dR^\infty$ be a point in
$\sigma(\overline{\pi(b_1)},\overline{\pi(b_2)},\dots).$ Again, let
$m\in\dN$ and $p_0\in\dR[t_1,\dots,t_m]$ be such that
$p_0(b_1,\dots,b_m)=0.$ Then we obtain
$$
\overline{\pi(p_0(b_1,\dots,b_m))}=0.
$$
Assume to the contrary that $p_0(x_1,x_2,\dots,x_m)\neq 0.$ Then for
every open neighborhood $O(x)$ we have $E_\pi(O(x))\neq 0.$ Using
(\ref{eq_spec_int_prop}) we get
$$
0=\overline{\pi(p_0(b_1,\dots,b_m))}=\overline{p_0(\pi(b_1),\dots,\pi(b_m))}=\int
p_0(\lambda_1,\dots,\lambda_m)dE_\pi(\lambda_1,\lambda_2,\dots)\neq
0,
$$
which is a contradiction. That is, we have
$p_0(x_1,x_2,\dots,x_m)=0.$ Thus we have shown that $\chi(b_k):=x_k$
defines a character and $E_\pi$ is supported on $\cBd.$ The
uniqueness of the spectral measure $E_\pi$ follows at once from the
corresponding assertion in Theorem \ref{thm_spec_countable}.

$(ii)$ Since $\cBd$ is a closed subset of the separable space
$\dR^\infty,$ $\cBd$ is also separable. Similar arguments as used in
the proof of $(i)$, show that $\cBp$ is closed in $\cBd.$

Assume to the contrary that $E_\pi(\cBd\backslash\cBp)\neq 0.$ Since
$\cBd$ is separable and $\cBp$ is a closed subset of $\cBd$, there
exists a countable dense subset $\set{\chi_i}_{i\in\dN}$ of
$\cBd\backslash\cBp.$ For every $\chi_i$ there exists an element
$c_i$ of $\cC$ such that $\chi_i(c_i)<0.$ Since
$\set{\chi_i}_{i\in\dN}$ is dense in $\cBd\backslash\cBp,$ the open
sets $f_{c_i}^{-1}((-\infty,0))$ cover $\cBd\backslash\cBp.$ From
the latter it follows that there exists a $k\in\dN$ such that
$E_\pi(f_{c_k}^{-1}((-\infty,0)))\neq 0.$ Hence there exists a
vector $\varphi\in\Ran E_\pi(f_{c_k}^{-1}((-\infty,0)))\cap\cD(\pi)$
such that $\langle\pi(c_k)\varphi,\varphi\rangle<0$ which
contradicts our assumption. \hfill $\Box$\mn

\begin{defn}
If $\cB,\ \pi$ and $E_\pi$ are as in the previous theorem, we shall say that
the integrable representation $\pi$ and the spectral measure $E_\pi$ are
\textit{associated} with each other.
\end{defn}

\bibliographystyle{amsalpha}

\end{document}